\documentclass[12pt]{amsart}
\usepackage{amssymb}
\usepackage{amsthm}
\usepackage{amscd}

\usepackage{amsmath}
\usepackage{epic,eepic}
\usepackage{graphicx}
\DeclareMathAlphabet{\mathpzc}{OT1}{pzc}{m}{it}

\usepackage{mathrsfs}
\usepackage{latexsym}

\usepackage{pifont}

\theoremstyle{plain}
\newtheorem{lemma}{Lemma}[subsection]
\newtheorem{prop}[lemma]{Proposition}
\newtheorem{thm}[lemma]{Theorem}
\newtheorem{cor}[lemma]{Corollary}
\newtheorem{aplemma}{Lemma~A.\hspace{-1.5mm}}
\newtheorem{approp}{Proposition~A.\hspace{-1.5mm}}
\newtheorem{apthm}{Theorem~A.\hspace{-1.5mm}}
\newtheorem{apcor}{Corollary~A.\hspace{-1.5mm}}
\newtheorem{intthm}{Theorem}

\usepackage[all]{xy}
\usepackage {cancel}

\newtheorem{conj}[lemma]{Conjecture}

\theoremstyle{definition}

\newtheorem{rema}[lemma]{Remark}

\newtheorem{remb}{Remark}

\newtheorem{defi}[lemma]{Definition}
\newtheorem{exa}[lemma]{Example}
\newtheorem{aprem}{Remark~A.\hspace{-1.5mm}}
\newtheorem{apdefi}{Definition~A.\hspace{-1.5mm}}
\newcommand{\bde}{\begin{defi}}
\newcommand{\ede}{\end{defi}}
\newcommand{\ble}{\begin{lemma}}
\newcommand{\ele}{\end{lemma}}
\newcommand{\bpr}{\begin{prop}}
\newcommand{\epr}{\end{prop}}
\newcommand{\bt}{\begin{thm}}
\newcommand{\et}{\end{thm}}
\newcommand{\bco}{\begin{cor}}
\newcommand{\eco}{\end{cor}}
\newcommand{\bre}{\begin{rema}}
\newcommand{\ere}{\end{rema}}
\newcommand{\brea}{\begin{rema}}
\newcommand{\erea}{\end{rema}\vspace{1mm}}
\newcommand{\breb}{\begin{remb}}
\newcommand{\ereb}{\end{remb}\vspace{1mm}}
\newcommand{\bex}{\begin{exa}}
\newcommand{\eex}{\end{exa}}
\newcommand{\bpf}{\begin{proof}}
\newcommand{\epf}{\end{proof}\vspace{1mm}}

\newcommand{\bade}{\begin{apdefi}}
\newcommand{\eade}{\end{apdefi}}
\newcommand{\bale}{\begin{aplemma}}
\newcommand{\eale}{\end{aplemma}}
\newcommand{\bapr}{\begin{approp}}
\newcommand{\eapr}{\end{approp}}
\newcommand{\bat}{\begin{apthm}}
\newcommand{\eat}{\end{apthm}}
\newcommand{\baco}{\begin{apcor}}
\newcommand{\eaco}{\end{apcor}}
\newcommand{\bare}{\begin{aprem}}
\newcommand{\eare}{\end{aprem}}

\newcommand{\be}{\begin{enumerate}}
\newcommand{\ee}{\end{enumerate}}
\newcommand{\bcd}{\[\begin{CD}}
\newcommand{\ecd}{\end{CD}\]}
\newcommand{\bit}{\begin{itemize}}
\newcommand{\eit}{\end{itemize}}
\newcommand{\bq}{\begin{quote}}
\newcommand{\eq}{\end{quote}}
\newcommand{\ba}{\begin{array}}
\newcommand{\ea}{\end{array}}

\newcommand{\mcD}{\mathcal{D}}
\newcommand{\mcE}{\mathcal{E}}
\newcommand{\mcF}{\mathcal{F}}

\newcommand{\mcI}{\mathcal{I}}

\newcommand{\mcL}{\mathcal{L}}
\newcommand{\mcM}{\mathcal{M}}
\newcommand{\mcN}{\mathcal{N}}
\newcommand{\mcO}{\mathcal{O}}
\newcommand{\mcP}{\mathcal{P}}

\newcommand{\mcT}{\mathcal{T}}

\newcommand{\mbB}{\mathbb{B}}

\newcommand{\mbF}{\mathbb{F}}
\newcommand{\mbG}{\mathbb{G}}

\newcommand{\mbP}{\mathbb{P}}

\newcommand{\mbZ}{\mathbb{Z}}

\newcommand{\mfM}{\mathfrak{M}}

\newcommand{\mfS}{\mathfrak{S}}

\newcommand{\mfd}{\mathfrak{d}}

\newcommand{\mfg}{\mathfrak{g}}

\newcommand{\mfl}{\mathfrak{l}}

\newcommand{\mfo}{\mathfrak{o}}
\newcommand{\mfp}{\mathfrak{p}}

\newcommand{\msE}{\mathscr{E}}
\newcommand{\msF}{\mathscr{F}}

\newcommand{\msN}{\mathscr{N}}

\newcommand{\msP}{\mathscr{P}}

\newcommand{\msR}{\mathscr{R}}

\newcommand{\msX}{\mathscr{X}}

\newcommand{\migi}{\rightarrow}
\newcommand{\longmigi}{\longrightarrow}

\newcommand{\isom}{\stackrel{\sim}{\migi}}

\newcommand{\migiincl}{\hookrightarrow}

\newcommand{\migisurj}{\twoheadrightarrow}

\newcommand{\mr}{\mathrm}
\newcommand{\hidden}[1]{\,}

\newcommand{\SSP}{\vspace{3mm}}
\newcommand{\LSP}{\vspace{5mm}}

\newcommand{\N}{N}
\newcommand{\M}{m}

\newcommand{\BB}{\blacktriangledown}

\newcommand{\EE}{E}

\newcommand{\Dual}{\rotatebox[origin=c]{180}{$D$}\hspace{-0.5mm}}

\begin{document}

\title[Differential modules   and dormant opers   of higher level]{Differential modules   and \\ dormant opers   of higher level}
\author{Yasuhiro Wakabayashi}
\date{}
\markboth{Yasuhiro Wakabayashi}{}
\maketitle
\footnotetext{Y. Wakabayashi: Department of Mathematics, Tokyo Institute of Technology, 2-12-1 Ookayama, Meguro-ku, Tokyo 152-8551, JAPAN;}
\footnotetext{e-mail: {\tt wkbysh@math.titech.ac.jp};}
\begin{abstract}
In the first half of the present paper, we study higher-level generalizations of differential modules in positive characteristic. These objects may be regarded   as ring-theoretic counterparts of vector bundles on a curve equipped with an action of the ring of (logarithmic) differential operators of finite level introduced by P. Berthelot (and C. Montagnon). The existence assertion for a cyclic vector of a differential module  is generalized to higher level under mild conditions. In the second half, we introduce  (dormant) opers of  level $\N > 0$ on a pointed smooth curve whose structure group is either $\mathrm{GL}_n$ or $\mathrm{PGL}_n$. Some of  the  results on higher-level differential modules are applied to  prove a duality theorem  between dormant $\mathrm{PGL}_n$-opers of level $N$ and dormant $\mathrm{PGL}_{p^N-n}$-opers of level $N$. Finally, in the case where the underlying curve is a $3$-pointed projective line, we establish a bijective correspondence between dormant $\mr{PGL}_2$-opers of  level $\N$ and certain tamely ramified coverings.

\end{abstract}
\tableofcontents 

\section*{Introduction}

\LSP
\subsection{} \label{SS002}
A {\it differential module} over a differential ring $(R, \partial : R \migi R)$ is an $R$-module $\EE$ equipped with an additive map $\nabla : \EE \migi \EE$ satisfying the Leibniz rule: $\nabla (a \cdot v) = \partial (a) \cdot v + a \cdot \nabla (v)$ ($a \in R$, $v \in \EE$). 
If $R$ is of characteristic $0$,   a differential module   may be regarded as a ring-theoretic counterpart  of
a $\mcD$-module (or equivalently, a sheaf equipped with a flat connection) on an algebraic  curve.
On the other hand, differential modules and $\mcD$-modules  in characteristic $p>0$ have many different features from those in characteristic $0$ (see, e.g.,  ~\cite{And}, ~\cite{Hon},  ~\cite{Kat3}, ~\cite{Kat4}).
For example,  unlike the case of characteristic $0$, Picard-Vessiot theory fails for  differential modules in characteristic $p$.
Despite this problem, such objects have attracted a lot of attention for various reasons, including in relation to the Grothendieck-Katz $p$-curvature conjecture.
We should note  that there are variations of the sheaf ``$\mcD$" defined on an algebraic  variety $X$ in characteristic $p$.
One is the ring of crystalline differential operators (following  the wording in  ~\cite{BMR}), which we denote by $\mcD^{(0)}_X$.
 Giving a   $\mcD^{(0)}_{X}$-module is equivalent to giving 
 an $\mcO_X$-module  together with  a flat connection.
This means that the notion of a $\mcD^{(0)}_{X}$-module corresponds exactly to  a differential module in the usual sense.
  Another variant is the ring of differential operators $\mcD^{(\infty)}_X$ in the sense of Grothendieck (cf. ~\cite[Section 16.8.1]{Gro}).
 A $\mcD^{(\infty)}_X$-module is often called a {\it stratified sheaf}
   and  interpreted as  
 an $\mcO_X$-module admitting  infinite Frobenius descent.
The ring-theoretic counterpart
 is known as   an {\it iterative differential module} (cf. ~\cite[\S\,1.2]{Oku}, ~\cite[Chap.\,13, \S\,13.3]{vdPS}).
 The notion of a stratified sheaf was introduced in D.  Gieseker's paper 
  (cf. ~\cite[Definition 1.1]{Gie}) and   discussed in, e.g., ~\cite{dSa}, ~\cite{Esn}, ~\cite{EsMe}, ~\cite{Kin1}, and ~\cite{Kin2}.
We also can find descriptions of iterative differential modules  in, e.g., ~\cite{Ern}, ~\cite{MavdP},  and ~\cite{Ros}. 

 Next, let us recall 
the ring of differential operators $\mcD_X^{(\M)}$ of  level $\M \in \mbZ_{\geq 0}$, as  introduced by P. Berthelot (cf. ~\cite{PBer1}, ~\cite{PBer2}); this kind of sheaf is an essential ingredient in defining  arithmetic $\mcD$-modules,  and 
it  may be positioned between $\mcD^{(0)}_X$ and $\mcD_{X}^{(\infty)}$.
In fact, 
the ring of crystalline differential operators  coincides with  Berthelot's $\mcD_X^{(0)}$ (i.e., $\mcD_X^{(\M)}$ for $\M =0$), and
there exists an inductive system 
\begin{align}
\mcD_X^{(0)} \migi \mcD_X^{(1)} \migi \mcD_X^{(2)} \migi \cdots  \migi  \mcD_X^{(\M)} \migi \cdots
\end{align}
 satisfying  $\varinjlim_{\M}\mcD_X^{(\M)} = \mcD^{(\infty)}_X$.
Moreover, C.  Montagnon generalized $\mcD_X^{(\M)}$ to the case where the underlying scheme is equipped with a log structure (cf. ~\cite[Chap.\,2, \S\,2.3, Definition 2.3.1]{Mon}).
This generalization allows us to deal with $\mcD^{(\M)}_X$-modules (in the logarithmic sense) for  (possibly singular) pointed curves $X$. 

In the first half of the present paper, 
we consider ring-theoretic counterparts of  (both non-logarithmic and logarithmic versions of)  $\mcD_X^{(\M)}$-modules, in other words, 
 higher-level generalizations  of differential modules.
The central character is 
an {\it $\M$-differential ring} (resp.,  an {\it $\M$-log differential ring}), which is defined  as  a ring $R$  in characteristic $p$ equipped with a collection of 
 certain additive endomorphisms $\partial_{\langle \bullet \rangle} := \{ \partial_{\langle j \rangle} \}_{j=0}^{p^\M}$ (cf. Definition \ref{Ee442}, (ii)).
Each such data $\msR := (R,  \partial_{\langle \bullet \rangle})$ yields a possibly noncommutative ring $\mcD_\msR^{(\M)}$ generated by the elements of $R$ and the set of  abstract symbols $\{ \partial_{\langle j \rangle} \}_{j \in \mbZ_{\geq 0}}$  
extending  the endomorphisms $\{ \partial_{\langle j \rangle} \}_{j =0}^{p^\M}$. 
(The definition of the related  ring can be found in ~\cite[Chap.\,1, \S\,1.1, Definition 1.1.1]{Kin1}.)
In particular, we obtain the notion of a $\mcD_\msR^{(\M)}$-module, which corresponds to the notion of a $\mcD_X^{(\M)}$-module.

An important ingredient  in the theory of  differential modules
is  the concept of a {\it cyclic vector}.
A cyclic vector  of a differential module $(\EE, \nabla)$  is an element $v \in \EE$ such that the elements $\nabla^0 (v) \left(= v \right), \nabla^1 (v),  \cdots, \nabla^l (v)$ (for some $l \geq  0$)  form a basis.
The choice of such an element enables  
$(\EE, \nabla)$ to be interpreted  as a higher-order linear differential operator on $R$;
accordingly, 
each element in $\EE$ that is horizontal with respect to $\nabla$ 
 can be described as a root function  of  that operator.
  A fundamental result is the one proving the existence of a cyclic vector in a general situation (cf.  ~\cite{ChKo}, ~\cite{Kat2}). 
We refer the reader to, e.g.,  ~\cite{Adj}, ~\cite{Ber}, ~\cite{Del}, and ~\cite{Kov}, for 
various  results concerning cyclic vectors.

Given an $\M$-differential ring  $\msR$ and $\mcD_\msR^{(\M)}$-module $(\EE, \nabla)$, 
we can describe the notion of   an {\it $\M$-cyclic vector} of  $(\EE, \nabla)$  (cf. Definition \ref{Ee440})  as a higher-level generalization of a cyclic vector in the classical sense.
An  $\M$-cyclic vector of $(\EE, \nabla)$ is, by definition,  an element of $\EE$ such that 
$\nabla_{\langle 0 \rangle} (v) \left(=v\right), \nabla_{\langle 1 \rangle} (v), \cdots,   \nabla_{\langle l \rangle} (v)$ (for some $l \geq 0$) forms a basis of $\EE$.
Also, by a {\it pinned $\mcD_{\msR}^{(\M)}$-module}, we mean a $\mcD_\msR^{(\M)}$-module together with an $\M$-cyclic vector.
Our study of $\M$-cyclic vectors is motivated by the observation that
a pinned $\mcD_{\msR}^{(0)}$-module
  is essentially the same as  a locally defined $\mr{GL}_n$-oper on a curve; 
   as such, 
   various properties of  $\M$-cyclic vectors can be used  to examine   higher-level generalizations of opers.
The first main result of  the present paper generalizes the classical assertion of the existence  of  a cyclic vector to higher level.

\SSP
\begin{intthm} [cf.  Theorem \ref{P023}] \label{E456}
Let $\M$ be a nonnegative integer and $(R, \partial_{\langle \bullet \rangle})$  an $\M$-differential field over $\mbF_p := \mbZ/p \mbZ$.
Suppose that there exists a positive integer $n$ satisfying 
 the following conditions:
\begin{itemize}
\item
The morphism $\mcD_{\msR, < n}^{(\M)} \migi \mr{End}_{\mbF_p} (R)$ naturally induced  by  $\partial_{\langle \bullet \rangle}$ is injective;
\item
The subfield $R^{\M+1} \left(:= \bigcap_{j=1}^{p^\M}  \mr{Ker} (\partial_{\langle j \rangle})\right)$ of $R$ contains at least $n$ nonzero elements.
\end{itemize}
(These conditions are fulfilled  if $n \leq p^{\M+1}$ and $R$ is either $k(t)$ or $ k(\!(t)\!)$ for a perfect  field $k$ over $\mbF_p$ equipped with the $\M$-derivation $\partial_{\langle \bullet \rangle} := \{ \partial_{\langle j \rangle} \}_{j}$ given by (\ref{reow98}).)
Then, 
each $\mcD_\msR^{(\M)}$-module  $(\EE, \nabla)$ with $\mr{rk}(\EE) = n$   admits an $\M$-cyclic vector.
 \end{intthm}
\SSP

After proving the above theorem, we discuss the {\it $p^{\M+1}$-curvature} of each  $\mcD_{\msR}^{(\M)}$-module in the situation (cf. \S\,\ref{SS99}) that $\msR := (R, \partial_{\langle \bullet \rangle})$ (resp., $\msR := (R, \breve{\partial}_{\langle \bullet \rangle})$))  is  a certain type of $\M$-differential ring 
  (resp., $\M$-log differential ring); 
    we use the notation  $\mcD_R^{(\M)}$ (resp., $\breve{\mcD}_R^{(\M)}$) to denote the ring $\mcD_\msR^{(\M)}$ for convenience.
The $p^{\M+1}$-curvature of a $\mcD_R^{(\M)}$-module  (resp., a $\breve{\mcD}_R^{(\M)}$-module) $(\EE, \nabla)$
is  
defined as an invariant measuring 
the extent to which the element
  $\partial_{\langle p^{\M+1} \rangle}$ (resp., $\breve{\partial}_{\langle p^{\M+1} \rangle}$) via  $\nabla$ vanishes.
We say that   $(\EE, \nabla)$  is {\it dormant} (cf. Definition \ref{D098}) if it has vanishing $p^{\M+1}$-curvature.
Here, suppose   that $(\EE, \nabla)$ is dormant  and the $R$-module $\EE$ is free and of rank $n >0$.
In the non-logarithmic case,
the structure of $(\EE, \nabla)$ is not difficult because  a classical result of Cartier implies that
 $(\EE, \nabla)$ is isomorphic to the direct sum of  finitely many copies of  the trivial $\mcD_R^{(\M)}$-module (cf. Corollary \ref{eoap09}).
On the other hand, in the logarithmic case,
 the formal completion of $(\EE, \nabla)$ is isomorphic to that of the  direct sum $\bigoplus_{j=1}^n (R, \nabla_{d_j})$ 
  for various elements $d_j$ ($j=1, \cdots, n$) of $\mbZ/p^{\M+1}\mbZ$,
 where each $\nabla_{d_j}$ denotes a $\breve{\mcD}_R^{(\M)}$-action on $R$ defined in (\ref{eoa78}).
The resulting multiset $[d_1, \cdots, d_n]$  is called  the {\it exponent} of $(\EE, \nabla)$ (cf. Definition \ref{D033}).
We will examine 
the relation of the notion of exponent with  the {\it residue} described  in \S\,\ref{SS9}  (cf. Propositions \ref{P0045}, (i), \ref{Ecsk}, and Remark \ref{eao78}), as well as its relation with the existence of  an $\M$-cyclic vector (cf. Proposition \ref{C001}).  
In addition, we will establish a duality between dormant pinned  $\breve{\mcD}_R^{(\M)}$-modules of rank $n$ (with $0 < n < p^{\M+1}$) and dormant pinned $\breve{\mcD}_R^{(\M)}$-modules of rank  $p^\N -n$ (cf. Corollary \ref{C003}, Proposition \ref{P0h23}).

\LSP
\subsection{} \label{SS003}

In the second half of the present paper, we study (dormant) $\mr{GL}_n$-opers and (dormant) $\mr{PGL}_n$-opers of level $\N > 0$. 
Here, let $\msX := (X, \{ \sigma_i \}_{i=1}^r)$ (cf. \S\,\ref{SS042d}), where $r \geq 0$,  be  an $r$-pointed smooth curve over an algebraically closed field $k$ of characteristic $p$.
A $\mr{GL}_n$-oper of level $\N$ on $\msX$  is, roughly speaking, a rank $n$ vector bundle on $X$ equipped with  both a $\mcD_X^{(\N-1)}$-action and complete flag structure satisfying a strict form of Griffiths transversality.
$\mr{GL}_n$-opers of level $\N = 0$ have been investigated from various points of view  (cf. ~\cite{BD1}, ~\cite{BD2}, ~\cite{BeBi}, ~\cite{Fr}, ~\cite{Wak8}).
In addition,  dormant $\mr{PGL}_2$-opers of level $\N$ (i.e., $\mr{PGL}_2$-opers of level $\N$  with vanishing $p^\N$-curvature)   on an unpointed smooth curve were discussed  in ~\cite{Hos2},  ~\cite{Wak6}, and ~\cite{Wak7} under the identification with  {\it $F^\N$-projective structures}.
In the case where the set of marked points $\{ \sigma_i \}_i$ of $\msX$ is nonempty,
we introduce the {\it radius} of a $\mr{PGL}_n$-oper of level $\N$ at each marked point $\sigma_i$ (cf. Definition \ref{epaddd}); this is an element of  $\mfS_n \backslash (\mbZ/p^\N \mbZ)^{\times n} / \Delta$ (where $\Delta$ denotes the image of the diagonal embedding $\mbZ/p^\N\mbZ \migiincl (\mbZ / p^\N \mbZ)^{\times n}$ and $\mfS_n$ denotes the symmetric group of $n$ letters acting on $(\mbZ / p^\N \mbZ)^{\times n}$ by permutation) induced from  the exponent of the $\breve{\mcD}_{k[\![t]\!]}$-module obtained by restricting that oper to the formal neighborhood of $\sigma_i$.
 Given an element  $\vec{\rho} \in (\rho_i)_{i=1}^r \in (\mfS_n \backslash (\mbZ/p^\N \mbZ)^{\times n} / \Delta)^{\times r}$, we set 
 \begin{align}\label{Kmdk}
{^\N}\mcO p_{n, \msX}^{^\mr{Zzz...}} \ \left(\text{resp.,} \ {^\N}\mcO p_{n, \msX, \vec{\rho}}^{^\mr{Zzz...}}  \right)
\end{align}
to be the set  of dormant $\mr{PGL}_{n}$-opers of level $\N$ (resp., dormant $\mr{PGL}_n$-opers of level $\N$ and radii $\vec{\rho}$) on  $\msX$.
Then,  
we verify   that ${^\N}\mcO p_{n, \msX}^{^\mr{Zzz...}} = \emptyset$ if $n > p^\N$ (cf. Corollary \ref{Peiri2})  and that ${^\N}\mcO p_{n, \msX}^{^\mr{Zzz...}}$ consists of exactly  one element  if $n = 1$.
Also, by the duality  of differential modules
  established  in the first part  of the present paper,  we obtain  the following assertions, generalizing ~\cite[Theorem A, (i) and (ii)]{Wak3},   ~\cite[\S\,4, Corollary 4.3.3]{Wak3}, and ~\cite[Theorem A]{Hos1}.

\SSP
\begin{intthm} [cf. Theorem \ref{T5949}, Corollary \ref{C90e8}] \label{efs89}
Suppose that $0 < n < p^\N$.
Then, the following assertions hold:
\begin{itemize}
\item[(i)]
There exists  a canonical bijection of sets
\begin{align}
\Dual_{n, \msX} : {^\N}\mcO p_{n, \msX}^{^\mr{Zzz...}} \isom  {^\N}\mcO p_{p^\N -n, \msX}^{^\mr{Zzz...}}
\end{align}
satisfying $\Dual_{p^\N -n, \msX} \circ \Dual_{n, \msX} = \mr{id}$.
In particular,  there exists exactly one isomorphism class of dormant $\mr{PGL}_{p^\N -1}$-oper of level $\N$ on $\msX$; i.e., 
 the following equality holds:
\begin{align}
\sharp ({^\N}\mcO p_{p^\N-1, \msX}^{^\mr{Zzz...}}) =1.
\end{align}

\item[(ii)]
Suppose further that $r > 0$, and let us take $\vec{\rho} := (\rho_i)_{i=1}^r  \in (\mfS_n \backslash (\mbZ/p^\N \mbZ)^{\times n} / \Delta)^{\times r}$.
Then, $\Dual_{n, \msX}$ is restricted to a bijection
\begin{align}
\Dual_{n, \msX, \vec{\rho}} :  {^\N}\mcO p_{n, \msX, \vec{\rho}}^{^\mr{Zzz...}} \isom  {^\N}\mcO p_{p^\N -n, \msX, \vec{\rho}^{\,\BB}}^{^\mr{Zzz...}},
\end{align}
where $\vec{\rho}^{\,\BB} := (\rho^\BB_i)_i$ is the set  defined in (\ref{eww299}),   satisfying $\Dual_{p^\N -n, \msX, \vec{\rho}^{\,\BB}} \circ \Dual_{n, \msX, \vec{\rho}} = \mr{id}$ (under the equality $\vec{\rho}^{\, \BB \BB} = \vec{\rho}$).
 \end{itemize}
 \end{intthm}
\SSP

The next topic concerns    the classification  of tamely ramified coverings of the projective line  in characteristic  $p$ with specified ramification data and fixed branch points.
For  related work on this problem, we refer the reader to  ~\cite{BoOs}, ~\cite{BoZa}, ~\cite{Ebe},
  ~\cite{Fab}, ~\cite{Oss1}, ~\cite{Oss2}, ~\cite{Oss3},  and ~\cite{O2}.
We know (cf. ~\cite{Mzk2}, ~\cite{Oss2}, ~\cite{O2}) that certain tamely ramified coverings with ramification indices $<p$ can be described in terms of dormant $\mr{PGL}_2$-opers (i.e., dormant torally indigenous bundles, in the sense of ~\cite{Mzk2}).
In particular, that description
allows us to translate  dormant $\mr{PGL}_2$-opers on a $3$-pointed projective line
into simple combinatorial data; this result is the starting point of the enumerative geometry of dormant opers studied in ~\cite{Wak8}.
In the present paper,  the situation  is generalized to the higher level case in order to deal with tamely ramified coverings having large ramification indices. 

Let  us consider
 the $3$-pointed projective line
  $\msP := (\mbP, \{ [0], [1], [\infty] \})$ over $k$, where $[x]$ (for each $x \in \{ 0, 1, \infty \}$) denotes the point of the projective line $\mbP$ determined by the value $x$.
We shall  write
\begin{align}
\mr{Cov}^{\mr{tame}}_{\msP}
\end{align}
for the set of equivalence classes of  finite, separable,  and tamely ramified coverings $\phi : \mbP \migi \mbP$ satisfying the following conditions:
\begin{itemize}
\item
The set of ramification points of $\phi$ coincides with  $\{[0], [1], [\infty] \}$;
\item
If $\lambda_x$ ($x = 0,1,\infty$) denotes the ramification index of $\phi$ at $[x]$, then  
$\lambda_0, \lambda_1, \lambda_\infty$ are all odd and satisfy the inequality $\lambda_0 + \lambda_1 + \lambda_\infty < 2p^\N$.
\end{itemize}
Here, for two such coverings $\phi_1, \phi_2 : \mbP \migi \mbP$,  we say that $\phi_1$ and $\phi_2$  are  equivalent  if $\phi_2 = h \circ \phi_1$ for some $h \in \mr{PGL}_2 (k) = \mr{Aut}_k (\mbP)$. 
The  final result of the present paper establishes, as described below,  a bijective correspondence between $\mr{Cov}_\msP^{\mr{tame}}$ and ${^\N}\mcO p_{2, \msP}^{^\mr{Zzz...}}$,  generalizing 
  ~\cite[Introduction, Theorem 1.3]{Mzk2}.
In future research, we will apply the resulting correspondence to develop the enumerative geometry of dormant opers of higher level.

\SSP
\begin{intthm} [cf. Theorem \ref{C00d3e} for the full statement] \label{Re345}
We can construct a canonical bijection of sets
\begin{align}
\Upsilon_{\msP} : \mr{Cov}^{\mr{tame}}_{\msP} \isom {^\N}\mcO p_{2, \msP}^{^\mr{Zzz...}}
\end{align}
satisfying the following condition: 
if $\phi$ is a tamely ramified covering classified by $\mr{Cov}^{\mr{tame}}_{\msP}$ whose ramification index at $[x]$ ($x=0,1,\infty$) is $\lambda_x$,
then the radii of the dormant $\mr{PGL}_2$-oper on $\msP$ determined by $\Upsilon_\msP (\phi)$  coincides with the image of $(\frac{1}{2} \cdot \overline{\lambda}_0, \frac{1}{2} \cdot \overline{\lambda}_1, \frac{1}{2} \cdot \overline{\lambda}_\infty)$ via the natural quotient $(\mbZ/p^\N \mbZ)^{\times 3} \migisurj  (\mfS_2 \backslash (\mbZ/p^\N \mbZ)^{\times 2}/\Delta)^{\times 3}$.
In particular, the set ${^\N}\mcO p_{2, \msP}^{^\mr{Zzz...}}$ is finite.
 \end{intthm}

\LSP
\subsection*{Acknowledgements} 
We are grateful for the many constructive conversations we had with {\it dormant opers of higher level}, who live in the world of mathematics!
Our work was partially supported by Grant-in-Aid for Scientific Research (KAKENHI No. 18K13385, 21K13770).

\vspace{10mm}
\section{Differential modules and cyclic vectors of higher level} \label{SS234} \SSP

First, we study differential modules of level $\M \geq 0$ and generalize the notion of a cyclic vector to  such modules, i.e., $\M$-cyclic vectors.
(For convenience, we will occasionally  include the case of $\M = -1$.)
At the end of this section,   we prove the  existence of  an $\M$-cyclic vector under   mild conditions (cf. Theorem \ref{P023}).

Throughout the present paper, we will fix a prime $p$.
Unless stated otherwise, 
all the rings appearing in the present paper are assumed to be  unital, associative, and commutative.

\LSP
\subsection{Modified binomial coefficients} \label{SgS5}

Let us fix an integer $\M \geq 0$.
For each  $l \in \mbZ_{\geq 0}$,  let $(q_l^{(\M)}, r_l^{(\M)})$ be the pair of nonnegative integers uniquely determined by the condition that $l = p^\M \cdot q^{(\M)}_l + r_l^{(\M)}$ and $0 \leq r^{(\M)}_l < p^\M$.
For each  pair of nonnegative integers $(j, j')$ with $j \geq j'$,  we set
\begin{align}
\begin{Bmatrix} j  \\ j' \end{Bmatrix}_{(\M)}
:= \frac{q^{(\M)}_{j}!}{q^{(\M)}_{j}! \cdot q^{(\M)}_{j-j'} !}, \hspace{10mm}
\left\langle\begin{matrix} j  \\ j' \end{matrix} \right\rangle_{(\M)} := \left(\begin{matrix} j \\ j' \end{matrix} \right)\cdot \left\{\begin{matrix} j  \\ j' \end{matrix} \right\}^{-1}_{(\M)}
\end{align}
(cf. ~\cite[\S\,1.1.2]{PBer1}).
Moreover, if $j''$ is an integer with $\mr{max}\{j, j -j' \} \leq j'' \leq j$, then we set
\begin{align}
\left\langle\begin{matrix} j  \\ j' \end{matrix} \right\rangle_{(\M)}^{[j'']} :=  \frac{j''!}{(j-j'')! \cdot  (j''-j)!  \cdot (j''+j'-j)!} \cdot \frac{q^{(\M)}_{j}! \cdot q^{(\M)}_{j-j'}!}{q^{(\M)}_{j''}!}.
\end{align}
In particular, we have $\left\langle\begin{matrix} j  \\ j' \end{matrix} \right\rangle_{(\M)}^{[j]} = \left\langle\begin{matrix} j  \\ j' \end{matrix} \right\rangle_{(\M)}$.
Note that all these values lie in $\mbZ_{\geq 0}$ and hence induce elements of $\mbF_p := \mbZ/p\mbZ$ via the natural quotient $\mbZ \migisurj \mbF_p$ even when  the integer  $j''$ is divisible by $p$.
When there is no fear of confusion, we will omit the notation  ``$(\M)$"; i.e., we will write
$q_l := q^{(\M)}_l$, $r_l := r_l^{(\M)}$, and 
\begin{align}
\begin{Bmatrix} j  \\ j' \end{Bmatrix} := \begin{Bmatrix} j  \\ j' \end{Bmatrix}_{(\M)},
\hspace{5mm}
\left\langle\begin{matrix} j \\ j' \end{matrix} \right\rangle := \left\langle\begin{matrix} j  \\ j' \end{matrix} \right\rangle_{(\M)},
\hspace{5mm}
\left\langle\begin{matrix} j  \\ j' \end{matrix} \right\rangle^{[j'']} := \left\langle\begin{matrix} j  \\ j' \end{matrix} \right\rangle_{(\M)}^{ [j'']}. 
\end{align}

\LSP
\subsection{Differential rings of higher level} \label{SS5}

Let us fix an integer $\M \geq -1$ and a ring $R_0$ over $\mbF_p$.
In the following discussion,  the non-resp'd portion deals with the non-logarithmic case and the resp'd portion deals with the logarithmic case.

\bde \label{Ee442}
\begin{itemize}
\item[(i)]
Let $R$ be a ring over $R_0$.
An {\bf $\M$-derivation} (resp., {\bf $\M$-log derivation}) on $R$ relative to $R_0$ is a collection 
\begin{align}
\partial_{\langle \bullet \rangle} := \{ \partial_{\langle j \rangle} \}_{0 \leq j \leq p^\M}
\end{align}
consisting of $R_0$-linear endomorphisms $\partial_{\langle j \rangle}$ of $R$, regarded  as an $R_0$-module, 
satisfying the following conditions:
\begin{itemize}
 \item[(a)]
If $j = 0$, then  $\partial_{\langle j \rangle} =\mr{id}_{R}$;
 \item[(b)]
 If $0 < j  \leq p^\M$, then the following equalities hold:
 \begin{itemize}
 \item
 $\displaystyle{\partial_{\langle j \rangle} ( a \cdot b) =\sum_{j' + j'' = j} \begin{Bmatrix} j \\ j'\end{Bmatrix} \cdot  \partial_{\langle j' \rangle}(a) \cdot \partial_{\langle j'' \rangle}} (b)$  for any elements  $a, b\in R$;
 \item
  $\displaystyle{\partial_{\langle j' \rangle} \circ \partial_{\langle j-j' \rangle} = 
 \left\langle\begin{matrix} j  \\ j' \end{matrix} \right\rangle  \cdot \partial_{\langle j\rangle}} \left(\text{resp.,} \ \partial_{\langle j' \rangle} \circ \partial_{\langle j-j'\rangle} = \sum_{j'' = \mr{max}\{j', j-j' \}}^{j} \left\langle\begin{matrix} j  \\ j' \end{matrix} \right\rangle^{[j'']} \cdot \partial_{\langle j''\rangle} \right)$ for
  any integer $j'$ with  $0 \leq j' \leq  j$. 
 \end{itemize}
 \end{itemize}

\item[(ii)]
By an {\bf $\M$-differential ring} (resp., {\bf $\M$-log differential ring}) over $R_0$, we mean the pair 
\begin{align}
\msR := (R, \partial_{\langle \bullet \rangle})
\end{align}
 consisting of a ring  $R$ over $R_0$ and an $\M$-derivation (resp., $\M$-log derivation)  $\partial_{\langle \bullet \rangle}$ on $R$.
(In particular, a $(-1)$-differential ring is  the same as an $R$-module.)
Finally, by  an {\bf $\M$-differential field}, we mean an $\M$-differential ring $(R, \partial_{\langle \bullet \rangle})$ such that $R$ is a  field.
\end{itemize}
\ede

\SSP
\begin{rema} \label{R034}
\begin{itemize}
\item[(i)]
If we are given a $0$-derivation $\{ \partial_{\langle j \rangle} \}_{0 \leq j \leq 1}$  on $R$,
then $\partial_{\langle 1 \rangle}$ defines a derivation on $R$ (over $R_0$) in the usual sense.
Under the correspondence $(R, \{ \partial_{\langle j \rangle} \}_{0 \leq j \leq 1}) \leftrightarrow (R, \partial_{\langle 1 \rangle})$,
  the notion of a $0$-differential ring  coincides with the usual notion   a differential ring.

\item[(ii)]
Let $R$ be a ring over  $R_0$.
The second equality in condition (b) above implies that
an $\M$-derivation  (resp., an $\M$-log derivation)  $\{ \partial_{\langle j \rangle} \}_{0 \leq j \leq p^\M}$ on $R$
 is uniquely determined by its subset 
$\{ \partial_{\langle p^j\rangle} \}_{0 \leq j \leq \M}$.
\end{itemize}
\end{rema}
\SSP

  Let us fix an $\M$-differential ring  (resp., $\M$-log differential ring) $\msR := (R, \partial_{\langle \bullet \rangle})$  over $R_0$, where  $\partial_{\langle \bullet \rangle} := \{ \partial_{\langle j \rangle} \}_{0 \leq j \leq p^\M}$.
If $\M \geq 0$, then $\partial_{\langle \bullet \rangle}$ can be extended uniquely, by adding  $\partial_{\langle j \rangle}$'s for various large $j$'s,  to a collection of $R_0$-linear endomorphisms $\{ \partial_{\langle j\rangle} \}_{j \in \mbZ_{\geq 0}}$  of $R$  indexed by   $\mbZ_{\geq 0}$ so that
the two equalities described in  condition (b) above  are fulfilled   for all $j >0$.
We will refer to  such a collection $\{ \partial_{\langle j\rangle} \}_{j \in \mbZ_{\geq 0}}$ as an {\bf extended $\M$-derivation} (resp., {\bf extended $\M$-log derivation}) on $R$.
 By using the resulting  collection, 
we obtain  the possibly noncommutative ring 
\begin{align}
\mcD^{(\M)}_{\msR}
\end{align}
over $R_0$ generated by the collection of symbols $\{ \partial_{\langle j \rangle}\}_{j \in \mbZ_{\geq 0}}$  subject to the following relations:
 \begin{itemize}
 \item
 $\partial_{\langle 0 \rangle} =1$;
 \item
 $\displaystyle{\partial_{\langle j \rangle} \cdot a =\sum_{j' + j'' = j} \begin{Bmatrix} j \\ j'\end{Bmatrix} \cdot  \partial_{\langle j' \rangle}(a) \cdot \partial_{\langle j'' \rangle}}$ for any $j \in \mbZ_{\geq 0}$ and $a\in R$;
 \item
   $\displaystyle{\partial_{\langle j' \rangle} \cdot \partial_{\langle j - j' \rangle} = 
 \left\langle\begin{matrix} j \\ j' \end{matrix} \right\rangle  \cdot \partial_{\langle j\rangle}  \left(\text{resp.,} \ \partial_{\langle j' \rangle} \cdot \partial_{\langle j-j'\rangle} = \sum_{j'' = \mr{max}\{j', j-j' \}}^{j} \left\langle\begin{matrix} j  \\ j' \end{matrix} \right\rangle^{[j'']} \cdot \partial_{\langle j'' \rangle}\right)}$ for any 
 integers $j, j'$ with  $0 \leq j' \leq j$.
 \end{itemize}
We shall set  $\mcD_\msR^{(-1)} := R$.
The ring $\mcD^{(\M)}_{\msR}$
  admits two  $R$-module structures given by left and right multiplications.
For each $l \in \mbZ_{\geq 0}$, 
we shall denote by $\mcD_{\msR, <l}^{(\M)}$ the two-sided $R$-submodule of $\mcD_{\msR}^{(\M)}$ generated by the products $\partial^{a_1}_{\langle j_1 \rangle} \cdots  \partial^{a_s}_{\langle j_s \rangle}$ ($s \geq 1$) with $j_i \leq p^{\M}$ ($i=1, \cdots, s$) and  $\sum_{i=1}^s a_i j_i < l$.
The collection $\{ \mcD_{\msR, < l}^{(\M)} \}_l$
forms an increasing filtration on $\mcD_\msR^{(\M)}$ with  $\bigcup_{l} \mcD_{\msR, < l}^{(\M)} = \mcD_\msR^{(\M)}$.
Also, one may verify that the $R_0$-algebra   $\mcD_\msR^{(\M)}$ is generated by the elements of $R$ and the set  $\{ \partial_{\langle p^j\rangle} \}_{0 \leq j \leq \M}$, and it forms a left and right noetherian ring  if $R$ is noetherian (cf. ~\cite[\S\,1, Proposition 1.2.4, (i)]{PBer1}, ~\cite[\S\,2, Proposition 2.3.2, (b)]{Mon}).

\LSP
\subsection{Differential modules of higher level} \label{SSf5}

Let us take an $R$-module $\EE$.
By a {\bf  (left) $\mcD_{\msR}^{(\M)}$-module structure} on $\EE$, we mean 
a  left $\mcD_{\msR}^{(\M)}$-action (i.e., an $R_0$-algebra homomorphism)  $\nabla : \mcD_{\msR}^{(\M)} \migi \mr{End}_{R_0}(\EE)$   on $\EE$ extending its  $R$-module structure.
An $R$-module equipped with a $\mcD_\msR^{(\M)}$-module structure is called a {\bf (left) $\mcD_\msR^{(\M)}$-module}, or a {\bf differential module over $\msR$}.
Moreover, we can define, in a natural manner,  the notion of an isomorphism between $\mcD_\msR^{(\M)}$-modules.
Given a $\mcD_{\msR}^{(\M)}$-module structure  $\nabla$ on $\EE$ and an integer $j \in \mbZ_{\geq 0}$,
we shall write $\nabla_{\langle j \rangle} := \nabla (\partial_{\langle j \rangle})$.
If $\msR$ is non-logarithmic (resp., logarithmic), then   the assignment $\nabla \mapsto \{ \nabla_{\langle j \rangle} \}_{j \in \mbZ_{\geq 0}}$ determines a bijective correspondence between the set of  $\mcD_{\msR}^{(\M)}$-module structures on $\EE$ and the set of collections of $R_0$-linear endomorphisms $\{ \nabla_{\langle j \rangle}\}_{j \in \mbZ_{\geq 0}}$  of $\EE$ 
  satisfying  the following conditions:
   \begin{itemize}
\item
$\nabla_{\langle 0 \rangle} = \mr{id}_\EE$;
\item
$\displaystyle{\nabla_{\langle j \rangle} (a \cdot v) = \sum_{j' + j'' = j} \begin{Bmatrix} j\\ j'\end{Bmatrix}  \cdot \partial_{\langle j' \rangle} (a) \cdot \nabla_{\langle j'' \rangle} (v)}$ for any integer $j >0$
  and any elements  $a \in R$, $v \in \EE$;
\item
$\displaystyle{\nabla_{\langle j' \rangle} \circ \nabla_{\langle j-j' \rangle} = \left\langle \begin{matrix} j \\ j' \end{matrix}  \right\rangle \cdot   \nabla_{\langle j' \rangle} \left(\text{resp.,} \ \nabla_{\langle j' \rangle} \circ \nabla_{\langle j-j' \rangle} = \sum_{j'' = \mr{max}\{j', j-j' \}}^{j} \left\langle\begin{matrix} j  \\ j' \end{matrix} \right\rangle^{[j'']} \cdot \nabla_{\langle j'' \rangle}\right)}$ for any integers $j, j'$ with $0 \leq j' \leq j$.
\end{itemize}
Because of this correspondence, we will not distinguish these two additional structures on $\EE$.

\SSP
\begin{rema} \label{R093}
Let us consider the case of $\M = 0$.
Suppose that $(\EE, \nabla)$ is  a $\mcD_{\msR}^{(0)}$-module.
 Then,  the $R$-module $\EE$ together with  the endomorphism $\nabla_{\langle 1 \rangle} \left(:= \nabla (\partial_{\langle 1 \rangle}) \right)$ specifies a differential module, in the classical sense,  over the differential ring corresponding to $\msR$ (cf. Remark \ref{R034}, (i)). 
\end{rema}

\SSP

\begin{rema} \label{R090}
Let  $\{ \partial_{\langle j \rangle}\}_{j \in \mbZ_{\geq 0}}$ be an extended $\M$-derivation (resp., extended $\M$-log derivation) on $R$.
To make the integer $\M$ explicit, 
 we here write $\partial_{\langle j \rangle}^{(\M)} := \partial_{\langle j \rangle}$ ($j \in \mbZ_{\geq 0}$).
For  an integer $\M'$ with $0 \leq \M' \leq \M$,
the endomorphism $\partial_{\langle j \rangle}^{(\M')} :=  \frac{q_j^{(\M')}!}{q_{j}^{(\M)}!}\cdot \partial_{\langle j \rangle}^{(\M)}$ of $R$  is well-defined, and
 the collection   $\{ \partial_{\langle j \rangle}^{(\M')} \}_{j \in \mbZ_{\geq 0}}$
 forms an extended $\M'$-derivation (resp., extended $\M$-log derivation) on $R$.
Let us set  $\mcD_{\msR}^{(\M')} := \mcD_{(R, \{ \partial_{\langle j \rangle}^{(\M')}\}_j)}^{(\M')}$.
Then,  
 the assignment $\partial_{\langle j \rangle}^{(\M')} \mapsto \frac{q_j^{(\M')}!}{q_{j}^{(\M)}!}\cdot \partial_{\langle j \rangle}^{(\M)}$ ($j \in \mbZ_{\geq 0}$) determines an $R_0$-algebra homomorphism  
$\mcD_{\msR}^{(\M')} \migi \mcD_\msR^{(\M)}$.
  This homomorphism allows us to construct  a $\mcD_{\msR}^{(\M')}$-module by means of 
   each $\mcD_\msR^{(\M)}$-module.
\end{rema}
\SSP

We shall denote by 
\begin{align}
\mfM \mfo \mfd (\mcD^{(\M)}_\msR)
\end{align}
the category of  $\mcD^{(\M)}_\msR$-modules.
(In particular, $\mfM \mfo \mfd (\mcD_\msR^{(-1)})$ coincides with the category of $R$-modules.)
This category has the structure of a tensor product: 
 given two  $\mcD_\msR^{(\M)}$-modules $(\EE', \nabla')$ and $(\EE'', \nabla'')$, we set 
 \begin{align}
 (\EE', \nabla') \otimes (\EE'', \nabla'') := (\EE' \otimes_R \EE'', \nabla' \otimes \nabla''),
 \end{align}
  where 
$\nabla' \otimes \nabla''$
denotes the $\mcD_\msR^{(\M)}$-module structure on the tensor product $\EE' \otimes_R \EE''$ determined by
\begin{align}
(\nabla' \otimes \nabla'')_{\langle j \rangle} (v' \otimes v'') := \sum_{j'+ j''= j} \begin{Bmatrix} j \\ j' \end{Bmatrix} \cdot \nabla'_{\langle j' \rangle} (v') \otimes \nabla''_{\langle j'' \rangle}(v'')
\end{align}
for any $j \in \mbZ_{\geq 0}$, $v' \in \EE'$, and $v'' \in \EE''$.
Similarly,  we can construct a $\mcD_{\msR}^{(\M)}$-module structure  on $\mr{Hom}_{R}(\EE', \EE'')$ arising  from $\nabla'$ and $\nabla''$.
In particular, for a $\mcD_\msR^{(\M)}$-module $(E, \nabla)$,
we can define the dual $(\EE^\vee, \nabla^\vee)$ of $(\EE, \nabla)$.
In this way, $\mfM \mfo \mfd (\mcD^{(\M)}_\msR)$ is equipped with a structure of  closed monoidal  category.

\LSP
\subsection{Varying levels} \label{SS56}

Fix an integer $\M \geq 0$ and 
an  integer  $l$ with $0  \leq l \leq \M +1$.
Furthermore, we set 
 \begin{align}
 R^{l} := \bigcap_{j=1}^{p^{l-1}}\mr{Ker}(\partial_{\langle j \rangle}) \left(= \bigcap_{j=0}^{l-1}\mr{Ker}(\partial_{\langle p^{j} \rangle})\right), 
 \end{align}
which is an  $R_0$-subalgebra of $R$.
Here, we obtain a sequence of inclusions between $R_0$-algebras
\begin{align}
R^{\M+1} \subseteq R^{\M} \subseteq \cdots \subseteq R^1 \subseteq R^0 = R.
\end{align}
Note that, if $\msR$ is an $\M$-differential field, then  $R^l$  (for every $l$) forms a subfield of $R$.
For a nonnegative integer $j$ with $j + l \leq \M+1$ and an $R^{j+l}$-module $\EE$, we shall set
$F^{(l)*}_j (E)$ to be the $R^j$-module defined as 
\begin{align} \label{E308}
F^{(l)*}_j (\EE) := R^j \otimes_{R^{j+l}} \EE.
\end{align}
For simplicity, 
we write $F^{(l)*} (\EE) := F^{(l)*}_0 (\EE)$.

Also, for each $j \in \mbZ_{\geq 0}$,
the endomorphism $\partial_{\langle j p^l \rangle}$ is restricted to an $R_0$-linear endomorphism of $R^l$; we will abuse the notation by writing  $\partial_{\langle j p^l \rangle}$ for this restriction.
Then,  the collection 
\begin{align} \label{Eq32}
\msR^l := (R^l, \{ \partial_{\langle j p^l \rangle} \}_{0 \leq j \leq p^{\M - l}})
\end{align}
 forms an $(\M-l)$-differential ring over $R_0$.
In particular, we obtain the $R_0$-algebra  $\mcD_{\msR^l}^{(\M-l)}$  and
 a sequence of inclusions
\begin{align} \label{Eq33}
\left( R^{\M+1}=\right)\mcD_{\msR^{\M+1}}^{(-1)} \migiincl \mcD_{\msR^\M}^{(0)} \migiincl  \mcD_{\msR^{\M-1}}^{(1)} \migiincl  \cdots \migiincl  \mcD_{\msR^1}^{(\M-1)} \migiincl \mcD_{\msR}^{(\M)}.
\end{align}

Next, 
given   a $\mcD_\msR^{(\M)}$-module $(\EE, \nabla)$, 
 we shall  write
  \begin{align} \label{Eq34}
 \EE^{l} := \bigcap_{j=1}^{p^{l-1}}\mr{Ker}(\nabla_{\langle j \rangle}) \left(= \bigcap_{j=0}^{l-1}\mr{Ker}(\nabla_{\langle p^{j} \rangle})\right),
  \end{align}
  where $\EE^{0} := E$.
In particular, we obtain a sequence of inclusions between modules
\begin{align} \label{Eq35}
E^{\M+1} \subseteq E^{\M} \subseteq \cdots \subseteq E^1 \subseteq E^0 = E.
\end{align}
Note that $\EE^l$ forms an $R^l$-module via the natural inclusion $R^l \migiincl R$,  and $\nabla_{\langle j p^l\rangle}$ (for each $j$) is restricted to an $R^l$-linear endomorphism of $\EE^l$; we will abuse the notation by writing  $\nabla_{\langle j p^l \rangle}$ for  this restriction.
One may verify that  the collection 
\begin{align}
\nabla^l :=  \{ \nabla_{\langle j p^l \rangle} \}_{j}
\end{align}
 forms a  $\mcD_{\msR^l}^{(\M-l)}$-module structure on $\EE^l$.
Also, if $f : (\EE, \nabla) \migi (\EE', \nabla')$ is a morphism of $\mcD_\msR^{(\M)}$-modules, then it is restricted  to a morphism of $\mcD_{\msR^l}^{(\M-l)}$-modules  $f^l : (\EE^l, \nabla^l) \migi (\EE'^l, \nabla'^l)$.
The resulting assignments $(\EE, \nabla) \mapsto (\EE^l, \nabla^l)$ and $f \mapsto f^l$ define a functor
\begin{align} \label{E323}
\Xi^{\downarrow (l)} : \mfM \mfo \mfd (\mcD^{(\M)}_{\msR}) \migi \mfM \mfo \mfd (\mcD_{\msR^l}^{(\M-l)}).
\end{align}
 
 Conversely,    
 given a $\mcD_{\msR^l}^{(\M-l)}$-module $(\EE, \nabla)$,
we can construct a $\mcD_{\msR}^{(\M)}$-module structure $F^{(l)*}(\nabla)$
  on $F^{(l)*}(\EE)$
   given by 
\begin{align} \label{E310}
F^{(l)*}(\nabla)_{\langle j \rangle} (a \otimes v) := \sum_{j' + j'' = j}\begin{Bmatrix} j \\ j'\end{Bmatrix} \cdot \partial_{\langle j' \rangle} (a) \otimes \nabla_{\langle j''/p^j \rangle} (v) 
\end{align}
for any $a \in R$ and $v \in \EE$, 
where $\nabla_{\langle s \rangle} := 0$ if $s \notin \mbZ_{\geq 0}$ (cf. ~\cite[Chap.\,2, \S\,2.2, Proposition 2.2.4, (ii)]{PBer2}, ~\cite[Chap.\,3, \S\,3.3, Proposition 3.4.1, (ii)]{Mon}).
In particular,
for an $R^{\M+1}$-module  $E'$, 
we obtain 
\begin{align} \label{dsskw}
\nabla^\mr{can}_{\msR, \EE'} := F^{(\M)*}(\nabla) : \mcD_\msR^{(\M+1)} \migi \mr{End}_{R_0}(F^{(\M+1)*}(\EE')).
\end{align}
Each morphism of $\mcD_{\msR^l}^{(\M-l)}$-modules $f: (\EE, \nabla) \migi (\EE', \nabla')$  induces  a morphism of $\mcD_{\msR}^{(\M)}$-modules $F^{(l)*}(f) : (F^{(l)*}(\EE), F^{(l)*}(\nabla)) \migi (F^{(l)*}(\EE'), F^{(l)*}(\nabla'))$.
The  resulting assignments $(\EE, \nabla) \mapsto (F^{(l)*}(\EE), F^{(l)*}(\nabla))$ and $f \mapsto F^{(l)*}(f)$ define a functor 
\begin{align} \label{E309}
\Xi^{\uparrow (l)} :  \mfM  \mfo \mfd (\mcD_{\msR^l}^{(\M-l)}) \migi \mfM \mfo \mfd (\mcD_\msR^{(\M)}).
\end{align}
This functor is compatible with the formation of the tensor product and  left adjoint to  $\Xi^{\downarrow (l)}$.

Here, let us describe the unit and counit morphisms for the adjunction  ``$\Xi^{\uparrow (l)}  \dashv \Xi^{\downarrow (l)}$".
If $(\EE, \nabla)$ is  a $\mcD_{\msR}^{(\M)}$-module, then
the  natural   inclusion $\EE^l \migiincl \EE$, which is $\mcD_{\msR^l}^{(\M-l)}$-linear, 
 extends  to
 a morphism of 
 $\mcD_\msR^{(\M)}$-modules
\begin{align} \label{Err4}
\tau_{(\EE, \nabla)}^{\downarrow \uparrow (l)} :
 \left((\Xi^{\uparrow (l)} \circ \Xi^{\downarrow (l)}) ((E, \nabla)) =\right) (F^{(l)*}(\EE^l), F^{(l)*}(\nabla^l))
 \migi (\EE, \nabla).
\end{align}
On the other hand, if $(\EE, \nabla)$ is a $\mcD^{(\M-l)}_{\msR^l}$-module,
then the morphism  $\EE \migi F^{(l)*}_{}(\EE)$ given by $v \mapsto 1 \otimes v$ is restricted to a morphism of  $\mcD_{\msR^l}^{(\M-l)}$-modules
\begin{align} \label{ewo08}
\tau^{\uparrow \downarrow (l)}_{(\EE, \nabla)} : (\EE, \nabla) \migi (F^{(l)*}(\EE)^l,  F^{(l)*}(\nabla)^l) \left(= (\Xi^{\downarrow (l)} \circ \Xi^{\uparrow (l)}) ((E, \nabla)) \right).
\end{align}
The formation of $\tau_{(\EE, \nabla)}^{\downarrow \uparrow (l)}$ (resp., $\tau^{\uparrow \downarrow (l)}_{(\EE, \nabla)}$) is functorial with respect to $(\EE, \nabla)$.

\LSP
\subsection{$\M$-cyclic vectors} \label{SS45}

Let $\M$ be a nonnegative integer.
In this subsection, we introduce the notion of an $\M$-cyclic vector, as a higher-level generalization  of a cyclic vector.

\SSP
\bde \label{Ee440}
Let $(\EE, \nabla)$ be a $\mcD_\msR^{(\M)}$-module. 
An element $v$ of $\EE$ is called   an {\bf $\M$-cyclic vector} of $(\EE, \nabla)$ if 
there  exists a positive  integer $n$ such that
the collection 
\begin{align}
\nabla_{\langle 0 \rangle}(v) \left(= v \right),   \nabla_{\langle 1 \rangle}(v), \cdots, \nabla_{\langle n-1 \rangle}(v)
\end{align}
  forms a basis of the $R$-module $\EE$.
 (In particular, a $\mcD_\msR^{(\M)}$-module admitting an $\M$-cyclic vector must be a  free $R$-module of finite rank.)
\ede
\SSP

\begin{rema} \label{R89989}
Suppose that $\M = 0$ and $\msR$ is non-logarithmic.
Then, 
each $\mcD_\msR^{(0)}$-module structure $\nabla$ on an $R$-module satisfies
$\nabla_{\langle j \rangle} = \nabla_{\langle 1 \rangle}^j$ for every $j$. 
It follows that a $0$-cyclic vector is the same as a cyclic vector in the classical sense.
\end{rema}

\SSP

\bde \label{Ee44dd0}
\begin{itemize}
\item[(i)]
A {\bf pinned $\mcD_{\msR}^{(\M)}$-module} is a triple 
\begin{align}
(\EE, \nabla, v)
\end{align}
 consisting of a $\mcD_{\msR}^{(\M)}$-module $(\EE, \nabla)$
and an $\M$-cyclic vector $v$ of it.
For a pinned $\mcD_{\msR}^{(\M)}$-module $(\EE, \nabla, v)$,  the {\bf rank} of  $(\EE, \nabla, v)$ is defined as the rank of the free $R$-module $\EE$.
\item[(ii)]
Let $(\EE, \nabla, v)$ and $(\EE', \nabla', v')$ be pinned $\mcD_{\msR}^{(\M)}$-modules.
An {\bf isomorphism of pinned $\mcD_{\msR}^{(\M)}$-modules} from $(\EE, \nabla, v)$ to $(\EE', \nabla', v')$ is a morphism of $\mcD_\msR^{(\M)}$-modules $f : (\EE, \nabla) \migi (\EE', \nabla')$ with $f (v) = v'$.
\end{itemize}
\ede
\SSP

Let us describe several basic properties on pinned $\mcD_\msR^{(\M)}$-modules:

\SSP
\bpr \label{edaP}
\begin{itemize}
\item[(i)]
Any morphism of pinned $\mcD_{\msR}^{(\M)}$-modules is surjective.
\item[(ii)]
Suppose that we are given a pinned $\mcD_\msR^{(\M)}$-module  $(\EE, \nabla, v)$  and a $\mcD_\msR^{(\M)}$-module  $(\EE', \nabla')$.
Denote by $\mr{Hom}((\EE, \nabla), (\EE', \nabla'))$ the set of  morphisms of $\mcD^{(\M)}_\msR$-modules from $(\EE, \nabla)$ to $(\EE', \nabla')$.
Then, the map of sets
\begin{align}
\mr{Hom}((\EE, \nabla), (\EE', \nabla')) \migi E'
\end{align}
given by $f \mapsto f(v)$ is injective.
\item[(iii)]
Let $(E, \nabla)$ be a $\mcD_{\msR}^{(\M)}$-module, 
   $v$  an element of $\EE$, and $\M'$ an integer with $0 \leq \M' \leq \M$.
Denote by $\nabla^{(\M')}$ the $\mcD_{\msR}^{(\M')}$-module structure on $E$ induced from $\nabla$ via the $R_0$-algebra homomorphism $\mcD_\msR^{(\M')} \migi \mcD_\msR^{(\M)}$ (cf. Remark \ref{R090}).
Then,  $v$ forms an $\M$-cyclic vector of $(\EE, \nabla)$ if $v$  forms an  $\M'$-cyclic vector of  the $\mcD_\msR^{(\M')}$-module $(\EE, \nabla^{(\M')})$.
\end{itemize}
\epr
\begin{proof}
To prove assertion (i), 
let  us take  a morphism of pinned  $\mcD_\msR^{(\M)}$-modules  $f : (\EE, \nabla, v) \migi (\EE', \nabla', v')$.
Then, since
 $f (\nabla_{\langle j \rangle} (v)) = \nabla'_{\langle j \rangle} (f (v)) = \nabla'_{\langle j \rangle} (v')$ ($j =0,1,2, \cdots$),
 the assertion follows from the fact that the set of   elements $\{ \nabla'_{\langle j \rangle} (v')\}_{j \in \mbZ_{\geq 0}}$ generates $\EE'$.

The remaining assertions, i.e., (ii) and (iii),  can be verified from the definition of an $\M$-cyclic vector (we will omit the details).
\end{proof}
\SSP

Now, let us prove  the following theorem, asserting the existence of an $\M$-cyclic vector in a general situation.
The proof  is based on the one in ~\cite[\S\,3, Theorem 3.11]{ChKo}.

\SSP
\bt[cf. Theorem \ref{E456}]\label{P023}
Let $\msR := (R, \partial_{\langle \bullet \rangle})$ be an $\M$-differential field over $\mbF_p$.
Suppose that there exists a positive integer $n$ 
satisfying  the following conditions:
\begin{itemize}
\item
The morphism $\mcD_{\msR, < n}^{(\M)} \migi \mr{End}_{\mbF_p}(R)$ naturally induced   by  $\partial_{\langle \bullet \rangle}$  is injective.
(This means that, for each nonzero element $D \in \mcD_{\msR, < n}^{(\M)}$, there exists an element $a$ of $R$ with $D (a) \neq 0$.)
\item
The subfield $R^{\M +1}$ of $R$ contains at least $n$ nonzero elements.
\end{itemize}
Then,
each  $\mcD_\msR^{(\M)}$-module  $(\EE, \nabla)$ with $\mr{rk}(E) = n$ admits
  an $\M$-cyclic vector.
\et
\begin{proof}
It suffices to  consider the case of $n >1$.
Suppose that a nonzero element $v$ of $\EE$ is not an $\M$-cyclic vector.
Then,  there exists an integer $l$ with $1 \leq l < n$ such that
\begin{align} \label{E3342}
v \wedge  \nabla_{\langle 1 \rangle} (v) \wedge \cdots \wedge \nabla_{\langle l-1\rangle} (v) \neq 0 \hspace{5mm} \text{and} \hspace{5mm}  \nabla_{\langle l \rangle} (v)= \sum_{j=0}^{l-1}a_j  \cdot  \nabla_{\langle j \rangle} (v)
\end{align}
for some $a_0, \cdots, a_{l-1} \in R$.
For simplicity, we set
$v_i = \nabla_{\langle j \rangle }(v)$ ($i=0,1,  \cdots, l-1$).
Choose an element  $u$ of $E$ not in the span of 
$\{ v, \nabla_{\langle 1 \rangle}(v), \cdots, \nabla_{\langle l-1 \rangle}(v) \}$.
We  extend the $R$-lineary independent set $\{ v_0, v_1, \cdots, v_{l-1} \}$ to a basis of $\EE$ by first adjoining $u$, and then, if necessary, some elements $e_1, \cdots, e_{n-l-1}$ of $E$.
For each  integer $j$ with  $0 \leq j \leq l$, we shall write
\begin{align}
\nabla_{\langle j \rangle}(u) = \sum_{i=0}^{l-1} \alpha_{ji} \cdot v_{i} + \beta_{j} \cdot u + \sum_{i=1}^{n-l-1} \gamma_{ji} \cdot e_{i},
\end{align}
where $\alpha_{ji}, \beta_j, \gamma_{ji} \in R$.
In particular, $\alpha_{0i} = \gamma_{0i} = 0$ and $\beta_0 = 1$.

For each integer $r$ with $0 \leq r \leq l$,  write  $L_r$  for the  $\mbF_p$-linear endomorphism of $R$ given  by 
\begin{align}
L_r :=\sum_{i=0}^r \begin{Bmatrix} r\\ i\end{Bmatrix}\cdot \beta_i \cdot \partial_{\langle r-i \rangle} 
\end{align}
(hence  $L_0  = \mr{id}_R$).
Next,
let us define $L$ to be the $R_0$-linear endomorphism of $R$ given by 
\begin{align}
L  := L_l - \sum_{r=0}^{l-1}a_r \cdot L_r  = \partial_{\langle l \rangle} + c_{l-1}\cdot  \partial_{\langle l-1 \rangle} + \cdots + c_1 \cdot \partial_{\langle 1 \rangle}+c_0,
\end{align}
where 
\begin{align}
c_i :=  \begin{Bmatrix} l\\ l-i\end{Bmatrix}  \cdot  \beta_{l-i}-\sum_{r=i}^{l-1}  \begin{Bmatrix} r\\ r- i\end{Bmatrix} \cdot a_{r} \cdot  \beta_{r-i}
\end{align}
($i=0, \cdots, l-1$).
Since the operator $L$ defines an element of $\mcD_{\msR, <n}^{(\M)}$,
the injectivity assumption  of  the morphism  $\mcD_{\msR, < n}^{(\M)} \migi \mr{End}_{\mbF_p}(R)$
implies that there exists  an element $z$ of $R$ with $L (z) \neq 0$.

 Let us choose  an indeterminate $\lambda$ over $R$.
 Extend the $\M$-derivation $\partial_{\langle \bullet \rangle}$ on $R$ to  the rational function field $R (\lambda)$ by defining $\partial_{\langle j \rangle}(\lambda) = 0$ ($j = 1, 2, \cdots$).
 To be precise, this $\M$-derivation can be obtained  by first defining $\partial_{\langle \bullet \rangle}$ on $R [\lambda] = R \otimes_{R^{\M+1}} R^{\M+1} [\lambda]$ by $\partial_{\langle l \rangle}(a \otimes b) = \partial_{\langle l \rangle}(a) \otimes b$ and then extending via the quotient rule.
 The tensor product $\EE^\lambda := R (\lambda) \otimes_R \EE$ has the natural $\mcD_{(R (\lambda), \partial_{\langle \bullet \rangle})}^{(\M)}$-module structure  obtained by defining 
 \begin{align}
 \nabla_{\langle j \rangle} (a \otimes w) = \sum_{j' + j'' = j}\begin{Bmatrix} j \\ j'\end{Bmatrix} \cdot \partial_{\langle j' \rangle}(a) \otimes \nabla_{\langle j'' \rangle} (w)
 \end{align} 
($j =0, 1,2, \cdots$) for any $a \in R(\lambda)$ and  $w \in \EE$.
Here, 
 we set 
 \begin{align}
 \hat{v}:= v + \lambda \cdot z \cdot  u \left(= 1 \otimes v + z \cdot  \lambda \otimes u \right) \in \EE^\lambda.
 \end{align}
 For each integer $r$ with  $0 \leq r \leq l$, we have
 \begin{align} \label{Ew23}
 & \hspace{5mm} \nabla_{\langle r \rangle}(\hat{v}) \\
  &= \nabla_{\langle r \rangle}(v) + \lambda \cdot \sum_{j=0}^r \begin{Bmatrix} r\\j \end{Bmatrix} \cdot \partial_{\langle r-j \rangle}(z) \cdot \nabla_{\langle j \rangle} (u)
  \notag \\
 &  = v_r + \lambda \cdot \sum_{j=0}^r  \begin{Bmatrix} r\\j \end{Bmatrix} \cdot  \partial_{\langle r-j \rangle} (z) \cdot \left(\sum_{i=0}^{l-1} \alpha_{ji}\cdot v_i + \beta_{j} \cdot u + \sum_{i = 1}^{n-l-1} \gamma_{j i} \cdot e_i \right) \notag \\
 & = v_r + \lambda \cdot \sum_{i=0}^{l-1} \sum_{j=0}^r\begin{Bmatrix} r\\j \end{Bmatrix} \cdot \partial_{\langle r-j \rangle}(z) \cdot \alpha_{ji} \cdot  v_i \notag  \\
 & \hspace{5mm}+ \lambda \cdot L_r (z) \cdot  u + \lambda \cdot \sum_{i=1}^{n-l-1} \sum_{j=0}^r \begin{Bmatrix} r\\j \end{Bmatrix}  \cdot \partial_{\langle r-j \rangle}(z) \cdot  \gamma_{j i} \cdot e_i \notag \\
 & = v_r + \lambda  \cdot  \sum_{i=0}^{l-1}  \theta_{ri} \cdot v_i + \lambda \cdot  L_r (z) \cdot u + \lambda \cdot  \sum_{i=1}^{n-l-1} \theta'_{ri} \cdot e_i \notag
 \end{align} 
 for some $\theta_{ri}, \theta'_{ri} \in R$.
 A similar calculation, using (\ref{E3342}), gives
 \begin{align} \label{Ew24}
 \nabla_{\langle l \rangle}(\hat{v}) 
  = \sum_{i=0}^{l-1}a_i \cdot v_i + \lambda \cdot \sum_{i=0}^{l-1}\theta_{li} \cdot v_i + \lambda \cdot  L_l (z) \cdot u + \lambda \cdot \sum_{i=1}^{n-l-1} \theta'_{li} \cdot e_i
 \end{align}
 for some $\theta_{l i}, \theta'_{l i} \in R$.
 
 Now, let us consider the vector 
 \begin{align}
 \hat{w} := \hat{v} \wedge \nabla_{\langle 1 \rangle}(\hat{v}) \wedge \cdots \wedge \nabla_{\langle l \rangle} (\hat{v}) \in {\bigwedge^{l+1}}_{R (\lambda)} \EE^\lambda.
 \end{align}
 Under the natural identification   $\bigwedge^{l+1}_{R (\lambda)} \EE^\lambda = \left(\bigwedge^{l+1}_R \EE \right) \otimes_R R (\lambda)$,
 we can write 
 \begin{align}
 \hat{w} =  w_0 +   w_1 \cdot  \lambda +  w_2  \cdot \lambda^2 + \cdots + w_{l+1} \cdot \lambda^{l+1} 
 \end{align}
for some  $w_0, \cdots, w_{l+1} \in \bigwedge^{l+1}_R\EE$.
Since $\hat{w} |_{\lambda = 0} = 0$, we have $w_0 = 0$.
  By (\ref{Ew23}) and (\ref{Ew24}), 
  the coefficient of $v_0 \wedge v_1 \wedge \cdots \wedge v_{l-1} \wedge u$ for $w_1$
 is given by 
  \begin{align}
 &\hspace{5mm}\left(\sum_{r=0}^{l-1} v_0 \wedge \cdots \wedge v_{r-1} \wedge \left( L_r (z)\cdot u\right) \wedge v_{r+1} \wedge \cdots \wedge  v_{l-1} \wedge \sum_{j=0}^{l-1} a_j  \cdot v_j \right) \\
 & \hspace{7mm}  +  v_0 \wedge \cdots \wedge v_{l-1} \wedge \left( L_l (z) \cdot u\right) 
\notag  \\
 & = \left( \sum_{r=0}^{l-1}a_r \cdot  L_r (z) \cdot v_0 \wedge \cdots \wedge v_{r-1} \wedge u \wedge v_{r+1} \wedge \cdots \wedge v_{l-1} \wedge v_r  \right) \notag \\
 & \hspace{7mm}+ L_l (z) \cdot  v_0 \wedge \cdots \wedge v_{l-1} \wedge u \notag \\
 & = \left(L_l (z)-\sum_{r=0}^{l-1}a_r \cdot  L_r (z)\right) \cdot  v_0 \wedge \cdots \wedge v_{l-1} \wedge u \notag  \\
 & = L (z) \cdot v_0 \wedge \cdots \wedge v_{l-1} \wedge u \left(\neq 0 \right). \notag
 \end{align}
Hence,  the coefficient of $v_0 \wedge v_1 \wedge \cdots \wedge v_{l-1} \wedge u$ for $\hat{w}$ is a nonzero  polynomial in $R [\lambda]$ whose degree is at most $l+1$ and whose constant term is $0$.
 Since $R^{\M+1}$ has at least $n \left(>l \right)$ nonzero elements, there exists an element $\lambda_0 \in R^{\M+1}$ which is not a zero of that polynomial.
  Then, the element  $\overline{v} := v + \lambda_0 \cdot  z \cdot  u \in \EE$
 satisfies  $\overline{v} \wedge \nabla_{\langle 1 \rangle} (\overline{v}) \wedge \cdots \wedge \nabla_{\langle l \rangle}(\overline{v}) \neq 0$.
 By repeating the procedure for constructing $\overline{v}$ using $v$,
 we obtain an $\M$-cyclic vector  of $(\EE, \nabla)$.
 This completes the proof of this assertion. 
 \end{proof}

\vspace{10mm}
\section{Dormant differential modules} \label{SS444} \SSP

This section deals with 
the $p^{\M+1}$-curvature of a differential module of  level $\M \geq 0$.
In particular, we focus on differential modules with vanishing $p^{\M+1}$-curvature, which will be  called 
{\it dormant} differential modules.
At the end of this section, we provide a functorial construction of duality between 
dormant pinned differential modules of rank $n$ (with $0 < n < p^{\M+1}$) and  those of rank  $p^{\M+1} -n$ (cf. Theorem \ref{T4}, Corollary \ref{C003}).

\LSP
\subsection{$p^{\M+1}$-curvature and  dormant differential modules} \label{SS99}

 Let us fix an integer $\M\geq 0$ and  
 a  field $K$ of characteristic $p$.
  In the following discussion, let us consider  an $\M$-derivation on $K$.
Since a perfect field of characteristic $p$
 has only the zero derivation, we should impose the condition  that  $K \neq K^{(1)} := \{ a^p \, | \, a \in K \}$.
 In particular, suppose here 
 that $[K: K^{(1)}] = p$ and there exists a discrete valuation ring  $R$
whose fraction field coincides with $K$.
 Examples of  fields $K$ satisfying this condition
   are 
  $k (t)$ and $k (\!(t)\!)$ with $k$ a perfect field of characteristic $p$.

For each integer $l \geq 0$, we shall  denote by  $R^{(l)}$  
(resp., $K^{(l)}$)
 the subring  of  $R$
  (resp., the subfield of $K$)  
consisting of elements $a^{p^l}$ for  $a \in R$.
  Let us choose a uniformizer $t \in R$.
 Then, $K$
    has basis $1, t, \cdots, t^{p^{\M+1}-1}$ over $K^{(\M+1)}$,
 and it  admits an  extended $\M$-derivation $\partial_{\langle \bullet \rangle} := \{ \partial_{\langle j \rangle} \}_{l \in \mbZ_{\geq 0}}$ (resp., extended $\M$-log derivation $\breve{\partial}_{\langle \bullet \rangle}:= \{ \breve{\partial}_{\langle j \rangle} \}_{j \in \mbZ_{\geq 0}}$) relative to $K^{(\M+1)}$ given  by 
 \begin{align} \label{reow98}
 \partial_{\langle j \rangle}(t^n) := q_j ! \cdot \binom{n}{j} \cdot t^{n-j} \ \left(\text{resp.,} \  \breve{ \partial}_{\langle j \rangle}(t^n) := q_j ! \cdot \binom{n}{j} \cdot t^{n}  \right)
 \end{align}
 for every $j, n \in \mbZ_{\geq 0}$.
For each $l \in \mbZ_{\geq 0}$,  $R^{(l)}$
  coincides with
 $R^l \left(= \bigcap_{j=0}^{l-1}\mr{Ker}(\partial_{\langle p^j \rangle}) =  \bigcap_{j=0}^{l-1} \mr{Ker}(\breve{\partial}_{\langle p^j \rangle})\right)$, 
 and the collection $(R^{(l)}, \{ \partial_{j p^l} \}_{0 \leq j \leq p^{\M-l}})$ forms an $(\M-l)$-differential ring over $R^{(\M+1)}$. 
 Moreover,  the  collection $\partial_{\langle \bullet \rangle} := \{ \partial_{\langle j \rangle} \}_j$ (resp., $\breve{\partial}_{\langle \bullet \rangle} := \{ \breve{\partial}_{\langle j \rangle} \}_j$) is restricted   to an extended $\M$-derivation (resp., extended $\M$-log derivation)  on $R$ relative to $R^{(\M+1)}$, which  we express in the same notation. 
For simplicity, we write 
\begin{align} \label{wqo9}
\mcD_{S}^{(\M)} := \mcD^{(\M)}_{(S, \partial_{\langle \bullet \rangle})} \hspace{3mm}  \text{and} \hspace{3mm}
\breve{\mcD}_{S}^{(\M)} := \mcD^{(\M)}_{(S, \breve{\partial}_{\langle \bullet \rangle})},
\end{align}
 where $S \in \{ R, K \}$.

Now, let ``$\dot{(-)}$" denote either the absence or presence of ``$\breve{(-)}$".
Each $\dot{\mcD}_R^{(\M)}$-module structure $\nabla$ on an $R$-module  $\EE$ naturally extends  to 
a $\dot{\mcD}_K^{(\M)}$-module  structure $\nabla_{\otimes K}$ on  $K \otimes_R \EE$.
The assignment $(\EE, \nabla) \mapsto (K \otimes_R \EE, \nabla_{\otimes K})$ defines a functor
\begin{align} \label{efps39}
\dot{\kappa} : \mfM \mfo \mfd (\dot{\mcD}_R^{(\M)}) \migi \mfM \mfo \mfd (\dot{\mcD}_K^{(\M)}).
\end{align}

Next, let $S \in \{ R, K \}$.
For a $\mcD_S^{(\M)}$-module $(\EE,  \nabla)$, 
 the collection $\{ t^j  \cdot \nabla_{\langle j \rangle} \}_j$ determines  
 a structure of $\breve{\mcD}_{S}^{(\M)}$-module on $\EE$.
Conversely, suppose that we are given a 
$\breve{\mcD}^{(\M)}_{S}$-bundle  $(\EE, \breve{\nabla})$
such that, for every $j \in \mbZ_{\geq 0}$, 
$\breve{\nabla}_{\langle j \rangle}$ may be expressed as 
$\breve{\nabla}_{\langle j \rangle} = t^j \cdot \nabla_{\langle j \rangle}$ for some
 $\nabla_{\langle j \rangle} \in \mr{End}_{S^{(\M+1)}} (\EE)$.
 Then,  the collection $\{ \nabla_{\langle j \rangle} \}_j$ determines  
a structure of $\mcD_{S}^{(\M)}$-module on $\EE$.
The assignment $(\EE, \nabla) \mapsto (\EE, \{ t^j \cdot \nabla_{\langle j \rangle} \}_j)$
defines a functor
\begin{align} \label{ee89d}
\eta_S :\mfM \mfo \mfd (\mcD_S^{(\M)}) \migi \mfM \mfo \mfd (\breve{\mcD}_S^{(\M)}),
\end{align}
and it becomes an equivalence of categories when $S = K$.
Moreover, the following square diagram of categories  is $1$-commutative:
\begin{align} \label{E02dd}
\vcenter{\xymatrix@C=46pt@R=36pt{
\mfM \mfo \mfd (\mcD^{(\M)}_R)\ar[r]^-{\kappa} \ar[d]_-{\eta_R} & \mfM \mfo \mfd (\mcD_K^{(\M)})\ar[d]_-{\wr}^{\eta_K} 
\\
\mfM \mfo \mfd (\breve{\mcD}_R^{(\M)}) \ar[r]_-{\breve{\kappa}} & \mfM \mfo \mfd (\breve{\mcD}_K^{(\M)}).
}}
\end{align}

\begin{rema} \label{roafos}
Let  $a$  be an integer and $(\EE, \nabla)$  a $\breve{\mcD}_{R}^{(\M)}$-module such that the $R$-module $\EE$ is  free.
Then, 
$\nabla$ naturally induces  a $\breve{\mcD}_{R}^{(\M)}$-module  structure on the $R$-module $t^a \cdot  \EE \subseteq K \otimes_R \EE$, which will be denoted by $\nabla |_{t^a \cdot  \EE}$.
 \end{rema}
\SSP

\begin{rema} \label{roeeafos}
Let $\EE$ be a $S^{(\M+1)}$-module.
Then, since $\Xi^{\uparrow (\M+1)}$ (cf. (\ref{E309})) is compatible with   $\eta_S$ (cf. (\ref{ee89d})),
 $\nabla^\mr{can}_{(S, \breve{\partial}_{\langle \bullet \rangle}), \EE}$ comes from the  $\mcD_S^{(\M+1)}$-module structure $\nabla^\mr{can}_{(S, \partial_{\langle \bullet \rangle}), \EE}$ via $\eta_S$.
 This  means  that the following equality holds:
\begin{align} \label{gsoa3}
\nabla^\mr{can}_{(S, \breve{\partial}_{\langle \bullet \rangle}), \EE} = \{ t^j (\nabla^\mr{can}_{(S, \partial_{\langle \bullet \rangle}), \EE})_{\langle j \rangle} \}_j.
\end{align}
 \end{rema}
\SSP

Let $(\EE, \nabla)$
 be a $\mcD_S^{(\M)}$-module (resp., a $\breve{\mcD}_S^{(\M)}$-module), where $S \in \{ R, K \}$.
 The {\bf $p^{\M+1}$-curvature} of $(\EE, \nabla)$
 is defined as 
 \begin{align}
 {^p}\psi_{(\EE, \nabla)} := \nabla_{\langle p^{\M+1}\rangle} 
    \in \mr{End}_{S^{(\M+1)}}(\EE).
 \end{align}
It can immediately be seen that $ {^p}\psi_{(\EE, \nabla)}$
  belongs to $\mr{End}_{R}(E)$.

\SSP
\bde \label{D098}
With the above notation, 
we shall say that 
$(\EE, \nabla)$
  is  {\bf dormant} if ${^p}\psi_{(\EE, \nabla)} = 0$.
   Also,  a pinned $\dot{\mcD}_S^{(\M)}$-module $(\EE, \nabla, v)$ is  called {\bf dormant} if $(\EE, \nabla)$ is dormant.
\ede
\SSP

\begin{rema}
\label{uin9}
Let us consider the case of $\M=0$.
 Then, the   $p^1$-curvature $ {^p}\psi_{(\EE, \nabla)}$ of a $\mcD_S^{(0)}$-module $(\EE, \nabla)$ coincides with  the  $p$-curvature of a differential module over the differential ring $(S, \partial_{\langle 1 \rangle})$  (cf. Remark \ref{R093} and ~\cite[\S\,13.1]{vdPS}).
By  an argument similar to  the proof of ~\cite[\S\,5,  Theorem (5.1)]{Kal},
the functors  $\Xi^{\downarrow (1)}$ and $\Xi^{\uparrow (1)}$ define an equivalence of categories
   \begin{align} \label{Efer22}
\begin{pmatrix}
\text{the category of } \\
\text{dormant $\mcD_S^{(0)}$-modules} 
\end{pmatrix}
\isom \begin{pmatrix}
\text{ the category of } \\
\text{$S^{(1)}$-modules} \\
\end{pmatrix}.
\end{align}
This equivalence for 
 $S = K$ can be found in  ~\cite[\S\,13.1, Lemma 13.2]{vdPS}.
 
\end{rema}
\SSP

The following assertion may be regarded as a version of ~\cite[\S\,2.3, Th\'{e}or\`{e}m 2.3.6]{PBer2} for higher-level differential modules.
In particular, 
from the equivalence  (\ref{Eq38})  in that proposition,
$\mcD_R^{(\M)}$-modules are essentially equivalent to $\mcD_R^{(0)}$-modules, i.e., differential modules in the classical sense.
(But, as we will  see in the next section,  this is not true for the logarithmic case, i.e.,  $\breve{\mcD}_R$-modules.)

\SSP
\bpr \label{P099}
Let $l$ be an integer with $0 \leq l \leq \M$ and let $S \in \{ R, K \}$.
Then, the following assertions hold:
\begin{itemize}
\item[(i)]
The functors $\Xi^{\downarrow (l)}$ (cf. (\ref{E323}))
and $\Xi^{\uparrow (l)}$ (cf. (\ref{E309})) 
define an equivalence of categories
\begin{align} \label{Eq38}
\mfM \mfo \mfd (\mcD_{S}^{(\M)}) \isom \mfM \mfo \mfd (\mcD^{(\M-l)}_{S^{(l)}}).
\end{align}
\item[(ii)]
Let $(\EE, \nabla)$ be a $\mcD_{S^{(l)}}^{(\M-l)}$-module.
Then, the $p^{\M-l+1}$-curvature ${^p}\psi_{(\EE, \nabla)}$ of $(\EE, \nabla)$ and  
the $p^{\M+1}$-curvature ${^p}\psi_{(F^{(l)*}(\EE), F^{(l)*}(\nabla))}$ of $(F^{(l)*}(\EE), F^{(l)*}(\nabla))$ $\left(= \Xi^{\uparrow (l)} ((\EE, \nabla)) \right)$ satisfy
the equality
\begin{align}
{^p}\psi_{(F^{(l)*}(\EE), F^{(l)*}(\nabla))} = \mr{id}_{R} \otimes {^p}\psi_{(\EE, \nabla)} \in \mr{End}_R (F^{(l)*}(\EE)).
\end{align}
In particular, $(\EE, \nabla)$ is dormant if and only if $(F^{(l)*}(\EE), F^{(l)*}(\nabla))$ is dormant.
\end{itemize}
\epr
\begin{proof}
First,  we shall prove  assertion (i)
 by induction on $l$.
The base step, i.e., $l =0$, is trivial.
For the induction step, 
let us take a $\mcD_{S}^{(\M)}$-module $(\EE, \nabla)$.
This $\mcD_{S}^{(\M)}$-module induces  a  $\mcD_{S^{(1)}}^{(\M-1)}$-module of the form $(\EE^1, \nabla^1) \left(= \Xi^{\downarrow (1)} ((\EE, \nabla)) \right)$.
Since $(\EE^1)^{l-1} = \EE^l$, 
 the induction hypothesis implies that 
the morphism 
$\tau^{\downarrow \uparrow (l-1)}_{(\EE^1, \nabla^1)} : F_1^{(l-1)*}(\EE^l) \migi E^1$
(cf. (\ref{Err4})) is an isomorphism.
On the other hand, 
 the $p$-curvature of the differential module $(\EE, \nabla_{\langle 1 \rangle})$ vanishes, so it follows from the equivalence of categories (\ref{Efer22}) that
the morphism $\tau_1 : F^{(1)*} (\EE^1) \isom \EE$ extending the inclusion $\EE^1 \migiincl \EE$
is an isomorphism.
Hence, $\tau^{\downarrow \uparrow (l)}_{(\EE, \nabla)}$ turns out to be an isomorphism because it coincides with   the composite  isomorphism
\begin{align}
F^{(l)*}(\EE^l) \left(= F^{(1)*}(F^{(l-1)*}_1 (E^l)) \right) \xrightarrow{\Xi^{\uparrow (1)*}(\tau^{\downarrow \uparrow (l-1)}_{(\EE^1, \nabla^1)})} F^{(1)*} (E^1) \xrightarrow{\tau_1} \EE.
\end{align}
Since $\tau^{\downarrow \uparrow (l)}_{(\EE, \nabla)}$ is functorial with respect to $(\EE, \nabla)$, we see that the composite functor $\Xi^{\uparrow (l)}\circ \Xi^{\downarrow (l)}$ is  isomorphic to the identity functor of $\mfM \mfo \mfd (\mcD_S^{(\M)})$.

Next, let  $(E, \nabla)$ be a $\mcD_{S^{l}}^{(\M-l)}$-module.
By applying $\Xi^{\uparrow (l-1)}$ to  $(E, \nabla)$,
we obtain
 a $\mcD_{S^{(1)}}^{(\M-1)}$-module of the form  $(F^{(l-1)*}_1 (\EE), F^{(l-1)*}_1 (\nabla)) \left(= \Xi^{\uparrow (l-1)} ((\EE, \nabla)) \right)$.
It follows from 
(\ref{Efer22})
 again that the natural morphism
$\tau_2 : F^{(l-1)*}_1 (\EE) \migi F^{(1)*}(F^{(l-1)*}_1 (\EE))^1$ is an isomorphism.
Also, the induction hypothesis implies that
the morphism  $\tau_{(\EE, \nabla)}^{\uparrow \downarrow (l-1)} : \EE \migi F^{(l-1)*}_1 (\EE)^{l-1}$ is an isomorphism.
Hence, $\tau_{(\EE, \nabla)}^{\uparrow \downarrow (l)}$ turns out to be an isomorphism because it coincides with the composite isomorphism
\begin{align}
\EE  \xrightarrow{\tau_{(\EE, \nabla)}^{\uparrow \downarrow (l-1)}}
F^{(l-1)*}_1 (\EE)^{l-1}  \xrightarrow{\Xi^{\downarrow (l-1)}(\tau_2)}
 \left( (F^{(1)*}(F^{(l-1)*}_1 (\EE))^1)^{l-1} = \right) F^{(l)*}(\EE)^l.
\end{align}
Since $\tau_{(\EE, \nabla)}^{\uparrow \downarrow (l)}$ is functorial with respect to $(\EE, \nabla)$,  the composite functor $\Xi^{\downarrow (l)}\circ \Xi^{\uparrow (l)}$ is  isomorphic to the identity functor of $\mfM \mfo \mfd (\mcD_{S^{(l)}}^{(\M-l)})$.
This completes the proof of assertion (i).

Assertion (ii) can be verified immediately from  the definitions of $p^{(-)}$-curvature and the functor $\Xi^{\uparrow (l)}$.
\end{proof}
\SSP

Moreover, by applying the above theorem for $l=\M$ and the equivalence of categories (\ref{Efer2}),
we obtain the following assertion (cf. ~\cite[\S\,3, Corollary 3.2.4]{LeQu}):

\SSP
\bco \label{eoap09}
Let $S \in \{ R, K \}$.
Then, the functors $\Xi^{\downarrow (\M+1)}$ and $\Xi^{\uparrow (\M+1)}$ (i.e., the assignments $(\EE, \nabla) \mapsto E^{\M+1}$ and $\EE' \mapsto (F^{(\M+1)*}(\EE'), \nabla^\mr{can}_{(S, \partial_{\langle \bullet \rangle}), \EE'})$) induce 
an equivalence of categories
  \begin{align} \label{Efer2}
\begin{pmatrix}
\text{the category of } \\
\text{dormant $\mcD_S^{(\M)}$-modules} 
\end{pmatrix}
\isom \begin{pmatrix}
\text{ the category of } \\
\text{$S^{(\M+1)}$-modules} \\
\end{pmatrix}.
\end{align}
\eco

\LSP
\subsection{Dormant pinned $\mcD_R^{(\M)}$-module of rank $p^{\M+1}$} \label{SSd45}

Let  $S \in \{ R, K \}$, and let ``$\dot{(-)}$" denote either the absence or presence of ``$\breve{(-)}$".
Here, let us construct  an example of a dormant pinned $\mcD_S^{(\M)}$-module of rank $p^{\M+1}$.
The $S$-module
  $\dot{P}_S := \dot{\mcD}_{S}^{(\M)}/\dot{\mcD}_{S}^{(\M)} \cdot \dot{\partial}_{\langle p^{\M+1}\rangle}$ has the a $\dot{\mcD}_S^{(\M)}$-module structure $\nabla_{\dot{P}_S}$ induced from the left $\dot{\mcD}_S^{(\M)}$-module structure of $\dot{\mcD}_{S}^{(\M)}$ itself.
  One can verify that $(\dot{P}_S, \nabla_{\dot{P}_S})$ is dormant.
If $\dot{\delta}_{\langle l \rangle}$ ($l =0, \cdots, p^{\M+1} -1$) is the image of $\dot{\partial}_{\langle l \rangle}$ via the quotient $\dot{\mcD}_S^{(\M)} \migisurj \dot{P}_S$, then we have 
$\dot{P}_S = \bigoplus_{l=0}^{p^{\M+1}-1} S \cdot \dot{\delta}_{\langle l \rangle}$.
Hence, the triple 
\begin{align} \label{E4567}
(\dot{P}_S, \nabla_{\dot{P}_S}, v_{\dot{P}_S}),
\end{align}
where $v_{\dot{P}_S} := \dot{\delta}_{\langle 0 \rangle}$,
forms a dormant pinned  $\dot{\mcD}_S^{(\M)}$-module of rank $p^{\M+1}$.

Next, let $(\EE, \nabla, v)$ be 
a dormant pinned $\breve{\mcD}_{S}^{(\M)}$-module. 
The $S$-linear  injection  $S \migiincl  \EE$ given by $a \mapsto a \cdot v$ (for any $a \in S$) extends to a  $\dot{\mcD}_S^{(\M)}$-linear morphism $\widetilde{\nu}_{(\EE, \nabla, v)} : \left(\dot{\mcD}_{S}^{(\M)} \otimes_S S = \right)\dot{\mcD}_{S}^{(\M)} \migi \EE$.
This morphism preserves the $\dot{\mcD}_{S}^{(\M)}$-action.
Since $(\EE, \nabla)$ has vanishing $p^{\M+1}$-curvature,
$\widetilde{\nu}_{(\EE, \nabla, v)}$
 factors through the quotient $\dot{\mcD}_S^{(\M)} \migisurj \dot{P}_S$.
Thus, $\widetilde{\nu}_{(\EE, \nabla, v)}$ induces  a  morphism 
\begin{align} \label{Eppop}
\nu_{(\EE, \nabla, v)} : \dot{P}_S \migisurj \EE,
\end{align}
which forms a morphism of pinned $\dot{\mcD}_{S}^{(\M)}$-modules $(\dot{P}_S, \nabla_{\dot{P}_S}, v_{\dot{P}_S}) \migi (\EE, \nabla, v)$.
It follows from   Proposition \ref{edaP}, (i), that $\nu_{(\mcE, \nabla, v)}$ is surjective.
By combining Theorem \ref{P023} with  the following proposition (in the case where $S = K$ and ``$\dot{(-)}$" denotes the absence of $\breve{(-)}$), we see that
the existence of an $\M$-cyclic vector for a dormant $\mcD_K^{(\M)}$-module  depends  only on the rank of the underlying $K$-vector space.

\SSP
\bpr \label{er45gj}
Let  $n$ be a positive integer and $(\EE, \nabla)$  a dormant $\dot{\mcD}_S^{(\M)}$-module  such that the $S$-module $\EE$ is free and  of rank $n$.
Then, the following assertions hold:
\begin{itemize}
\item[(i)]
If the inequality $n > p^{\M+1}$ holds, then there are no $\M$-cyclic  vectors of $(\EE, \nabla)$.
\item[(ii)]
If the equality $n = p^{\M +1}$ holds and there exists an $\M$-cyclic vector of $(\EE, \nabla)$, then $\nu_{(\EE, \nabla, v)}$ is an isomorphism.
In particular, the isomorphism class of a  dormant  $\dot{\mcD}_S^{(\M)}$-module whose underlying $S$-module is free and  of rank $p^{\M+1}$  is uniquely    determined, i.e., the class represented by $(\dot{P}_S, \nabla_{\dot{P}_S}, v_{\dot{P}_S})$.
\end{itemize}
\epr
\begin{proof}
The assertions follow immediately from the surjectivity of the morphism $\nu_{(\mcE, \nabla, v)}$ for each $\M$-cyclic vector  $v$ of $(\EE, \nabla)$.
\end{proof}

\SSP
\begin{exa} \label{E3432}
Let us consider the dual of $\dot{P}_S \left( = \bigoplus_{l=0}^{p^{\M+1}-1} S \cdot \dot{\delta}_{\langle l \rangle}\right)$.
Denote the dual basis of $\dot{\delta}_{\langle 0 \rangle}, \cdots, \dot{\delta}_{\langle p^{\M+1} -1 \rangle}$  by $\dot{\delta}^\vee_{\langle 0 \rangle}, \cdots, \dot{\delta}^\vee_{\langle p^{\M+1}-1\rangle}$.
 From the definition of $\nabla_{\dot{P}_S}$,  
 the element $v_{\dot{P}_S^\vee} := \dot{\delta}^\vee_{\langle p^{\M+1}-1\rangle}$
   defines an $\M$-cyclic vector of 
 the dual $\dot{\mcD}_S^{(\M)}$-module $(\dot{P}_S^\vee, \nabla_{\dot{P}_S}^\vee)$ of $(\dot{P}_S, \nabla_{\dot{P}_S})$.
Hence, we obtain the dormant pinned $\dot{\mcD}_S^{(\M)}$-module 
 \begin{align} \label{fao332}
  (\dot{P}_S^\vee, \nabla_{\dot{P}_S}^\vee, v_{\dot{P}_S^\vee}).
 \end{align}
Since  the free  $S$-module $\dot{P}_S^\vee$ is of rank $p^{\M+1}$,  it follows from Proposition \ref{er45gj}, (ii), that the induced morphism
\begin{align} \label{fao333}
\nu_{(\dot{P}_S^\vee, \nabla_{\dot{P}_S}^\vee, v_{\dot{P}_S^\vee})} : (\dot{P}_S, \nabla_{\dot{P}_S}, v_{\dot{P}_S})  \migi  (\dot{P}_S^\vee, \nabla_{\dot{P}_S}^\vee, v_{\dot{P}_S^\vee})
\end{align}
defines  an isomorphism of pinned $\dot{\mcD}_S^{(\M)}$-modules.
\end{exa}

\LSP
\subsection{Duality of dormant pinned $\mcD_R^{(\M)}$-modules} \label{SS45dd}

Let us keep the above notation.
Also, let us denote by $\nabla_{\mr{Ker}}$ the $\dot{\mcD}_S^{(\M)}$-module structure on 
$\mr{Ker}(\nu_{(\EE, \nabla, v)})$ obtained by restricting $\nabla_{\dot{P}_S}$.
  If $(\EE^\BB, \nabla^\BB)$
  is  the dual of the $\dot{\mcD}_{S}^{(\M)}$-module $(\mr{Ker}(\nu_{\msE^\circledast}), \nabla_{\mr{Ker}})$, then it
   has 
 vanishing $p^{\M+1}$-curvature.
We shall write  $v^\BB$ for the element of $\EE^\BB$ defined to be  the image of $1$ via the dual of 
the composite
\begin{align}
\mr{Ker}(\nu_{(\EE, \nabla, v)}) \xrightarrow{\mr{inclusion}} \dot{P}_S \left(= \bigoplus_{l=0}^{p^{\M+1} -1}  S \cdot  \dot{\delta}_{\langle l \rangle}\right) \migi \left( S \cdot \dot{\delta}_{\langle p^{\M+1} -1\rangle}=\right) S,
\end{align}
where the second arrow  denotes the projection to the last factor.

\SSP
\ble \label{L004}
The  element $v^\BB$ forms an $\M$-cyclic vector of $(\EE^\BB, \nabla^\BB)$.
\ele
\begin{proof}
For each $j \in \mbZ_{\geq 0}$, we write 
\begin{align}
\dot{P}_{S, j} := \mr{Im}\left( \dot{\mcD}_{S, < j}^{(\M)} \migiincl \dot{\mcD}_{S}^{(\M)} \migisurj \dot{P}_S\right) \left(= \bigoplus_{l=0}^{j-1} S \cdot \dot{\delta}_{\langle l \rangle} \right).
\end{align}
Let us set $h$ to be  the composite
\begin{align}
h : \dot{P}_{S, n} \xrightarrow{\mr{inclusion}} \dot{P}_S\xrightarrow{\nu_{(\EE, \nabla, v)}} \EE.
\end{align}
By the definition of $\nu_{(\EE, \nabla, v)}$, 
we have
\begin{align}
h ( \dot{\delta}_{\langle j \rangle}) =  h ((\nabla_{\dot{P}_S})_{\langle j \rangle} (\dot{\delta}_{\langle 0 \rangle})) = \nabla_{\langle j \rangle} (h (\dot{\delta}_{\langle 0 \rangle})) = \nabla_{\langle j \rangle} (v)
\end{align}
 for every $j = 0, \cdots, n-1$.
On the other hand,  recall that  $\dot{P}_{S, n}$ and $\EE$  are generated by the sets  $\{ \dot{\delta}_{\langle j \rangle} \}_{j=0}^{n-1}$ and $\{ \nabla_{\langle j \rangle} (v) \}_{j=0}^{n-1}$ respectively.
Hence,  
$h$ turns out to be   an isomorphism.
This implies that the composite 
\begin{align}
\lambda : \mr{Ker}(\nu_{(\EE, \nabla, v)}) \xrightarrow{\mr{inclusion}} \dot{P}_S \migisurj  \dot{P}_S/\dot{P}_{S, n}
\end{align}
 is an isomorphism.
 Since the element $v_{\dot{P}_S^\vee}$ of $\dot{P}^\vee_S$ (cf. Example \ref{E3432}) forms an $\M$-cyclic vector,
 the  elements $\left(v_{\dot{P}_S^\vee} =  \right)\dot{\partial}_{\langle p^{\M+1}-1 \rangle}^\vee, \cdots, \dot{\partial}^\vee_{\langle n\rangle}$ generate an $S$-submodule $(\dot{P}_S/\dot{P}_{S, n})^\vee \subseteq \dot{P}^\vee_S$.
In particular, 
the elements $\lambda^\vee (\dot{\partial}_{\langle p^{\M+1}-1 \rangle}), \cdots, \lambda^\vee (\dot{\partial}_{\langle n \rangle})$ generate $\EE^\BB$,   where $\lambda^\vee$ denotes the dual $(\dot{P}_S/\dot{P}_{S, n})^\vee \isom \EE^\BB  \left(=\mr{Ker}(\nu_{(\EE, \nabla, v)})^\vee  \right)$ of $\lambda$.
 On the other hand, it follows from the various definitions involved that the equality $\nabla^\BB_{\langle j \rangle} (v^\BB) = \lambda^\vee (\dot{\partial}^\vee_{\langle p^{\M+1}-1-j \rangle})$  holds for every $j =0, \cdots, p^{\M+1}-n-1$.
 Thus,  $v^\BB$ turns out to  form an $\M$-cyclic vector of $(\EE^\BB, \nabla^\BB)$.
 This completes the proof of the assertion.
\end{proof}
\SSP

By the above lemma, 
 we obtain a dormant pinned $\dot{\mcD}_{S}^{(\M)}$-module 
\begin{align} \label{E456d} 
(\EE^\BB, \nabla^\BB, v^\BB)
\end{align}
of rank $p^{\M}-n$, which we call  the {\bf dual} of $(\EE, \nabla, v)$.

Next, let us take a morphism  $f : (\EE, \nabla, v) \migi (\EE', \nabla', v')$ between dormant  pinned $\dot{\mcD}_S^{(\M)}$-modules.
The construction of $\nu_{(-)}$ yields  the equality $f \circ \nu_{(\EE, \nabla, v)} = \nu_{(\EE', \nabla', v')}$.
Hence, $f$ is restricted  to 
 the  inclusion $f_{\mr{Ker}} : \mr{Ker}(\nu_{(\EE, \nabla, v)}) \migiincl \mr{Ker}(\nu_{(\EE', \nabla', v')})$.
Taking its dual gives  a morphism of pinned $\dot{\mcD}_{S}^{(\M)}$-modules
\begin{align}
f^\BB : (\EE'^\BB, \nabla'^\BB, v'^\BB) \migi (\EE^\BB, \nabla^\BB, v^\BB).
\end{align}

Here, we shall  write
\begin{align}
\mfM \mfo \mfd (\dot{\mcD}_{S}^{(\M)})^\circledast
\end{align}
for the category of dormant pinned $\dot{\mcD}_{S}^{(\M)}$-modules.

\SSP
\bt \label{T4}
\begin{itemize}
\item[(i)]
The assignments $(\EE, \nabla, v) \mapsto (\EE^\BB, \nabla^\BB, v^\BB)$ and $f \mapsto f^\BB$ constructed above define a self-equivalence 
\begin{align} \label{bamoz9}
\dot{\Dual}_S : \mfM \mfo \mfd (\dot{\mcD}_{S}^{(\M)})^\circledast \isom \mfM \mfo \mfd (\dot{\mcD}_{S}^{(\M)})^\circledast
\end{align}
of the category $\mfM \mfo \mfd (\dot{\mcD}_{S}^{(\M)})^\circledast$ with $\dot{\Dual}_S \circ \dot{\Dual}_S = \mr{id}$.
In particular,  for each pinned $\dot{\mcD}_S^{(\M)}$-module $(\EE, \nabla, v)$,  there exists an isomorphism of pinned $\dot{\mcD}_S^{(\M)}$-modules
\begin{align} \label{bmk2}
(\EE, \nabla, v) \isom  (\EE^{\BB\BB}, \nabla^{\BB \BB}, v^{\BB \BB}).
\end{align}
\item[(ii)]
The following diagram of functors is $1$-commutative:
\begin{align} \label{E002fg7}
\vcenter{\xymatrix@C=46pt@R=36pt{
 \mfM \mfo \mfd (\mcD_{R}^{(\M)})^\circledast  \ar[d]_-{\wr}^{\tiny{\Dual}_R} \ar[r]^{\eta_R^\circledcirc}  & \mfM \mfo \mfd (\breve{\mcD}_{R}^{(\M)})^\circledast   \ar[d]_-{\wr}^-{\breve{\tiny{\Dual}}_R} \ar[r]^{(\eta_K^{-1}\circ \breve{\kappa})^\circledcirc} & \mfM \mfo \mfd (\mcD_{K}^{(\M)})^\circledast \ar[d]_-{\wr}^-{\tiny{\Dual}_K} 
 \\
  \mfM \mfo \mfd (\mcD_{R}^{(\M)})^\circledast  \ar[r]_-{\eta_R^\circledcirc}  &\mfM \mfo \mfd (\breve{\mcD}_{R}^{(\M)})^\circledast  \ar[r]_-{(\eta_K^{-1}\circ \breve{\kappa})^\circledcirc} &   \mfM \mfo \mfd (\mcD_{K}^{(\M)})^\circledast,
  }}
\end{align}
where $\eta_R^\circledcirc$  and $(\eta_K^{-1}\circ \breve{\kappa})^\circledcirc$ denote the functors induced, via restriction, from  $\eta_R$  and $\eta_K^{-1}\circ \breve{\kappa}$, respectively.
\end{itemize}
\et
\begin{proof}
Let us prove assertion (i).
Let $(\EE, \nabla, v)$ be  
a  dormant pinned $\mcD_{S}^{(\M)}$-module.
Then,  it induces  
the following  short exact sequence: 
\begin{align} \label{e45}
0 \longmigi \mr{Ker}(\nu_{(\EE, \nabla, v)}) \xrightarrow{\iota}\dot{P}_S \xrightarrow{\nu_{(\EE, \nabla, v)}} \EE \longmigi 0, 
\end{align}
where we write  $\iota$  for the natural inclusion.
The dual of this sequence fits into the following   diagram:
\begin{align} \label{E0023d7}
\vcenter{\xymatrix@C=56pt@R=36pt{
0 \ar[r]  & \mr{Ker}(\nu_{(\EE^\BB, \nabla^\BB, v^\BB)}) \ar[r]^-{}  & \dot{P}_S \ar[d]_-{\nu_{(\dot{P}_S, \nabla_{\dot{P}_S}, v_{\dot{P}_S})}} \ar[r]^-{\nu_{(\EE^\BB, \nabla^\BB, v^\BB)}} & \EE^\BB\ar[d]_-{\wr}^-{\mr{id}} \ar[r] & 0
 \\
 0  \ar[r] & \EE^\vee \ar[r]_-{\nu_{(\EE, \nabla, v)}^\vee} & \dot{P}_S^\vee\ar[r]_-{\iota^\vee} &  \EE^\BB  \ar[r] & 0,
  }}
\end{align}
where the right-hand square is commutative because of Proposition \ref{edaP}, (ii), together with the equality $(\iota^\vee \circ \nu_{(\dot{P}_S, \nabla_{\dot{P}_S}, v_{\dot{P}_S})}) (v_{\dot{P}_S}) = (\mr{id} \circ \nu_{(\EE^\BB, \nabla^\BB, v^\BB)}) (v_{\dot{P}_S})$.
Hence, this diagram induces an isomorphism $\mr{Ker}(\nu_{(\EE^\BB, \nabla^\BB, v^\BB)}) \isom \EE^\vee$.
By taking the dual of this isomorphism, we obtain 
 an isomorphism $\xi : \EE \isom \EE^{\BB\BB}$.
Note that
$\nu_{(\dot{P}_S, \nabla_{\dot{P}_S}, v_{\dot{P}_S})}$ preserves the $\dot{\mcD}_S^{(\M)}$-action and 
 it is compatible with the respective projections onto $S$, i.e., $v_{\dot{P}_S^\vee}$ and $v_{\dot{P}_S}^\vee$.
This implies 
 that
$\xi$ defines an isomorphism  of pinned $\dot{\mcD}_S^{(\M)}$-modules $(\EE, \nabla, v)\isom  (\EE^{\BB}, \nabla^{\BB \BB}, v^{\BB \BB})$.
Moreover, the formation of this isomorphism is functorial with respect to $(\EE, \nabla, v)$, so it induces an isomorphism  $\mr{id}_{\mfM \mfo \mfd (\dot{\mcD}^{(\M)}_S)^\circledast} \isom \dot{D}_S \circ \dot{D}_S$ of functors.
This completes the proof of  assertion (i).

Assertion (ii) follows  immediately from the construction of  $\dot{\Dual}_S$.
\end{proof}
\SSP

For each integer $n$ with $1\leq n \leq p^\M-1$, we shall write 
\begin{align} \label{ggie8}
\mr{Mod} (\dot{\mcD}_{S}^{(\M)})^\circledast_{n}
\end{align}
 for the set of isomorphism classes of dormant  pinned $\dot{\mcD}_S^{(\M)}$-modules of rank $n$.
 In the case of $n =1$, the set  $\mr{Mod} (\dot{\mcD}_{S}^{(\M)})^\circledast_{1}$  forms a group with a binary operation given by $((\EE, \nabla, v), (\EE', \nabla', v')) \mapsto (\EE \otimes \EE, \nabla \otimes \nabla', v \otimes v')$ (cf. Example \ref{eora921} below). 

The  following assertion is a direct consequence of the above theorem.

\SSP
\bco \label{C003}
The assignment $(\EE, \nabla, v) \mapsto (\EE^\BB, \nabla^\BB, v^\BB)$ defines a bijection of sets
\begin{align}\label{ew234}
[\dot{\Dual}_{S}]_n : \mr{Mod} (\dot{\mcD}_{S}^{(\M)})^\circledast_{n} \isom \mr{Mod} (\dot{\mcD}_{S}^{(\M)})^\circledast_{p^{\M+1} -n}
\end{align}
satisfying $[\dot{\Dual}_{S}]_{p^{\M+1}-n} \circ [\dot{\Dual}_{S}]_n = \mr{id}$.
\eco
\SSP

\begin{exa} \label{eora921}
Let us  examine  the group $\mr{Mod} (\mcD_{S}^{(\M)})^\circledast_{1}$.
For each  element $a \in S^\times$,
denote by $\mu_a$ the automorphism of the trivial $S$-module $S$ 
  given by multiplication by $a$.
There exists a unique $\mcD_{S}^{(\M)}$-module structure $\mu_{a *}(\partial_{\langle  \bullet \rangle})$ on $S$ such that
$\mu_a$ becomes an isomorphism of $\mcD_S^{(\M)}$-modules $(S, \partial_{\langle \bullet \rangle}) \isom (S, \mu_{a *}(\partial_{\langle  \bullet \rangle}))$.
Then, the triple
 \begin{align} \label{e23}
 (S, \mu_{a *}(\partial_{\langle  \bullet \rangle}), 1)
 \end{align}
 forms a dormant pinned $\mcD_{S}^{(\M)}$-module of rank $1$.
Since $\Xi^{\uparrow (\M+1)}(S^{(\M+1)}) = (S, \partial_{\langle \bullet \rangle})$, the equivalence of categories (\ref{Efer2}) shows that
 any  dormant pinned $\mcD_S^{(\M)}$-module of rank $1$ is isomorphic to $(S, \mu_{a *}(\partial_{\langle  \bullet \rangle}), 1)$ for some $a \in S$. 
Also,
 for $a, b \in S^\times$,  
 $ (S, \mu_{a *}(\partial_{\langle  \bullet \rangle}), 1) \cong  (S, \mu_{b *}(\partial_{\langle  \bullet \rangle}), 1)$ ($a, b \in S^\times$) if and only if there exists $c \in (S^{(\M+1)})^\times$ satisfying $b = a \cdot c$.
Hence,  the assignment $a \mapsto (S, \mu_{a *}(\partial_{\langle  \bullet \rangle}), 1)$ gives a well-defined  bijection of sets
\begin{align} \label{eaph}
S^\times/(S^{(\M + 1)})^\times \isom  \mr{Mod} (\mcD_{S}^{(\M)})^\circledast_{1}.
\end{align}
Since $\mu_a \circ \mu_b = \mu_{a \cdot b}$,  we have $(S, \mu_{a *}(\partial_{\langle  \bullet \rangle}), 1) \otimes (S, \mu_{b *}(\partial_{\langle  \bullet \rangle}), 1) \cong (S, \mu_{a \cdot  b *}(\partial_{\langle  \bullet \rangle}), 1)$.
This implies that (\ref{eaph}) becomes an isomorphism of groups.

Moreover, by composing  with 
$[\Dual_S]_1$,
we obtain a  bijection
\begin{align}
S^\times/(S^{(\M+1)})^\times \isom  \mr{Mod} (\mcD_{S}^{(\M)})^\circledast_{p^{\M+1} -1}.
\end{align}
\end{exa}

\vspace{10mm}
\section{Residues and exponents of log differential modules}\label{Swe46} \SSP

Let $\M$, $K$, $R$, $t$, and $\partial_{\langle \bullet \rangle}$ be as in the previous section.
This section discusses the residue and  exponent   of a  dormant $\breve{\mcD}_R^{(\M)}$-module.
In particular, we examine the exponent of a dormant $\breve{\mcD}_R^{(\M)}$-module admitting an $\M$-cyclic vector (cf. Proposition \ref{C001}).

\LSP
\subsection{Residue of a $\breve{\mcD}_R^{(\M)}$-module} \label{SS9}

Let $(\EE, \nabla)$
 be a dormant   $\breve{\mcD}^{(\M)}_{R}$-bundle such that the $R$-module $\EE$ is free and of rank $n > 0$.
For each $l=0, \cdots, \M$,
the pair 
\begin{align} \label{fmxo}
(E^l, \overline{\nabla}^l),
\end{align}
 where $\overline{\nabla}^l := \nabla_{\langle p^{l}\rangle} |_{E^l}$,  determines 
a dormant $\breve{\mcD}_{R^{(l)}}^{(0)}$-module (cf. Remark \ref{R093}).

From the equivalence of categories (\ref{Efer2}), 
the injective morphism $\tau^{\downarrow \uparrow (\M+1)}_{(\EE, \nabla)} : F^{(\M +1)*}(\EE^{\M+1}) \migi \EE$ (cf (\ref{Err4})) becomes bijective after tensoring with $K$.
Hence, the cokernel 
\begin{align}
\mr{Res}(\nabla) := \mr{Coker}(\tau^{\downarrow \uparrow (\M +1)}_{(\EE, \nabla)})
\end{align}
of $\tau^{\downarrow \uparrow (\M+1)}_{(\EE, \nabla)}$  is of finite length.
 We will refer to $\mr{Res}(\nabla)$ as the {\bf residue} of $\nabla$.

Note that there exists a natural sequence of inclusions
\begin{align} \label{eaoo458}
F^{(\M +1)*}(\EE^{\M+1}) \subseteq F^{(\M)*}(\EE^{\M}) \subseteq \cdots \subseteq F^{(1)*}(\EE^1)  \subseteq F^{(0)*} (\EE^0) = E.
\end{align}
For each integer $l$ with $0 \leq l \leq m +1$,
we shall write
\begin{align}
\mr{Res}(\nabla)^l := \mr{Im} \left( F^{(l) * }(\EE^l) \xrightarrow{\tau^{\downarrow \uparrow (l)}_{(\EE, \nabla)}} \EE \xrightarrow{\mr{quotient}} \mr{Res}(\nabla) \right) \left( \subseteq \mr{Res}(\nabla) \right).
\end{align}
The natural surjection   $F^{(l)*}(\EE^l) \migisurj \mr{Res}(\nabla)^l$ induces  an isomorphism of $R$-modules
\begin{align} \label{Ek32}
F^{(l)*}(\EE^l) / F^{(\M+1)*}(\EE^{\M+1}) \isom \mr{Res}(\nabla)^l.
\end{align}

\SSP
\bpr \label{P021}
The collection $\{ \mr{Res}(\nabla)^l \}_{0 \leq l \leq \M+1}$ forms a decreasing filtration of $\mr{Res}(\nabla)$ satisfying  that $\mr{Res}(\nabla)^0 = \mr{Res}(\nabla)$, $\mr{Res}(\nabla)^{\M+1} = 0$, and 
\begin{align}
\mr{Res}(\nabla)^l / \mr{Res}(\nabla)^{l+1} \cong F^{(l)*}(\mr{Res}(\overline{\nabla}^l))
\end{align}
for every $l = 0, \cdots, \M$.
In particular, the following equality holds:
\begin{align}
\mr{length}_{R}(\mr{Res}(\nabla)) =  \sum_{l=0}^\M p^l \cdot \mr{length}_{R} (\mr{Res}(\overline{\nabla}^l).
\end{align}
\epr
\begin{proof}
Let us prove the former assertion.
It is clear  that $\mr{Res}(\nabla)^0 = \mr{Res}(\nabla)$ and  $\mr{Res}(\nabla)^{\M+1} = 0$.
Next, 
let us take $l \in \{ 0,  \cdots, \M \}$.
The short exact sequence of $R^{(l)}$-modules
\begin{align} \label{E112}
0 \longmigi F^{(1)*}_{l} (E^{l+1}) \longmigi  \EE^l  \longmigi  \mr{Res}(\overline{\nabla}^l)\longmigi 0
\end{align}
induces, via application of the functor  $F^{(l)*}(-)$,   
a short exact sequence
\begin{align} \label{E1g12}
0 \longmigi F^{(l+1)*}(\EE^{l+1}) \left( = F^{(l)*}(F^{(1)*}_{l} (\EE^{l+1}))\right) \longmigi  F^{(l)*}(\EE_l)  \longmigi  F^{(l)*}(\mr{Res}(\overline{\nabla}^l))\longmigi 0.
\end{align}
Hence, the assertion follows from this sequence together with (\ref{Ek32}).

The latter assertion follows from the former one because 
\begin{align}
\mr{length}_{R}(\mr{Res}(\nabla)) &=  \sum_{l=0}^{\M} \mr{length}_R(\mr{Res}(\nabla)^l /\mr{Res}(\nabla)^{l+1}) \\
&= \sum_{l=0}^{\M} \mr{length}_R(F^{(l)*}(\mr{Res}(\overline{\nabla}^l))) \notag \\
&  = \sum_{l=0}^\M p^l \cdot \mr{length}_{R} (\mr{Res}(\overline{\nabla}^l)). \notag
\end{align}
This completes the proof of this proposition.
\end{proof}
\SSP

\begin{exa} \label{Eqqs}
Let $\EE$ be a free $R^{(\M+1)}$-module of finite rank.
Then, the  natural morphism  $\EE \migi (F^{(\M+1)*}(\EE))^{\M+1}$  with respect to the $\breve{\mcD}_R^{(\M+1)}$-module structure $\nabla^\mr{can}_{(R, \breve{\partial}_{\langle \bullet \rangle}), \EE}$ (cf. (\ref{dsskw})) is an isomorphism.
This implies  that  the residue $\mr{Res}(\nabla^\mr{can}_{(R, \breve{\partial}_{\langle \bullet \rangle}), \EE})$ of $\nabla^\mr{can}_{(R, \breve{\partial}_{\langle \bullet \rangle}), \EE}$ vanishes.
\end{exa}

\LSP
\subsection{Dormant $\breve{\mcD}_R^{(\M)}$-modules of rank $1$} \label{SS19}

Let us take an element $a$ of $\mbZ/p^{\M+1} \mbZ$.
 Denote by $\widetilde{a}$ the integer defined as  the unique  lifting of $a$ via the natural surjection $\mbZ \migisurj \mbZ /p^\M \mbZ$ satisfying  $0 \leq \widetilde{a} < p^{\M+1}$.
Then, there exists a unique $\breve{\mcD}_{R}^{(\M)}$-module structure  
\begin{align} \label{eoa78}
\nabla_a := \{ \nabla_{a, \langle j \rangle} \}_{j \in \mbZ_{\geq 0}}
\end{align}
 on $R$ satisfying $\nabla_{a, \langle j \rangle}  (t^n) = q_j !  \cdot \binom{n- \widetilde{a}}{j} \cdot t^n$ (where $\binom{n- \widetilde{a}}{j} := \frac{(n- \widetilde{a}) \cdots (n- \widetilde{a} -j +1)}{j!}$) for every $j, n \in \mbZ_{\geq 0}$.
 In particular,   we have  $\nabla_0  = \breve{\partial}_{\langle \bullet \rangle}$.
The $\breve{\mcD}_{R}^{(\M)}$-module $(R, \nabla_a)$ is isomorphic to 
the unique extension of  
the trivial $\breve{\mcD}_{R}^{(\M)}$-module $(R, \breve{\partial}_{\langle \bullet \rangle})$   to   $t^{-\widetilde{a}}  \cdot R \left( \subseteq K \right)$ (cf. Remark \ref{roafos}).
It follows  immediately that $\nabla_a$ has vanishing $p^{\M+1}$-curvature.

\SSP
\bpr \label{P0045}
\begin{itemize}
\item[(i)]
Let us express the integer $\widetilde{a}$ as  $\widetilde{a} = \sum_{l=0}^{\M}p^l \cdot \widetilde{a}_l$, where $0 \leq \widetilde{a}_l < p$ ($l=0, \cdots, \M$).
Then, the following equalities hold: 
\begin{align}
\mr{length}_R (\mr{Res}(\nabla_a)) = \widetilde{a}, \hspace{3mm} 
\mr{length}_R (\mr{Res}(\nabla_a)^l/\mr{Res}(\nabla_a)^{l+1}) = p^l \cdot  \widetilde{a}_l
\end{align}
($l = 0, \cdots, \M$).
\item[(ii)]
Let $a, b \in \mbZ/ p^{\M +1} \mbZ$.
Then, we have 
\begin{align} \label{E1123}
\mr{Hom}
((R, \nabla_a), (R, \nabla_b)) = \begin{cases} \{ \mu_s \, | \, s \in R^{(\M+1)}\} & \text{if $a =b$}, \\  \{ 0 \} & \text{if $a \neq b$}, \end{cases}
\end{align}
where $\mu_s$ (for each $s \in R^{(\M+1)}$) denotes the endomorphism of $R$ given by multiplication by $s$.
\item[(iii)]
For each $a, b \in \mbZ/p^{\M+1}\mbZ$, we have 
\begin{align} \label{Eww3}
(R, \nabla_a) \otimes (R, \nabla_b) \cong (R, \nabla_{a+b}),  \hspace{3mm} (R^\vee, \nabla_a^\vee) \cong (R, \nabla_{-a}).
\end{align}
\item[(iv)]
Let $(\EE, \nabla)$ be a dormant $\breve{\mcD}_R^{(\M)}$-module with $\EE = R$.
Then,  it is isomorphic to $(R, \nabla_a)$ for a unique element  $a$ of $\mbZ/p^{\M+1}\mbZ$.
\end{itemize}
\epr
\begin{proof}
Assertion (i) follows from (\ref{Ek32}) together with the equality
$\bigcap_{j=0}^{p^l} \mr{Ker} (\nabla_{a, \langle j \rangle}) = t^{\sum_{j =0}^{l}p^{j} \cdot \widetilde{a}_j} \cdot R^{l+1} \left(\subseteq R \right)$, which can be verified by induction on $l$.
Assertions (ii) and (iii) follow immediately from the definition of $\nabla_a$.

Here, we shall prove  assertion (iv).
Let us write $\widetilde{a} := \mr{length}_R (\mr{Res}(\nabla))$ and write $a$ for the image  of $\widetilde{a}$ via the quotient  $\mbZ \migisurj \mbZ/p^{\M+1}\mbZ$.
From  Proposition \ref{P021} and ~\cite[\S\,2, Proposition 2.8]{O2}, we have 
\begin{align} \label{eapfk}
\widetilde{a} = \sum_{l=0}^{\M} p^l \cdot \mr{length}_R (\mr{Res}(\overline{\nabla}^l)) < \sum_{l=0}^{\M}p^l \cdot (p-1) \leq p^{\M+1}.
\end{align}
Since $\tau_{(\EE, \nabla)}^{\downarrow \uparrow (\M+1)} : F^{(\M+1)*}(\EE^{\M+1}) \migi \EE$ is injective and $\EE$ is a free $R$-module  of rank $1$, the $R^{(\M+1)}$-module $\EE^{\M+1}$ may be identified with $R^{(\M+1)}$.
This identification gives an identification $F^{(\M+1)*}(\EE^{\M+1}) = \left(F^{(\M+1)*}(R^{(\M+1)})  =  \right) R$, by which 
the $\breve{\mcD}_R^{(\M+1)}$-module structure  $\nabla_{(R, \breve{\partial}_{\langle \bullet \rangle}), E^{\M+1}}^\mr{can}$ (cf. (\ref{dsskw})) corresponds to  
the trivial  one $\breve{\partial}_{\langle \bullet \rangle}$.
Hence, since $\tau_{(\EE, \nabla)}^{\downarrow \uparrow (\M+1)}$
is $\breve{\mcD}_R^{(\M+1)}$-linear, 
 $\nabla$ may be identified, via $\tau_{(\EE, \nabla)}^{\downarrow \uparrow (\M+1)}$,  with a unique $\breve{\mcD}_R^{(\M+1)}$-module structure on $t^{-\widetilde{a}} \cdot R \left(\subseteq K \right)$ extending  $\breve{\partial}_{\langle \bullet \rangle}$.
 This implies  $\nabla = \nabla_a$ by (\ref{eapfk}),  which completes the proof.
\end{proof}
\SSP

\LSP
\subsection{Exponent of a dormant $\breve{\mcD}_R^{(\M)}$-module} \label{SS070}

Denote by $k$ the residue field of the discrete valuation ring $R$.
Since $k$ is  perfect and $R/(t) = k$, the $t$-adic completion $\widehat{R}$ of $R$ is naturally isomorphic to $k[\![t]\!]$, i.e.,   the ring of formal power series with coefficients in $k$ (cf. ~\cite[Chap.\,I, \S\,4, Theorem 2]{Ser}).
Now, let  $(\EE, \nabla)$ be 
a dormant  $\breve{\mcD}_{R}^{(\M)}$-module such that the $R$-module $\EE$ is  free. 
The $t$-adic completion of $(\EE, \nabla)$ defines a $\breve{\mcD}_{k[\![t]\!]}^{(\M)}$-module $(\widehat{\EE}, \widehat{\nabla})$.
According to ~\cite[Chap.\,1, \S\,1.1, Proposition 1.1.12]{Kin1}, 
there exists an isomorphism  of $\breve{\mcD}_{k[\![t]\!]}^{(\M)}$-modules
\begin{align} \label{L090}
\xi : (\widehat{\EE}, \widehat{\nabla}) \isom \bigoplus_{i=1}^n (k[\![t]\!], \widehat{\nabla}_{d_i}),
\end{align}
where $d_1, \cdots, d_n \in \mbZ/p^{\M+1} \mbZ$ and each $\widehat{\nabla}_{d_i}$ ($i=1, \cdots, n$) denotes
the $\breve{\mcD}_{k[\![t]\!]}^{(\M)}$-module structure on $k[\![t]\!]$ defined as 
 the $t$-adic completion of $\nabla_{d_i}$.
It follows from Proposition \ref{P0045}, (ii),  that 
the multiset $[d_1, \cdots, d_n]$
  depends only on the isomorphism class of   $(\EE, \nabla)$.
(For the definition and various descriptions  concerning a {\it multiset}, we refer the reader  to  ~\cite{SIYS}.)

\SSP
\bde  \label{D033}
In  the above situation, 
the multiset $[d_1, \cdots, d_n]$ is called the {\bf exponent} of $(\EE, \nabla)$.
\ede
\SSP

\begin{exa} \label{Ex989}
Let $(\EE, \nabla)$ be as above and $(L, \nabla_L)$
a dormant $\breve{\mcD}_R^{(\M)}$-module with  $L \cong R$.
 According to Proposition \ref{P0045}, (iv),  $(L, \nabla_L)$ is isomorphic to $(R, \nabla_a)$ for some $a \in \mbZ /p^{\M+1}\mbZ$.
Then, the tensor product $(L \otimes \EE,  \nabla_L \otimes \nabla)$
forms a dormant $\breve{\mcD}_{R}^{(\M)}$-module whose exponent is  $[d_1 + a, d_2 + a, \cdots, d_n + a]$.
\end{exa}
\SSP

\begin{rema} \label{eao78}
For each $\breve{\mcD}_{R}^{(0)}$-module $(\EE', \nabla')$, 
the {\bf monodromy (operator)} of $(\EE', \nabla')$ is 
the element $\mu_{(\EE', \nabla')}$ of $\mr{End}_k (k \otimes_R \EE')$ naturally  induced  by $\nabla'_{\langle 1 \rangle}$ via reduction modulo $(t)$ (cf. ~\cite[Chap.\,6, \S\,1.6, Definition 1.46]{Wak8}).
If $\EE' = R$, we have
 \begin{align} \label{eqgj}
 \mu_{(\EE', \nabla')} \equiv - \mr{length}_R(\mr{Res}(\nabla')) \ \text{mod} \  p.
 \end{align}
Now, suppose that the exponent of a dormant $\breve{\mcD}_R^{(\M)}$-module  $(\EE, \nabla)$ as above  is  $[d_1, \cdots, d_n]$.
For each $i =1, \cdots, n$, we express the integer $\widetilde{d}_i$ (cf. \S\,\ref{SS19}) as $\widetilde{d}_i = \sum_{l=0}^\M p^l \cdot \widetilde{d}_{il}$, where $0 \leq \widetilde{d}_{il}  <  p$.
Then, it follows from (\ref{eqgj}) together with Propositions \ref{P021} and \ref{P0045}, (i), that,  for each $l =0, \cdots, \M$,  the monodromy $\mu_{(\EE^l, \overline{\nabla}^l)}$ of the $\breve{\mcD}_{R^{(l)}}^{(0)}$-module $(\EE^l, \overline{\nabla}^l)$
(cf. (\ref{fmxo})) is diagonalized and  conjugate to the diagonal matrix with diagonal entries $-\widetilde{d}_{1l}, \cdots, -\widetilde{d}_{nl}$ mod $p$.
\end{rema}
\SSP

\bpr \label{Ecsk}
Let $(\EE, \nabla)$ be a dormant $\breve{\mcD}_R^{(\M)}$-module such that the $R$-module  $\EE$ is free.
Then, the following three conditions (a)-(c) are equivalent to each other:
\begin{itemize}
\item[(a)]
The residue $\mr{Res}(\nabla)$ of $\nabla$ vanishes;
\item[(b)]
The exponent of $(\EE, \nabla)$ coincides with $[0, 0, \cdots, 0]$;
\item[(c)]
$(\EE, \nabla)$ comes from a  $\mcD_R^{(\N-1)}$-module via $\eta_R$ (cf. (\ref{ee89d})), meaning that
there exists a  $\mcD_R^{(\N-1)}$-module structure 
$\nabla' \left( = \{ \nabla'_{\langle j \rangle} \}_{j}\right)$ with
$\nabla = \{ t^j \cdot \nabla'_{\langle j \rangle} \}_{j}$.
\end{itemize}

\epr
\begin{proof}
The equivalence (a) $\Longleftrightarrow$ (b) follows from  Proposition \ref{P0045}, (i), and  the existence of the decomposition (\ref{L090}).
The equivalence  (a) $\Longleftrightarrow$ (c) follows from the equivalence of categories asserted in Corollary \ref{eoap09} and the comments in Remark \ref{roeeafos} and  Example \ref{Eqqs} (together with the fact that $\tau_{(\EE, \nabla)}^{\downarrow \uparrow (\M+1)}$ is compatible with the respective $\mcD_R^{(\M+1)}$-module structures, i.e, $\nabla^\mr{can}_{(R, \breve{\partial}_{\langle \bullet \rangle})}$ and $\nabla$).
\end{proof}

\LSP
\subsection{Relationship with the existence of an $\M$-cyclic vector} \label{SS75}

The following assertion helps us to  understand the exponents of  dormant pinned $\breve{\mcD}_R^{(\M)}$-modules.

\SSP
\bpr \label{C001}
Let $(\EE, \nabla)$ be
a dormant $\breve{\mcD}_{R}^{(\M)}$-module such that the $R$-module $\EE$ is free and of rank $n >0$.
Let $[d_1, \cdots, d_n]$ be the exponent of $(\EE, \nabla)$.
Then, the  following two  conditions (a), (b) are equivalent to each other:
\begin{itemize}
\item[(a)]
$(\EE, \nabla)$ admits an $\M$-cyclic vector;
\item[(b)]
The elements $d_1, \cdots, d_n$ of $\mbZ/p^{\M+1} \mbZ$ are mutually distinct.
\end{itemize}

\epr
\begin{proof}
Denote by $(\widehat{\EE}, \widehat{\nabla})$ the $t$-adic completion of $(\EE, \nabla)$, and
 choose an isomorphism $\xi : (\widehat{\EE}, \widehat{\nabla}) \isom \bigoplus_{i=1}^n  (k[\![t]\!], \widehat{\nabla}_{d_i})$ as in  (\ref{L090}).
Let us take an element $v$ of $\EE$, which  
determines an element $\widehat{v}$ of $\widehat{\EE}$ via the natural morphism $\EE \migi \widehat{\EE}$.
Write
    $(u_i)_{i=1}^n  := \xi (\widehat{v})   \in k[\![t]\!]^{\oplus n}$.
 Also, for each $i= 1, \cdots, n$, 
 we  set 
$u_i  = \sum_{s =0}^{\infty} u_{i, s} \cdot t^s$
(where $u_{i,s} \in k$).
Then,   $\widehat{\nabla}_{d_i, \langle j \rangle} (u_i)$ is  expressed as
\begin{equation} \label{10001}
\widehat{\nabla}_{d_i, \langle j \rangle} (u_i) = \sum_{s =0}^{\infty} q_j ! \cdot \binom{s- \widetilde{d}_i}{j}
\cdot u_{i, s} \cdot t^s.
\end{equation}
The isomorphism  $\xi$ preserves the $\breve{\mcD}_{k[\![t]\!]}^{(\M)}$-action, 
so the equality 
$\xi (\widehat{\nabla}_{\langle j \rangle}(\widehat{v} )) = (\widehat{\nabla}_{d_i, \langle j \rangle} (u_i))_{i=1}^n$
holds for every $j = 0, \cdots, n-1$.
Hence, $\{ \widehat{\nabla}_{\langle j \rangle}(\widehat{v} )\}_{0 \leq j \leq n-1}$ forms a basis of $\widehat{\EE}$
 if and only if
the collection $\{ (q_j! \cdot \binom{-\widetilde{d}_i}{j} \cdot u_{i, 0})_{i=1}^n \}_{0 \leq j \leq n-1} \left(= \{ (\nabla_{d_i, \langle j \rangle}(u_i)|_{t=0})_{i=1}^n\}_{0 \leq j \leq n-1} \right)$ forms a basis of $k^{\oplus n}$.
One may verify that this condition is  equivalent to the condition that $u_{i, 0} \in k^\times$ (or equivalently, $u_i \in k[\![t]\!]^\times$) for every $i=1, \cdots, n$ and  $d_1, \cdots, d_n$ are mutually distinct.
This shows the implication  (a) $\Rightarrow$ (b).

Conversely,  suppose that $d_1, \cdots, d_n$ are mutually distinct.
Then, since $k \otimes_R \EE = k \otimes \widehat{\EE} = k^{\oplus n}$,
we can take an element $v$ of $\EE$ such that the induced element $\widehat{v}$ of $\widehat{\EE}$ satisfies  $\xi (\widehat{v}) \in (k[\![t]\!]^\times)^{\oplus n}$.
According to  the above discussion,  $\widehat{v}$ defines an $\M$-cyclic vector of $(\widehat{\EE},  \widehat{\nabla})$.
By the faithful flatness of  the natural homomorphism $R \migi \left(\widehat{R}= \right) k[\![t]\!]$, we see that $v$ forms an $\M$-cyclic vector of $(\EE, \nabla)$.
This implies (b) $\Rightarrow$ (a), which completes  the proof of this proposition.
\end{proof}

\SSP
\bco \label{CC001}
The  exponent of the dormant $\breve{\mcD}_{R}^{(\M)}$-module  $(\breve{P}_R, \nabla_{\breve{P}_R})$ (cf. (\ref{E4567})) coincides with  $[0, 1, \cdots, p^{\M+1}-1]$.

\eco
\begin{proof}
Recall that  $(\breve{P}_R, \nabla_{\breve{P}_R})$  has an $\M$-cyclic vector (by the discussion in \S,\ref{SSd45}) and the free $R$-module $\breve{P}_R$ is of rank $p^{\M+1}$. 
Hence, the assertion follows  from the above proposition. 
\end{proof}
\SSP

Let $\delta := \{ d_1, \cdots, d_n \}$ (where $d_i \neq d_{i'}$ if $i \neq i'$) be 
a subset of $\mbZ/p^{\M+1} \mbZ$ whose cardinality equals $n$.
We shall set
\begin{align} \label{eww299}
\delta^\BB :=  \{ d^\BB_1, \cdots, d_{p^{\M+1}-n}^\BB\}
\end{align}
to  be the subset  of $\mbZ/p^{\M+1} \mbZ$ with  $\delta \sqcup  \{ -d^\BB_1, \cdots,  - d_{p^{\M+1}-n}^\BB\}  = \ \mbZ/p^{\M+1} \mbZ$.
It can immediately be seen  that $\delta^{\BB \BB} = \delta$.

If a dormant $\breve{\mcD}_R^{(\M)}$-module admits an $\M$-cyclic vector, then it follows from Proposition \ref{C001} that    its exponent may be regarded as a subset of $\mbZ/p^{\M+1} \mbZ$.
In particular, for each subset $\delta$ of $\mbZ/p^{\M+1} \mbZ$,  it makes sense to speak of a dormant pinned $\breve{\mcD}_R^{(\M)}$-module {\it of exponent $\delta$}. 
For each subset $\delta := \{ d_1, \cdots, d_n \}$ of $\mbZ/p^{\M+1} \mbZ$ whose cardinality equals $n$, we shall denote by 
\begin{align}
\mr{Mod}(\breve{\mcD}_R^{(\M)})^\circledast_{n, \delta}
\end{align}
 the subset of
 $\mr{Mod}(\breve{\mcD}_R^{(\M)})^\circledast_{n}$ consisting of dormant pinned  $\breve{\mcD}_{R}^{(\M)}$-modules of rank $n$ and exponent $\delta$.

\bpr \label{P0h23}
The bijection 
 $[\breve{\Dual}_R]_n : \mr{Mod}(\breve{\mcD}_R^{(\M)})^\circledast_{n} \isom  \mr{Mod}(\breve{\mcD}_R^{(\M)})^\circledast_{p^{\M+1} - n}$ (cf. (\ref{ew234}))
 is restricted to a bijection of sets
\begin{align} \label{eapgagopg}
[\breve{\Dual}_R]_{n, \delta} : \mr{Mod}(\breve{\mcD}_R^{(\M)})^\circledast_{n, \delta}  \isom \mr{Mod}(\breve{\mcD}_R^{(\M)})^\circledast_{p^{\M+1}-n, \delta^\BB}.
\end{align}
\epr
\begin{proof}
Let us take a dormant pinned $\breve{\mcD}^{(\M)}_R$-module $(\EE, \nabla, v)$ classified by $\mr{Mod} (\breve{\mcD}_R^{(\M)})^\circledast_{n, \delta}$.
Denote by $(\widehat{\EE}, \widehat{\nabla})$
  the $t$-adic completion of $(\EE, \nabla)$.
Note that the $t$-adic completion of $(\breve{P}_{R}, \nabla_{\breve{P}_{R}})$ may be identified with $(\breve{P}_{k[\![t]\!]}, \nabla_{\breve{P}_{k[\![t]\!]}})$.
Hence, 
the $t$-adic completion $\widehat{\nu}_{(\EE, \nabla, v)}$  of $\nu_{(\EE, \nabla, v)}$ (cf. (\ref{Eppop})) 
specifies a morphism $(\breve{P}_{k[\![t]\!]}, \nabla_{\breve{P}_{k[\![t]\!]}}) \migi (\widehat{\EE}, \widehat{\nabla})$.
We shall fix an isomorphism  
of $\breve{\mcD}_{k[\![t]\!]}^{(\M)}$-modules $\xi : (\widehat{E}, \widehat{\nabla}) \isom \bigoplus_{i=1}^n (k[\![t]\!], \widehat{\nabla}_{d_i})$ as in (\ref{L090}).
On the other hand, 
according to Corollary \ref{CC001}, there exists an isomorphism $\xi_{P} : (\breve{P}_{k[\![t]\!]}, \nabla_{\breve{P}_{k[\![t]\!]}}) \isom \bigoplus_{d \in \mbZ/p^{\M+1}\mbZ}^{} (k[\![t]\!], \widehat{\nabla}_{d})$.
In particular, we obtain a surjective  morphism
\begin{align}
\xi \circ \widehat{\nu}_{(\EE, \nabla, v)}\circ \xi_P^{-1} : \bigoplus_{d \in \mbZ/p^{\M+1}\mbZ}^{} (k[\![t]\!], \widehat{\nabla}_{d}) \migi \bigoplus_{i=1}^n (k[\![t]\!], \widehat{\nabla}_{d_i}).
\end{align}
From  Proposition \ref{P0045}, (ii), 
 we see that the kernel 
 of this morphism is isomorphic to $\bigoplus_{d \notin \delta}(k[\![t]\!], \widehat{\nabla}_d)$.
 Hence,  the  dual of this dormant $\breve{\mcD}_{k[\![t]\!]}^{(\M)}$-module, i.e., the $t$-adic completion of $(\EE^\BB, \nabla^\BB)$, is isomorphic to $\bigoplus_{d  \in \delta^\BB} (k[\![t]\!], \widehat{\nabla}_d)$ (cf. Proposition \ref{P0045}, (iii)).
 This means that the exponent of $(\EE^\BB, \nabla^\BB)$ coincides with $\delta^\BB$, which   completes the proof.
\end{proof}

\vspace{10mm}
\section{Dormant $\mr{PGL}_n$-opers of higher level} \label{LkJf} \SSP

In the second half of the present paper, 
 we apply  the results on higher-level differential modules proved so far to  discuss the corresponding objects defined on an algebraic  curve in characteristic $p >0$.
The notion of a $\mr{PGL}_n$-oper of level $\N >0$ will be defined in terms of  the ring of differential operators of  level $\N-1$ (cf. Definition \ref{D29}).
Also, we introduce the radius of a $\mr{PGL}_n$-oper of level $\N$ by using 
the local description  at each marked point of the underlying pointed curve (cf. Definition \ref{epaddd}).

\LSP
\subsection{Logarithmic differential operators} \label{SS07}

In the rest of the present paper, let us fix  a positive 
integer $\N$,  a nonnegative integer $r$, and 
   an algebraically closed field $k$  of characteristic $p$.
Also,  let us fix an $r$-pointed (possibly nonproper) smooth  curve  
 \begin{align} \label{ea34j}
 \msX := (f :X \migi \mr{Spec}(k), \{ \sigma_i \}_{1 \leq i \leq r})
 \end{align}
   over $k$, i.e., a  smooth curve $X$ over $k$
  together with $r$ marked points $\{\sigma_i \}_{1 \leq i \leq r} \left(\subseteq X(k)\right)$.
The  divisor on $X$ defined as the union of the marked points $\sigma_i$ determines 
a log structure on $X$;
we shall denote the resulting log scheme by $X^\mr{log}$.
Since $X^\mr{log}/k$ is log smooth over $k$, 
 the sheaf of logarithmic $1$-forms $\Omega_{X^\mr{log}/k}$ of $X^\mr{log}/k$, as well as its dual $\mcT_{X^\mr{log}/k} := \Omega_{X^\mr{log}/k}^\vee$,  is a line bundle.

Next, write $F_k$ (resp., $F_X$) for the absolute Frobenius endomorphism of $\mr{Spec}(k)$ (resp., $X$).
We shall denote by  $X^{(\N)}$ 
 the base-change $k \times_{F_k^\N, k} X$ of $X$ by  the $\N$-th iterate $F_k^{\N}$ of $F_k$; we will refer to it as  the {\bf $\N$-th Frobenius twist} of $X$ over $k$.
Also, the morphism $F^{(\N)}_{X/k} \left(:= (f, F_X^{\N}) \right) : X \migi X^{(\N)}$  is called the {\bf $\N$-th relative Frobenius morphism} of $X$ over $k$.

Denote by $\mcD^{(\N-1)}_{X^\mr{log}}$ (cf. ~\cite[Chap.\,2, \S\,2.3, Definition 2.3.1]{Mon})
  the ring of logarithmic differential operators on $X^\mr{log}/k$ (equipped with the trivial $(\N-1)$-PD structures) of level $\N-1$.
  For each integer $j$, we shall write  $\mcD_{X^\mr{log}, < j }^{(\N-1)}$  for the $\mcO_X$-submodule of $\mcD_{X^\mr{log}}^{(\N-1)}$ consisting of logarithmic differential operators of order $<  j$.
  If $\N'$ is  a positive integer with $\N' \geq \N$,
  then there exists a canonical morphism $\mcD_{X^\mr{log}}^{(\N-1)} \migi \mcD_{X^\mr{log}}^{(\N' -1)}$ (cf. ~\cite[Chap.\,2, \S\,2.5.1]{Mon}).

A {\bf (left) $\mcD_{X^\mr{log}}^{(\N-1)}$-module structure} on  an $\mcO_X$-module $\mcE$ is 
a left $\mcD_{X^\mr{log}}^{(\N-1)}$-action $\nabla : \mcD_{X^\mr{log}}^{(\N-1)} \migi \mcE nd_k (\mcE)$ on $\mcE$  extending its $\mcO_X$-module structure.  
An $\mcO_X$-module equipped with a $\mcD_{X^\mr{log}}^{(\N-1)}$-module structure is called a {\bf (left) $\mcD_{X^\mr{log}}^{(\N-1)}$-module}.
Given a $\mcD_{X^\mr{log}}^{(\N-1)}$-module $(\mcE, \nabla)$, we shall write 
$\mcE^\nabla$
for the subsheaf of $\mcE$ on which $\mcD_{X^\mr{log}, +}^{(\N-1)}$ acts  as zero, where
$\mcD_{X^\mr{log}, +}^{(\N-1)}$ denotes the kernel of the canonical projection $\mcD_{X^\mr{log}}^{(\N-1)} \migisurj \mcO_X$.
The sheaf $\mcE^\nabla$ may be regarded as an $\mcO_{X^{(\N)}}$-module via the underlying homeomorphism of $F^{(\N)}_{X/k}$.

Recall that giving a $\mcD_{X^\mr{log}}^{(0)}$-module structure on an $\mcO_X$-module $\mcE$ is equivalent to giving a logarithmic connection on $\mcE$, i.e., 
a $k$-linear morphism $\mcE \migi \Omega_{X^\mr{log}/k} \otimes \mcE$ satisfying 
the Leibniz rule.
For each $\mcD_{X^\mr{log}}^{(\N-1)}$-module structure $\nabla$ on an $\mcO_X$-module $\mcE$, we shall write 
\begin{align} \label{KKclo}
\overline{\nabla} : \mcE \migi \Omega_{X^\mr{log}/k} \otimes \mcE
\end{align}
 for the logarithmic connection on $\mcE$
corresponding to the $\mcD_{X^\mr{log}}^{(0)}$-module structure induced from $\nabla$ via the canonical morphism $\mcD_{X^\mr{log}}^{(0)} \migi \mcD_{X^\mr{log}}^{(\N-1)}$.

 Denote by ${^p}\psi_{X^\mr{log}}$  the $p^\N$-curvature map $\mcT_{X^\mr{log}/k}^{\otimes p^\N} \migi \mcD_{X^\mr{log}}^{(\N-1)}$ defined in ~\cite[\S\,3.2, Definition 3.10]{Ohk}.
Given a $\mcD_{X^\mr{log}}^{(\N-1)}$-module $(\mcE, \nabla)$,
we shall  set
\begin{align}
{^p}\psi_{(\mcE, \nabla)} := \nabla \circ {^p}\psi_{X^\mr{log}}  : \mcT^{\otimes p^\N}_{X^\mr{log}/k} \migi \mcE nd_k (\mcE),
\end{align} 
which will be called the {\bf $p^\N$-curvature} of $(\mcE, \nabla)$.
The $p^1$-curvature of  a $\mcD_{X^\mr{log}}^{(0)}$-module  is essentially the same as the $p$-curvature  of the corresponding logarithmic connection (cf. ~\cite[\S\,1.2]{Og}).

\SSP
\bde \label{epapDD}
Let $(\mcE, \nabla)$ be a $\mcD_{X^\mr{log}}^{(\N-1)}$-module.
Then, we shall say that $(\mcE, \nabla)$ is {\bf dormant} if ${^p}\psi_{(\mcE, \nabla)} = 0$.
\ede

\SSP
\begin{rema} \label{eRfi}
Let us review  a result in  the case where 
$r =0$,  or equivalently, the log structure of $X^\mr{log}$ is trivial.
 Then, our definition of $p^\N$-curvature coincides with the $p$-$(\N-1)$-curvature in the sense of ~\cite[\S\,3, Definition 3.1.1]{LeQu}.
 For each $\mcO_{X^{(\N)}}$-module $\mcE$,
 there exists a canonical $\mcD_{X}^{(\N-1)}$-module structure $\nabla^\mr{can}_\mcE$ on
 $F^{(\N)*}_{X/k}(\mcE)$ with vanishing $p^\N$-curvature.
According to ~\cite[\S\,3, Corollary 3.2.4]{LeQu}, the assignments $(\mcE, \nabla) \mapsto \mcE^\nabla$ and $\mcE \mapsto (F^{(\N)*}_{X/k}(\mcE, \nabla_\mcE^\mr{can}))$ determine an equivalence of categories
  \begin{align} \label{Efjk2}
\begin{pmatrix}
\text{the category of } \\
\text{dormant $\mcD_{X}^{(\N-1)}$-modules} 
\end{pmatrix}
\isom \begin{pmatrix}
\text{ the category of } \\
\text{$\mcO_{X^{(\N)}}$-modules} \\
\end{pmatrix}.
\end{align}
The ring-theoretic counterpart of this equivalence was already verified in  Corollary \ref{eoap09}.
\end{rema}

\SSP
Here, we shall consider  the local description of  $\mcD_{X^\mr{log}}^{(\N-1)}$-modules.
 Let $x$ be a $k$-rational point of $X$ and fix a local function $t$ on $X$ defining $x$.
 Suppose that  $x$ lies in $X \setminus \bigcup_{i=1}^n \{ \sigma_i \}$ (resp., $\bigcup_{i=1}^n \{\sigma_i\}$).
 Then, it follows from ~\cite[\S\,2.2, Proposition 2.2.4]{PBer1} (resp.,  ~\cite[Chap.\,2, Lemma 2.3.3]{Mon}) that (the sections of) the restriction of $\mcD_{X^\mr{log}}^{(\N-1)}$ to  
 $D_x := \mr{Spec}(\mcO_{X, x})$
  may be identified with $\mcD_{\mcO_{X, x}}^{(\N-1)}$ (resp., $\breve{\mcD}_{\mcO_{X, x}}^{(\N-1)}$) defined in (\ref{wqo9}) for  $R = \mcO_{X, x}$. 
 Each  $\mcD_{X^\mr{log}}^{(\N-1)}$-module $(\mcE, \nabla)$ induces, via restriction to $D_x$,  a $\mcD_{\mcO_{X, x}}^{(\N-1)}$-module  (resp., $\breve{\mcD}_{\mcO_{X, x}}^{(\N-1)}$-module) $(\mcE, \nabla) |_{D_x}$.
  According to ~\cite[\S\,3.2, Proposition 3.11]{Ohk},
  the $p^\N$-curvature of $(\mcE, \nabla) |_{D_x}$ may be regarded as
  the restriction to $D_x$ of the $p^\N$-curvature of $(\mcE, \nabla)$.
  In particular, if
  $(\mcE, \nabla)$ is dormant, then the restriction $(\mcE, \nabla) |_{D_x}$ is dormant.

\SSP
\bde \label{Ddfow}
Let $(\mcE, \nabla)$ be a dormant $\mcD_{X^\mr{log}}^{(\N-1)}$-module such that   $\mcE$  is a vector bundle on $X$ of rank $n>0$.
Suppose that $r >0$. 
Then,  for each  $ i \in \{ 1, \cdots, n \}$,
the {\bf (local) exponent} of  $(\mcE, \nabla)$ (or, of $\nabla$) at the marked point $\sigma_i$
is defined as  the exponent of  $(\mcE, \nabla) |_{D_{\sigma_i}}$ (cf. Definition \ref{D033}). 
\ede
\SSP

The following assertion will be applied in the proof of Proposition \ref{Pioe}.

\SSP
\bpr \label{Lkg09}
Let $l$ be a positive integer with $p \nmid l$, $\mcN$ a line bundle on $X$, and 
$\nabla_{\mcN^{\otimes l}}$  a $\mcD_{X^\mr{log}}^{(\N-1)}$-module structure on the $l$-th tensor product $\mcN^{\otimes l}$ of $\mcN$ with vanishing $p^\N$-curvature.
Then, there exists a unique  $\mcD_{X^\mr{log}}^{(\N -1)}$-module structure $\nabla_{\mcN}$ on $\mcN$ with vanishing $p^\N$-curvature
whose $l$-th tensor product 
   $\nabla_\mcN^{\otimes l}$ coincides with $\nabla_{\mcN^{\otimes l}}$.
   \epr
\begin{proof}
For each $i =1, \cdots, r$,
denote by $e_i$ the exponent of $\nabla_{\mcN^{\otimes l}}$ at $\sigma_i$.
Write 
$\widetilde{e}'_i$ for the unique integer with $0 \leq \widetilde{e}'_i < p^\N$ and 
$\widetilde{e}'_i \equiv e_i/l$ mod $p^\N$.
According to the discussion in Remark \ref{roafos},
the trivial $\mcD_{X^\mr{log}}^{(\N-1)}$-module structure on $\mcO_X$ extends uniquely to a $\mcD_{X^\mr{log}}^{(\N-1)}$-module structure $\nabla_+$ on 
$\mcO_{X}(\sum_{i=1}^r\widetilde{e}'_i \cdot \sigma_i) \left(\supseteq \mcO_X \right)$ with vanishing $p^\N$-curvature.
The exponent of 
$\nabla_+^{\otimes l} \otimes \nabla_{\mcN^{\otimes l}}$
at $\sigma_i$ is $\overline{l \cdot \widetilde{e}'_i} - e_i = l \cdot (e_i/l) -e_i =  0$.
Hence, 
by Proposition \ref{Ecsk} together with  the equivalence of categories (\ref{Efjk2}),
there exists a line bundle $\mcM$ on $X^{(\N)}$ with  
$(F^{(\N)*}(\mcM), \nabla_{\mcM}^\mr{can}) \cong (\mcN'^{\otimes l}, \nabla_{+}^{\otimes l} \otimes \nabla_{\mcN^{\otimes l}})$, where $\mcN' := \mcO_{X}(\sum_{i=1}^r\widetilde{e}'_i \cdot \sigma_i)\otimes \mcN$.
If  $\mcM'$ denotes  the line bundle on $X$ corresponding to $\mcM$ via base-change $X^{(\N)} \isom X$ by $F^\N_k$, then we have
$\mcM'^{\otimes p^{\N}} \cong \mcN'^{\otimes l}$.
Here,
let us take a pair of  integers $(a, b)$ with $a \cdot p^\N + b \cdot  l =1$.
Then, 
\begin{align} \label{eapc2}
\mcN' =  \mcN'^{\otimes (a \cdot p^\N + b \cdot  l)} \cong \mcN'^{\otimes ap^\N} \otimes \mcM'^{\otimes b p^\N}  = (\mcN'^{\otimes a} \otimes \mcM'^{\otimes b})^{\otimes p^\N}.
\end{align}
Let us define   $\mcL$  to be  the line bundle on $X^{(\N)}$ corresponding to $\mcN'^{\otimes a} \otimes \mcM'^{\otimes b}$ via base-change by  $F_k^{\N}$.
By (\ref{eapc2}),   we see that $F^{(\N)*}(\mcL) \cong \mcN'$, i.e.,  there exists an isomorphism $F^{(\N)*}(\mcL) \otimes \mcO_X (-\sum_i \widetilde{e}'_i \cdot  \sigma_i) \isom  \mcN$.
The line bundle $\mcN$ is equipped with the $\mcD_{X^\mr{log}}^{(\N-1)}$-module structure $\nabla_\mcN$ corresponding to $\nabla_\mcL^\mr{can} \otimes \nabla_+^\vee$ via this isomorphism.
The $l$-th tensor product  $\nabla_{\mcN}^{\otimes l}$ of $\nabla_\mcN$ has vanishing $p^\N$-curvature and coincides with $\nabla_{\mcN^{\otimes l}}$ by its construction.
This completes the proof of this proposition.
\end{proof}

\LSP
\subsection{$\mr{GL}_n$-opers and $\mr{PGL}_n$-opers  of level $\N$} \label{SS042d}

Let us fix a positive integer $n$.
We shall define the notion of a 
$\mr{GL}_n$-oper  of level $\N$, as follows.
(Note that a $\mr{GL}_n$-oper of level $1$ is the same as a $\mr{GL}_n$-oper in the classical sense.)

\SSP
\bde \label{D23}
\begin{itemize}
\item[(i)]
Let us consider a collection of data
\begin{align}
\msF^\heartsuit := (\mcF, \nabla, \{ \mcF^j \}_{0 \leq j \leq n}),
\end{align}
where
\begin{itemize}
\item
$\mcF$ is a vector bundle on $X$ of rank $n$;
\item 
$\nabla$ is a $\mcD_{X^{\mr{log}}}^{(\N-1)}$-module structure on $\mcF$;
\item
$\{ \mcF^j \}_{0 \leq j \leq n}$ is an $n$-step decreasing filtration
\begin{align}
0 = \mcF^n \subseteq \mcF^{n-1} \subseteq \cdots \subseteq \mcF^0 = \mcF
\end{align}
on $\mcF$ consisting of subbundles.
\end{itemize}
Then, we say that $\msF^\heartsuit$ is a {\bf $\mr{GL}_n$-oper of level $\N$} on $\msX$
if it satisfies the following three conditions:
\begin{itemize}
\item
For each $j=0, \cdots, n -1$, the subquotient $\mcF^j /\mcF^{j+1}$ is a line bundle;
\item
For each $j=1, \cdots, n-1$, $\overline{\nabla} (\mcF^j)$ (cf. (\ref{KKclo})) is contained in $\Omega_{X^\mr{log}/k} \otimes \mcF^{j-1}$;
\item
For each $j=1, \cdots, n-1$, the well-defined {\it $\mcO_X$-linear} morphism
\begin{align} \label{eraopaz98}
\mr{KS}_{\msF^\heartsuit}^j : \mcF^j /\mcF^{j+1} \migi \Omega_{X^\mr{log}/k} \otimes (\mcF^{j-1}/\mcF^j)
\end{align}
induced by $\overline{\nabla} |_{\mcF^j}  : \mcF^j \migi \Omega_{X^\mr{log}/k} \otimes \mcF^{j-1}$
 is an isomorphism.
\end{itemize}
The notion of an isomorphism between two $\mr{GL}_n$-opers of level $\N$ can be defined in a natural manner (so we omit the details).
\item[(ii)]
Let $\msF^\heartsuit := (\mcF, \nabla, \{ \mcF^j \}_j)$ be a $\mr{GL}_n$-oper of level $\N$.
Then, we shall say that $\mcF^\heartsuit$ is {\bf dormant} if ${^p}\psi_{(\mcF, \nabla)} = 0$.
\end{itemize}
\ede
\SSP

\begin{rema} \label{errc3}
\begin{itemize}
\item[(i)]
A  (dormant) $\mr{GL}_1$-oper of level $\N$ is the same as a (dormant) $\mcD_{X^\mr{log}}^{(\N-1)}$-module $(\mcF, \nabla)$ such that $\mcF$ is a line bundle.
\item[(ii)]
In  the case of $n=2$,
 a $\mr{GL}_2$-oper of level $\N$ on $\msX$ is given as a triple  $(\mcF, \nabla, \mcL)$ consisting of 
 \begin{itemize}
 \item
 a $\mcD_{X^\mr{log}}^{(\N-1)}$-module   $(\mcF, \nabla)$ such that
 $\mcF$ is a rank $2$ vector bundle, and 
 \item
  a line subbundle $\mcL$ of $\mcF$ such that
 the $\mcO_X$-linear  composite
 \begin{align} \label{fai752}
 \mcL \xrightarrow{\mr{inclusion}} \mcF \xrightarrow{\overline{\nabla}} \Omega_{X^\mr{log}/k} \otimes \mcF \xrightarrow{\mr{quotient}} \Omega_{X^\mr{log}/k} \otimes (\mcF/ \mcL)
 \end{align} 
 defines  an isomorphism between line bundles.
 \end{itemize}
 \end{itemize}
\end{rema}
\SSP

The following assertion implies that higher-level differential modules with a cyclic vector may be regarded as ring-theoretic counterparts of $\mr{GL}_n$-opers of  higher level.

\SSP
\bpr \label{Peiri3}
Let $(\mcF, \nabla)$ be a $\mcD_{X^\mr{log}}^{(\N-1)}$-module such that $\mcF$ is a vector bundle of rank $n$.
Also, let $x$ be 
a (possibly generic) point $x$ of $X \setminus \bigcup_{i=1}^r \{\sigma_i\}$ (resp., $\bigcup_{i=1}^r \{\sigma_i\}$).
Write $(\mcF_x, \nabla_x) := (\mcF, \nabla) |_{D_x}$ for the $\mcD_{\mcO_{X, x}}^{(\N-1)}$-module (resp., $\breve{\mcD}_{\mcO_{X, x}}^{(\N-1)}$-module)
obtained as the restriction of $(\mcF, \nabla)$ to $D_x \left(:= \mr{Spec}(\mcO_{X, x}) \right)$.
\begin{itemize}
\item[(i)]
Suppose that there exists an $n$-step decreasing  filtration $\{ \mcF^j \}_{0 \leq j \leq n}$ on $\mcF$ for which the collection $(\mcF, \nabla, \{ \mcF^j \}_{0 \leq j \leq n})$ forms a $\mr{GL}_n$-oper of level $\N$.
Then, each generator of the restriction of  $\mcF^{n-1}$ to $D_x$ defines an $(\N-1)$-cyclic vector of  $(\mcF_x, \nabla_x)$.
\item[(ii)]
Suppose that there exists an $(\N-1)$-cyclic vector of $(\mcF_x, \nabla_x)$.
For each $j=0, \cdots, n$, we shall write $\mcF_x^j$ for the $\mcO_{X, x}$-submodule of $\mcF_x$ generated by  the elements $\nabla_{x, \langle l \rangle}(v)$ for $l \leq n-j-1$.
Then, 
there exists  an open neighborhood $U$ of $x$ satisfying the following condition: the filtration  $\{ \mcF_{x}^j \}_j$ extends to a decreasing filtration $\{ (\mcF|_U)^j \}_{0 \leq j \leq n}$ of $\mcF |_U$ (i.e., $(\mcF|_U)^j |_{D_x} = \mcF_x^j$ for every $j$) for which the collection $(\mcF |_U, \nabla |_U, \{ (\mcF|_U)^j \}_{0 \leq j \leq n})$ forms a $\mr{GL}_n$-oper of level $\N$ on the pointed curve $\msX$ restricted to $U$.
\end{itemize}
\epr
\begin{proof}
The assertions 
follow immediately from the 
definitions of a  $\mr{GL}_n$-oper of level $\N$ and an $(\N-1)$-cyclic vector.
\end{proof}
\SSP

By applying results on differential modules proved in the previous discussion,
we can obtain the following assertion.

\SSP
\bco \label{Peiri2}
\begin{itemize}
\item[(i)]
Let $(\mcF, \nabla)$ be a $\mcD_{X^\mr{log}}^{(\N-1)}$-module such that $\mcF$ is a vector bundle of rank $n \leq p^\N$.
Also, let $x$ be a $k$-rational point of $X$.
Then,  there exists an  open neighborhood  $U$ of $x$ in  $X$ 
 and 
 an $n$-step  decreasing filtration $\{ (\mcF|_U)^j \}_{0 \leq j \leq n}$ on $\mcF |_U$
  such that
 the collection $(\mcF |_U, \nabla |_U, \{ (\mcF|_U)^j \}_{0 \leq j \leq n})$ forms a $\mr{GL}_n$-oper of level $\N$ on the pointed curve  $\msX$ restricted to $U$.
\item[(ii)]
If $n > p^{\N}$, then there are no dormant $\mr{GL}_n$-opers on $\msX$ of level $\N$.
\end{itemize}
\eco
\begin{proof}
Assertion (i) follows from Theorem \ref{P023} and Proposition \ref{Peiri3}, (ii).
Assertion (ii)
follows from Propositions  \ref{er45gj}, (i),  and \ref{Peiri3}, (i).
\end{proof}
\SSP

Next, we shall define an equivalence relation in the set of dormant $\mr{GL}_n$-opers of level $\N$.
Let $\msF^\heartsuit := (\mcF, \nabla, \{ \mcF^j \}_{j})$ be a $\mr{GL}_n$-oper of level $\N$ on $\msX$ and $(\mcN, \nabla_\mcN)$ 
a line bundle  on $X$ equipped with a  $\mcD_{X^\mr{log}}^{(\N-1)}$-module structure.
According to ~\cite[Chap.\,2, Corollary 2.6.1]{Mon},
there exists 
a canonical   $\mcD_{X^\mr{log}}^{(\N-1)}$-module structure $\nabla \otimes \nabla_\mcN$  on  the tensor product $\mcF \otimes \mcN$ naturally arising  from $\nabla$ and $\nabla_\mcN$.
One may verify that the collection
\begin{align}
\msF_{ \otimes (\mcN, \nabla_\mcN)}^\heartsuit := (\mcN \otimes\mcF, \nabla_\mcN \otimes \nabla,  \{ \mcN \otimes \mcF^j \}_{0 \leq j \leq n})
\end{align} 
forms a  $\mr{GL}_n$-oper of level $\N$.
If both $\msF^\heartsuit$ and $(\mcN, \nabla_\mcN)$
are dormant, 
 then  $\msF_{ \otimes \msN}^\heartsuit$ is  dormant.
Now, let us consider  the binary relation ``$\sim$" in the set of dormant $\mr{GL}_n$-opers of level $\N$ on $\msX$ defined by
$\msF^\heartsuit \sim \msF'^\heartsuit$ if and only if 
$\msF^\heartsuit_{\otimes (\mcN, \nabla_\mcN)}\cong \msF'^\heartsuit$ for some 
$(\mcN, \nabla_\mcN)$ as above;
 this relation in fact defines an equivalence relation.
 For each   $\msF^\heartsuit$  as above, we shall write
 \begin{align} \label{ffw78}
 \msF^{\heartsuit \Rightarrow \spadesuit}
 \end{align}
  for the equivalence class represented by  $\msF^\heartsuit$.
 
\SSP
\bde \label{D29}
A {\bf $\mr{PGL}_n$-oper of level $\N$} on $\msX$  is the equivalence class $\msF^\spadesuit \left(= \msF^{\heartsuit \Rightarrow \spadesuit}\right)$ of a $\mr{GL}_n$-oper  $\mcF^\heartsuit$ of level $\N$ on $\msX$.
A $\mr{PGL}_n$-oper of level $\N$ is called {\bf dormant} if it may be   represented by a dormant $\mr{GL}_n$-oper of level $\N$.
\ede
\SSP

We shall denote by
\begin{align}
{^\N}\mcO p_{n,  \msX}^{^\mr{Zzz...}}
\end{align}
the set of dormant $\mr{PGL}_n$-opers of level $\N$ on $\msX$.

\SSP
\begin{rema} \label{RRRw3}
It can immediately be seen that $\sharp ({^\N}\mcO p_{1,  \msX}^{^\mr{Zzz...}}) =1$ (cf. Remark \ref{errc3}, (i)).
Also, according to  
~\cite[Chap.\,II, \S\,2.3, Theorem 2.8]{Mzk2},
the set ${^\N}\mcO p_{2,  \msX}^{^\mr{Zzz...}}$ is nonempty when 
 $\msX$ is general in the moduli space $\mcM_{g, r}$ of $r$-pointed smooth curves of genus $g$ over $k$.
For a general $n$, we know that ${^1}\mcO p_{n,  \msX}^{^\mr{Zzz...}}$ is nonempty when $n$ is sufficiently small  relative to $p$ (cf. ~\cite[Chap.\,3, \S\,3.8, Theorem 3.38]{Wak8}).
On the other hand,  it follows from  Corollary \ref{Peiri2} that ${^\N}\mcO p_{n,  \msX}^{^\mr{Zzz...}} = \emptyset$ if $n > p^\N$.
\end{rema}

\LSP
\subsection{Radius  of  a dormant $\mr{PGL}_n$-oper of level $\N$} \label{SS0f46}

Denote by $\Delta$ the image of the diagonal embedding  $\mbZ/p^\N \mbZ \migiincl (\mbZ/p^\N \mbZ)^{\times n}$.
In particular, by regarding it as a group homomorphism, we obtain the quotient  $(\mbZ/p^\N \mbZ)^{\times n}/\Delta$.
Note that the set $(\mbZ/p^\N \mbZ)^{\times n}$ is equipped with the action of the symmetric group $\mfS_n$ of $n$ letters by permutation;
this action induces a well-defined $\mfS_n$-action on $(\mbZ/p^\N \mbZ)^{\times n}/\Delta$.
Hence, we obtain the sets $\mfS_n \backslash (\mbZ/p^\N \mbZ)^{\times n}$, $\mfS_n \backslash (\mbZ/p^\N \mbZ)^{\times n}/\Delta$, and moreover, obtain the natural projection
\begin{align} \label{ssak}
\pi : \mfS_n \backslash (\mbZ/p^\N \mbZ)^{\times n} \migisurj \mfS_n \backslash (\mbZ/p^\N \mbZ)^{\times n}/\Delta.
\end{align}
Each element of $\mfS_n \backslash (\mbZ/p^\N \mbZ)^{\times n}$ may be regarded as a multiset of $\mbZ/p^\N \mbZ$ whose cardinality equals $n$.

\SSP
\begin{rema} \label{ewgk}
Let us consider the case of $n=2$.
Denote by $(\mbZ/p^\N \mbZ)/\{ \pm 1 \}$ the set of equivalence classes of elements $a \in \mbZ/p^\N \mbZ$, in which 
 $a$ and $-a$ are identified.
 Then, the assignment $a \mapsto [a, -a]$ determines a well-defined bijection
\begin{align} \label{wpaid}
(\mbZ/p^\N \mbZ)/\{ \pm 1 \} \isom \mfS_2 \backslash (\mbZ/p^\N \mbZ)^{\times 2}/\Delta.
\end{align}
By using this bijection, we will identify $\mfS_2 \backslash (\mbZ/p^\N \mbZ)^{\times 2}/\Delta$ with  $(\mbZ/p^\N \mbZ)/\{ \pm 1 \}$ (cf. the discussion in \S\,\ref{S09p}).
\end{rema}
\SSP

Let $\msF^\spadesuit$ be a dormant $\mr{PGL}_n$-oper on $\msX$ of level $\N$, and choose
a dormant $\mr{GL}_n$-oper  $\msF^\heartsuit := (\mcF, \nabla, \{ \mcF^j \}_{j})$ of level $\N$ representing  $\msF^\spadesuit$.
Suppose that $r>0$.
For each $i =1, \cdots, r$, 
denote by $\delta_i$ the exponent of $(\mcF, \nabla)$ at $\sigma_i$.
Let us write
$\rho_{\msF^\spadesuit, i} := \pi (\delta_i) \in \mfS_n \backslash (\mbZ/p^\N \mbZ)^{\times n}/\Delta$.
It follows from the fact mentioned in Example \ref{Ex989} that the  element $\rho_{\msF^\spadesuit, i}$ does not depend on the choice of the representative $\msF^\heartsuit$ of $\msF^\spadesuit$.

\SSP
\bde \label{epaddd}
\begin{itemize}
\item[(i)]
We shall refer to $\rho_{\msF^\spadesuit, i}$ as the {\bf radius} of $\msF^\spadesuit$ at $\sigma_i$.
\item[(ii)]
Let $\vec{\rho} := (\rho_i )_{i=1}^r$ be an element of $(\mfS_n \backslash (\mbZ/p^\N \mbZ)^{\times n}/\Delta)^{\times r}$.
We shall say that $\msF^\spadesuit$ is {\bf of radii $\vec{\rho}$} if $\rho_i = \rho_{\msF^\spadesuit, i}$ for every $i =1, \cdots, r$.
\end{itemize}
\ede
\SSP

For each  $\vec{\rho} \in  (\mfS_n \backslash (\mbZ/p^\N \mbZ)^{\times n}/\Delta)^{\times r}$,
we shall denote by
\begin{align}
{^\N}\mcO p_{n,  \msX, \vec{\rho}}^{^\mr{Zzz...}} 
\end{align}
the subset of 
${^\N}\mcO p_{n,  \msX}^{^\mr{Zzz...}}$
consisting of dormant $\mr{PGL}_n$-opers of level $\N$ 
and  radii $\vec{\rho}$.

\SSP
\begin{rema} \label{peec}
Let us recall the previous study for  $\N=1$.
A $\mr{PGL}_2$-oper of level  $1$ is essentially the same as a {\it torally indigenous bundle} in the sense of
 ~\cite[Chap.\,I, \S\,4, Definition 4.1]{Mzk2}.
Also, the radii of a dormant  torally indigenous bundle on $\msX$ (which  belong to the set $\mbF_p$, as proved in   ~\cite[Chap.\,II, \S\,1, Proposition 1.5]{Mzk2}) is consistent   with the radii of the corresponding $\mr{PGL}_2$-oper via the quotient  $\mbF_p \migisurj \mbF_p/\{\pm1\} \left(= (\mbZ/p\mbZ)/\{\pm1\} \right)$.
According to ~\cite[Chap.\,II, \S\,1, Proposition 1.4]{Mzk2},
${^1}\mcO p_{2,  \msX, (\rho)_{i=1}^r}^{^\mr{Zzz...}}$ is empty  unless 
$\rho_i \in \mbF_p^\times/\{\pm 1 \}$ for every $i=1, \cdots, r$.

Moreover,   for a general $n$, 
 the radii of a dormant $\mr{PGL}_n$-oper introduced above
 coincides with the one  in the sense of ~\cite[Chap.\,2, \S\,2.8, Definition 2.32]{Wak8}
 under the identification of  each element in   $\mfS_n \backslash \mbF_p^{\times n} /\Delta$ with  an $\mbF_p$-rational point in  the adjoint quotient of the Lie algebra $\mfp \mfg \mfl_n$. 
 \end{rema}

\vspace{10mm}
\section{Duality of dormant $\mr{PGL}_n$-opers}\label{Sgh46} \SSP

In this section, we establish a duality between dormant $\mr{PGL}_n$-opers of level $\N$ and dormant $\mr{PGL}_{p^\N-n}$-opers of level $\N$ (cf. Theorem \ref{T5949}).
As a corollary, we will see that there is exactly one isomorphism class of dormant $\mr{PGL}_{p^\N -1}$-oper of level $\N$ (cf. Corollary \ref{C90e8}, (ii)).

We  keep the notation in the previous section.

\LSP
\subsection{Dormant $\mr{GL}_{p^\N}$-opers of level $\N$} \label{SS046}

Let $\mcL$ be a line bundle on $X$.
We equip $\mcD_{X^\mr{log}}^{(\N-1)} \otimes \mcL$ with  the  $\mcD_{X^\mr{log}}^{(\N-1)}$-module structure  given by left multiplication.
We shall write
$\mcP_\mcL$ for the quotient of the left $\mcD_{X^\mr{log}}^{(\N-1)}$-module $\mcD_{X^\mr{log}}^{(\N-1)} \otimes \mcL$ by the $\mcD_{X^\mr{log}}^{(\N-1)}$-submodule generated by the image of ${^p}\psi_{X^\mr{log}} \otimes \mr{id}_\mcL : \mcT_{X^\mr{log}/k}^{\otimes p^\N} \otimes \mcL \migi \mcD_{X^\mr{log}}^{(\N-1)} \otimes \mcL$.
Denote by $\nabla_{\mcP_\mcL}$ the resulting  $\mcD_{X^\mr{log}}^{(\N-1)}$-module structure of $\mcP_{\mcL}$;   by construction,  $(\mcP_\mcL, \nabla_{\mcP_\mcL})$ has  vanishing  $p^\N$-curvature.
Also, for each $j = 0, \cdots, p^\N$, we shall set $\mcP^j_\mcL$ to  be the subbundle of $\mcP_\mcL$ defines as 
\begin{align}
\mcP^j_\mcL  := \mr{Im} \left(\mcD_{X^\mr{log}/k, < p^\N -j }^{(\N-1)} \otimes \mcL\xrightarrow{\mr{inclusion}} \mcD_{X^\mr{log}}^{(\N-1)} \otimes \mcL \xrightarrow{\mr{quotient}} \mcP_\mcL \right).
\end{align}
The collection of data
 \begin{align}
 \msP^\heartsuit_\mcL := (\mcP_\mcL, \nabla_{\mcP_\mcL}, \{ \mcP^j_\mcL \}_{0 \leq j \leq p^\N})
 \end{align}
  forms a dormant 
  $\mr{GL}_{p^\N}$-oper of level $\N$ on $\msX$.
  Indeed, as discussed in \S\,\ref{SSd45},  the restriction of this data to $D_x \left(= \mr{Spec}(\mcO_{X, x}) \right)$ for  each $k$-rational point $x$ of $X \setminus \bigcup_{i=1}^r \{\sigma_i \}$ (resp., $\bigcup_{i=1}^r \{\sigma_i\}$) defines  a dormant pinned $\mcD_{\mcO_{X, x}}^{(\N-1)}$-module (resp., a dormant pinned $\breve{\mcD}_{\mcO_{X, x}}^{(\N-1)}$-module).
  In particular, we obtain  a dormant $\mr{PGL}_{p^\N}$-oper $\msP^{\heartsuit \Rightarrow \spadesuit}_\mcL$ of level $\N$ on $\msX$.

Next, let $(\mcN, \nabla_\mcN)$ be a dormant $\mcD_{X^\mr{log}}^{(\N-1)}$-module such that $\mcN$ is a line bundle.
Since the tensor product $\nabla_\mcN \otimes \nabla$ has vanishing $p^\N$-curvature, 
the composite
\begin{align}
\mcD_{X^\mr{log}}^{(\N-1)} \otimes (\mcN \otimes \mcL) \left(= \mcD_{X^\mr{log}}^{(\N-1)} \otimes (\mcN \otimes \mcP_\mcL^{p^\N-1} ) \right)
&\xrightarrow{\mr{inclusion}}
 \mcD_{X^\mr{log}}^{(\N-1)} \otimes (\mcN \otimes \mcP_\mcL) \\
 & \xrightarrow{\nabla_\mcN \otimes \nabla}
 \mcN \otimes \mcP_\mcL \notag
\end{align}
factors through the  quotient $\mcD_{X^\mr{log}}^{(\N-1)} \otimes (\mcN \otimes \mcL) \migisurj \mcP_{\mcN \otimes \mcL}$.
By considering the local description, we can see  that the  resulting morphism $\mcP_{\mcN \otimes \mcL} \migi \mcN \otimes \mcP_\mcL$
defines  an isomorphism
\begin{align} \label{eaopd}
\msP^\heartsuit_{\mcN \otimes \mcL} \isom (\msP^\heartsuit_{\mcL})_{\otimes (\mcN, \nabla_\mcN)}
\end{align}
 of dormant $\mr{GL}_n$-opers of level $\N$.

\LSP
\subsection{Duality for  dormant $\mr{GL}_n$ opers of  level $\N$} \label{SSggh046}

Next, let $\msF^\heartsuit := (\mcF, \nabla, \{ \mcF^j \}_{j})$ be a dormant $\mr{GL}_n$-oper of level $\N$ on $\msX$ with $\mcF^{n-1} = \mcL$.
The inclusion $\mcL \left(=\mcD_{X^\mr{log}, <1}^{(\N-1)} \otimes \mcL  \right) \migiincl \mcF$ extends uniquely  to a $\mcD_{X^\mr{log}}^{(\N-1)}$-linear morphism $\mcD_{X^\mr{log}}^{(\N-1)} \otimes \mcL \migi \mcF$.
Since $(\mcF, \nabla)$ has vanishing $p^\N$-curvature, 
this morphism factors through the quotient $\mcD_{X^{\mr{log}}}^{(\N-1)} \otimes \mcL \migisurj \mcP_\mcL$.
Thus, we obtain a morphism of $\mcD_{X^{\mr{log}}}^{(\N-1)}$-modules
\begin{align}
\nu_{\msF^\heartsuit} : (\mcP_\mcL, \nabla_{\mcP_\mcL}) \migi (\mcF, \nabla).
\end{align}
By Proposition \ref{Peiri3}, (i), the  restriction of this morphism  to $D_x$ for each point $x$ of $X \setminus \bigcup_{i=1}^r \{\sigma_i\}$ (resp., $\bigcup_{i=1}^r \{\sigma_i\}$) may be regarded as a morphism of 
 pinned $\mcD_{\mcO_{X, x}}^{(\N-1)}$-modules (resp., pinned $\breve{\mcD}_{\mcO_{X, x}}^{(\N-1)}$-modules). 
Hence, it follows from  Proposition \ref{edaP}, (i), that  $\nu_{\msF^\heartsuit}$ is verified to be surjective.

\SSP
\begin{exa} \label{eapa9}
Let us consider the dual $(\mcP_\mcL^\vee, \nabla_{\mcP_\mcL}^\vee)$ of $(\mcP_\mcL, \nabla_{\mcP_\mcL})$.
For each $j = 0, \cdots, p^\N$,
we shall write $\mcP_\mcL^{\vee j}$ for the $\mcO_X$-submodule of $\mcP_\mcL^\vee$ defined to be the image of the natural injection  $(\mcP_\mcL/\mcP_\mcL^{p^\N -j})^\vee \migiincl  \mcP^\vee_\mcL$.
Thus, we obtain a collection of data
\begin{align}
\msP^{\heartsuit \vee}_\mcL := (\mcP_\mcL^\vee, \nabla_{\mcP_\mcL}^\vee, \{ \mcP_\mcL^{\vee j} \}_{0 \leq j \leq p^\N}).
\end{align}
For each $k$-rational point $x$ of $X$, the restriction of $(\mcP_\mcL^\vee, \nabla_{\mcP_\mcL}^\vee)$ to $D_x$  together with a generator of $\mcP_\mcL^{\vee p^\N -1} |_{D_x}$ may be regarded as the data (\ref{fao332}).
It follows that  $\msP^{\heartsuit \vee}_\mcL$ forms a dormant $\mr{GL}_{p^\N}$-oper of level $\N$.
Also, the morphism
\begin{align} \label{eapa98z}
\nu_{\msP_\mcL^{\heartsuit \vee}} : (\mcP_\mcL, \nabla_{\mcP_\mcL}) \migi (\mcP_\mcL^\vee, \nabla_{\mcP_\mcL}^\vee)
\end{align}
is compatible, via restriction to $D_x$,  with the morphism  (\ref{fao333}).
It follows that  $\nu_{\msF^{\heartsuit \vee}}$ defines an isomorphism $\msP^\heartsuit_\mcL \isom \msP^{\heartsuit \vee}_\mcL$ of $\mr{GL}_{p^\N}$-opers of level $\N$.
\end{exa}
\SSP

Let us write $\mcF^\BB := \mr{Ker}(\nu_{\msF^\heartsuit})^\vee$.
Since $\nu_{\msF^\heartsuit}$ preserves the $\mcD_{X^\mr{log}}^{(\N-1)}$-action,
$\nabla_{\mcP_\mcL}$ is restricted to a $\mcD_{X^\mr{log}}^{(\N-1)}$-module structure $\nabla_{\mr{Ker}(\nu_{\msF^\heartsuit})}$ on 
  the kernel $\mr{Ker}(\nu_{\msF^\heartsuit})$; it induces 
  a $\mcD_{X^\mr{log}}^{(\N-1)}$-module structure on $\mcF^\BB$, which we denote by $\nabla^\BB$.
 For each $j = 0, \cdots, p^\N -n$, let us define $\gamma_j$ to  be the composite 
\begin{align}
\gamma_j : \mr{Ker}(\nu_{\msF^\heartsuit}) \xrightarrow{\mr{inclusion}} \mcP_\mcL \migisurj \mcP_\mcL/\mcP^{p^\N -n - j }_\mcL.
\end{align}
By letting  $\mcF^{\BB j}:= \mr{Im}(\gamma^\vee_j) \left( \subseteq \mcF^\BB \right)$, we  obtain  
a collection of data
\begin{align}
\msF^{\heartsuit \BB} := (\mcF^\BB, \nabla^\BB, \{ \mcF^{\BB j} \}_{0 \leq j \leq p^\N-n}).
\end{align}

\SSP
\bpr \label{P038}
Let us keep the above notation.
\begin{itemize}
\item[(i)]
$\msF^{\heartsuit \BB}$
 forms a dormant $\mr{GL}_{p^\N -n}$-oper of level $\N$ on $\msX$.
 Moreover,  there exists a canonical isomorphism 
$\msF^{\heartsuit} \isom   \msF^{\heartsuit \BB\BB}$ of $\mr{GL}_n$-opers of level $\N$.
\item[(ii)]
Let $(\mcN, \nabla_\mcN)$ be a dormant  $\mcD_{X^\mr{log}}^{(\N-1)}$-module  such that $\mcN$ is a line bundle.
Then, there exists a canonical isomorphism 
\begin{align}
(\msF^{\heartsuit \BB})_{\otimes (\mcN^\vee, \nabla_\mcN^\vee)} \cong (\msF^{\heartsuit}_{\otimes (\mcN, \nabla_\mcN)})^\BB
\end{align}
of $\mr{GL}_{p^\N -n}$-opers of level $\N$.
\end{itemize}
\epr
\begin{proof}
First, let us consider assertion (i).
The formation of $\msF^{\heartsuit \BB}$ is compatible with 
that of $(\EE^\BB, \nabla^\BB, v^\BB)$ (cf. (\ref{E456d})) via restriction to $D_x$ for every point $x$ of $X$.
This implies that $\msF^{\heartsuit \BB}$ forms a dormant $\mr{GL}_{p^\N -n}$-oper of level $\N$ (cf. Proposition \ref{Peiri3}).
Also, let us consider  the following diagram:
\begin{align} \label{Edd7}
\vcenter{\xymatrix@C=56pt@R=36pt{
0 \ar[r]  & \mr{Ker}(\nu_{\msF^{\heartsuit \BB}}) \ar[r]^-{\mr{inclusion}}  & \mcP_\mcL \ar[d]_-{\nu_{\msP_\mcL^{\heartsuit \BB}}} \ar[r]^-{\nu_{\msF^{\heartsuit \BB}}} & \mcF^\BB \ar[d]_-{\wr}^-{\mr{id}} \ar[r] & 0
 \\
 0  \ar[r] & \mcF^\vee \ar[r]_-{\nu_{\msF^\heartsuit}^\vee} & \mcP_\mcL^\vee \ar[r]_-{\mr{quotient}} &  \mcF^\BB  \ar[r] & 0.
  }}
\end{align}
The right-hand square is commutative because its restriction to $D_x$ for every point $x$ of $X$ is commutative, as observed in the proof of Theorem \ref{T4}, (i).
Hence, it induces a morphism $\mr{Ker}(\nu_{\msF^{\heartsuit \BB}}) \isom \mcF^\vee$.
By considering the local description again, one verifies that the dual of this isomorphism specifies an isomorphism $\msF^\heartsuit \isom \msF^{\heartsuit \BB \BB}$ of $\mr{GL}_{n}$-opers of level $\N$.
This completes the proof of assertion (i).

Next, we shall prove  assertion (ii).
Consider the following diagram:
\begin{align} \label{Eddd7}
\vcenter{\xymatrix@C=50pt@R=36pt{
0 \ar[r]  & \mcN^\vee \otimes \mcF^\vee \ar[r]^-{\mr{id}_{\mcN^\vee} \otimes \nu_{\msF^\heartsuit}^\vee} \ar[d]_-{\wr}  & \mcN^\vee \otimes \mcP_\mcL^\vee \ar[d]_-{\wr} \ar[r]^-{\mr{quotient}} & \mcN^\vee  \otimes \mr{Ker}(\nu_{\msF^\heartsuit})^\vee  \ar[r] & 0
 \\
 0  \ar[r] & (\mcN \otimes\mcF)^\vee \ar[r]_-{\nu_{\msF^\heartsuit}^\vee} & \mcP_{\mcN \otimes \mcL}^\vee \ar[r]_-{\mr{quotient}} &  \mr{Ker}(\nu_{\msF^\heartsuit_{\otimes (\mcN, \nabla_\mcN)}})^\vee  \ar[r] & 0,
  }}
\end{align}
where the middle vertical arrow denotes the dual of (\ref{eaopd}) and the left-hand vertical arrow denotes the canonical isomorphism.
Since the left-hand square diagram is commutative, this diagram induces an isomorphism
$\mcN^\vee  \otimes \mr{Ker}(\nu_{\msF^\heartsuit})^\vee \isom  \mr{Ker}(\nu_{\msF^\heartsuit_{\otimes (\mcN, \nabla_\mcN)}})^\vee$.
This isomorphism specifies an isomorphism $(\msF^{\heartsuit \BB})_{\otimes (\mcN^\vee, \nabla_\mcN^\vee)} \isom (\msF^{\heartsuit}_{\otimes (\mcN, \nabla_\mcN)})^\BB$ of $\mr{GL}_{p^\N-n}$-opers of level $\N$.
This completes the proof of assertion (ii).
\end{proof}
\SSP
\begin{rema} \label{R0ddd90}
In this remark,
we shall  examine  the determinant of
a $\mr{GL}_n$-oper.
Let $\msF^\heartsuit := (\mcF, \nabla, \{ \mcF^j \}_j)$ be as above.
Since $\mcD_{X^\mr{log}, < j+1}^{(\N-1)}/\mcD_{X^\mr{log}, < j}^{(\N-1)} \cong \mcT_{X^\mr{log}/k}^{\otimes j}$ ($j=0,1,2, \cdots$),
we obtain the composite of canonical isomorphisms 
\begin{align} \label{P09df}
\mr{det} (\mcP_\mcL) &\isom   \bigotimes_{j=0}^{p^\N-1} \mcP_\mcL^j / \mcP_\mcL^{j+1} \\
&\isom \bigotimes_{j=0}^{p^\N -1} \mcT_{X^\mr{log}/k}^{\otimes (p^\N - j -1)} \otimes \mcL  \notag\\
&\isom  \mcT_{X^\mr{log}/k}^{\otimes p^\N (p^\N-1)/2} \otimes \mcL^{\otimes p^\N} \notag \\
& \isom F^{(\N)*}_{X/k} (\mcN_\mcL),\notag
\end{align}
where $\mcN_\mcL$ denotes the line bundle on $X^{(\N)}$ corresponding to $\mcT_{X^\mr{log}/k}^{\otimes (p^\N-1)/2}\otimes \mcL$ via  base-change $X^{(\N)} \isom X$ by $F_k^\N$.
Similarly,
there exists an isomorphism 
\begin{align} \label{ew34}
\mr{det}(\mcF) \isom \mcT_{X^\mr{log}/k}^{\otimes n(n-1)/2} \otimes \mcL^{\otimes n}.
\end{align}
The $\mcD_{X^\mr{log}}^{(\N-1)}$-module structure on the determinant bundle $\mr{det}(\mcP_\mcL)$ induced by $\nabla_{\mcP_\mcL}$ corresponds to $\nabla_{\mcN_\mcL}^\mr{can}$ via (\ref{P09df}).
Hence, 
 the determinant of $\nabla^\BB$ corresponds to $(\nabla^\mr{can}_{\mcN_\mcL})^\vee \otimes \mr{det} (\nabla)$ via the following composite of natural isomorphisms:
\begin{align} \label{fsp98395}
\mr{det}(\mcF^\BB) &\isom  \mr{det}(\mr{Ker}(\nu_{\msF^\heartsuit}))^\vee \isom \mr{det}(\mcP_\mcL)^\vee  \otimes \mr{det}(\mcF)  \isom  F^{(\N)*}_{X/k}(\mcN_\mcL)^\vee \otimes \mr{det} (\mcF). 
\end{align}
\end{rema}
\SSP

\LSP
\subsection{Duality for  dormant $\mr{PGL}_n$ opers of  level $\N$} \label{Sfamh046}

By applying Proposition \ref{P038}, (i), 
we obtain a  bijective correspondence
 \begin{align} \label{Edee2}
\begin{pmatrix}
\text{the set of isomorphism classes of} \\
\text{dormant $\mr{GL}_n$-opers  of level $\N$ on $\msX$} 
\end{pmatrix}
\isom \begin{pmatrix}
\text{the set of isomorphism classes of} \\
\text{dormant $\mr{GL}_{p^\N -n}$-opers  of level $\N$ on $\msX$} \\
\end{pmatrix}. \hspace{-7mm}
\end{align}
Moreover, this correspondence  and  Proposition \ref{P038}, (ii), together   imply  the following assertion, which is a part of Theorem \ref{efs89}.

\SSP
\bt \label{T5949}
\begin{itemize}
\item[(i)]
The assignment $\msF^{\heartsuit \Rightarrow \spadesuit} \mapsto \msF^{\heartsuit \BB \Rightarrow \spadesuit}$ defines a well-defined  bijection of sets
\begin{align}
\Dual_{n, \msX} : {^\N}\mcO p^{^\mr{Zzz...}}_{n, \msX} \isom {^\N}\mcO p_{p^\N -n, \msX}^{^\mr{Zzz...}}
\end{align}
satisfying   $\Dual_{p^\N -n, \msX} \circ \Dual_{n, \msX} = \mr{id}$.
\item[(ii)]
Suppose further that $r>0$.
Let  $\vec{\rho} := (\rho_i)_{i=1}^r$ be an element of $(\mfS_n \backslash (\mbZ/p^\N \mbZ)^{\times n}/ \Delta)^{\times r}$.
Then, $\Dual_{n, \msX}$ is restricted to a bijection
\begin{align}
\Dual_{n, \msX, \vec{\rho}} : {^\N}\mcO p^{^\mr{Zzz...}}_{n, \msX, \vec{\rho}} \isom {^\N}\mcO p_{p^\M -n, \msX, \vec{\rho}^{\,\BB}}^{^\mr{Zzz...}},
\end{align}
where $\vec{\rho}^{\,\BB} := (\rho_i^\BB)_{i=1}^r$
(cf. (\ref{eww299})).
\end{itemize}
\et
\begin{proof}
Assertions follow from Propositions \ref{P0h23} and  \ref{P038}.
\end{proof}
\SSP

Moreover, the duality theorem established  above implies 
the following assertion, which is  the remaining portion of Theorem \ref{efs89}; note that this assertion generalizes  results proved by Y. Hoshi (cf. ~\cite[Theorem A]{Hos1}) and  the author (cf.  ~\cite[\S\,4, Corollary 4.3.3]{Wak3}).

\SSP
\bco \label{C90e8}
\begin{itemize}
\item[(i)]
Let $\mcL$ be a line bundle on $X$ and $\nabla_{\circledcirc}$ a dormant $\mcD_{X^\mr{log}}^{(\N-1)}$-module structure on $\mcT_{X^\mr{log}/k}^{\otimes n (n-1)/2} \otimes \mcL^{\otimes n}$.
Then,
there exists exactly  one isomorphism class of dormant $\mr{GL}_{p^\N -1}$-oper $\msF^\heartsuit := (\mcF, \nabla, \{ \mcF^j \}_j)$ of level $\N$  on $\msX$ such that  $\mcF^{n-1} = \mcL$ and $\mr{det}(\nabla) = \nabla_{\circledcirc}$ under the identification $\mr{det}(\mcF) = \mcT_{X^\mr{log}/k}^{\otimes n (n-1)/2} \otimes \mcL^{\otimes n}$ given by  (\ref{ew34}).
\item[(ii)]
There exists exactly one  isomorphism class of dormant $\mr{PGL}_{p^\N-1}$-oper of level $\N$, i.e., the following equality holds:
\begin{align}
\sharp (\mcO p_{p^\N-1, \msX}^{^\mr{Zzz...}}) = 1.
\end{align}
Moreover, if $r > 0$, then the
radius of 
 the unique  dormant $\mr{PGL}_{p^\N -1}$-oper of level $\N$  at any marked point  coincides with $\pi ([1,2, \cdots, p^\N-1])$ (cf. (\ref{ssak})).
\end{itemize}
\eco
\begin{proof}
Assertion (i) follows from the observation  that the desired $\mr{GL}_{p^\N -1}$-oper is the unique
one corresponding, via (\ref{Edee2}), to the level-$\N$ dormant 
$\mr{GL}_1$-oper 
\begin{align}
(F^{(\N)*}_{X/k}(\mcN_\mcL)^\vee \otimes \mcT_{X^\mr{log}/k}^{\otimes n (n-1)/2}\otimes \mcL^{\otimes n}, (\nabla^\mr{can}_{\mcN_\mcL})^\vee \otimes \nabla_\circledcirc)
\end{align}
 (cf. Remarks \ref{errc3}, (i), and  \ref{R0ddd90}).
Assertion (ii) follows from the bijection $\Dual_{1, \msX}$ asserted in Theorem \ref{T5949}, (i), and the equalities  $\sharp ({^\N}\mcO p^{^\mr{Zzz...}}_{1, \msX}) = \sharp ({^\N}\mcO p^{^\mr{Zzz...}}_{1, \msX, (\pi ([0]), \cdots, \pi ([0]))}) = 1$.
\end{proof}

\vspace{10mm}
\section{Tamely ramified coverings and dormant $\mr{PGL}_2$-opers}  \label{S09p}
\SSP

Recall (cf.  ~\cite{Mzk2}, ~\cite{O2}, ~\cite{Oss2}) that
certain tamely ramified coverings with ramification indices $<p$ between two copies of the projective line can be  described in terms of dormant $\mr{PGL}_2$-opers (of level $1$).
That description is the  starting point of the enumerative geometry of dormant opers because it allows us to translate   dormant $\mr{PGL}_2$-opers on a $3$-pointed projective line into simple combinatorial  data.
In this section, 
the situation is generalized to the case of higher level in order to deal with tamely ramified covering having large ramification indices.
Theorem \ref{Re345} will be proved at the end of this section.

\LSP
\subsection{Dormant $\mr{PGL}_2$-opers arising from   tamely ramified coverings} \label{SSd4h2}

Denote by 
$\mbP$ the projective line   over $k$, i.e., $\mbP := \mr{Proj} (k[x_1, x_2])$.
Let $\msX := (X, \{ \sigma_i \}_{i=1}^r)$ be as before, and take an $r$-tuple of integers $\vec{\lambda} := (\lambda_1, \cdots, \lambda_r)$
with $0< \lambda_i < p^\N$ ($i=1, \cdots, r$).
We shall  write $\vec{\rho} := (\rho_1, \cdots, \rho_r)$, where  $\rho_i := \frac{1}{2} \cdot \overline{\lambda}_i \in (\mbZ / p^\N \mbZ)/\{ \pm 1 \}$.
Under the identification $(\mbZ / p^\N \mbZ)/\{ \pm 1 \} = \mfS_2 \backslash (\mbZ/p^\N \mbZ)^{\times 2}/\Delta$ defined in  (\ref{wpaid})),
$\vec{\rho}$ may be regarded as an element of $(\mfS_2 \backslash (\mbZ/p^\N \mbZ)^{\times 2}/\Delta)^{\times r}$.
We shall denote by
\begin{align}
\mr{Cov}^\mr{tame}_{\msX, \vec{\lambda}}
\end{align}
the set of equivalence classes of  finite, separable, and  tamely ramified coverings  $\phi : X \migi \mbP$ that are ramified at 
$\sigma_i$ with index $\lambda_i$
 and \'{e}tale elsewhere.
 Here, the equivalence relation is defined in such a way that two coverings $\phi_1, \phi_2 : X \migi \mbP^1$ are equivalent if there exists an element $h \in \mr{PGL}_2 (k) \left(= \mr{Aut}_k (\mbP) \right)$ with $\phi_2 = h \circ \phi_1$.
 For each $\phi$ as above, we shall denote by $[\phi]$ the element of $\mr{Cov}^\mr{tame}_{\msX, \vec{\lambda}}$ (i.e., the equivalence class) represented  by $\phi$.

In what follows, let us construct a map of sets  $\mr{Cov}^\mr{tame}_{\msX, \vec{\lambda}} \migi {^\N}\mcO p_{2, \msX, \vec{\rho}}^{^\mr{Zzz...}}$.
Let us take 
 an element $[\phi]$  of $\mr{Cov}^\mr{tame}_{\msX, \vec{\lambda}}$, and choose 
a tamely  ramified covering $\phi$ representing $[\phi]$.
Let  $\sigma'_1, \cdots, \sigma'_{r'}$ ($0 < r' \leq r$) be the mutually distinct points of $\mbP$ such that
$\bigcup_{i=1}^r \{ \phi (\sigma_i)\}  = \{ \sigma'_j \}_{j=1}^{r'}$.
In particular, $\msP' := (\mbP, \{ \sigma'_j \}_{j=1}^{r'})$ defines an $r'$-pointed genus-$0$ curve; we denote the  induced  log curve by  $\mbP^{\mr{log}'}$. 
Since $\phi$ is tamely ramified,  the morphism  $\phi$ extends to a {\it log \'{e}tale} morphism $\phi^\mr{log} : X^\mr{log} \migi \mbP^{\mr{log}'}$.
Write $\mcL := \mcO_{\mbP}(-1) \otimes \mcO_{\mbP}(\sum_{j=1}^{r'}\sigma_j)$, and write $\tau_0$ for  the $\mcO_\mbP$-linear injection $\mcO_{\mbP} (-1) \migiincl \mcO_{\mbP}^{\oplus 2}$ given by $w \mapsto (w x_1, w x_2)$ for each local section $w \in \mcO_\mbP (-1)$.
Also, let $\mcF$ be a rank $2$ vector bundle on $\mbP$ which makes the following square diagram cocartesian:
\begin{align} \label{E0233f}
\vcenter{\xymatrix@C=46pt@R=36pt{
\mcO_{\mbP}(-1)\ar[r]^-{\tau_0} \ar[d]_-{\mr{inclusion}} & \mcO_\mbP^{\oplus 2}\ar[d] 
\\
\mcL \ar[r] & \mcF.
}}
\end{align}
The trivial $\mcD^{(\N-1)}_{\mbP^{\mr{log}'}}$-module structure  on $\mcO_{\mbP}^{\oplus 2}$ extends uniquely to 
a $\mcD^{(\N-1)}_{\mbP^{\mr{log}'}}$-module structure  $\nabla_\mcF$ on $\mcF$.
It follows from the various definitions involved that the composite
\begin{align}
\mcL \xrightarrow{\mr{inclusion}} \mcF \xrightarrow{\overline{\nabla}_\mcF} \Omega_{\mbP^{\mr{log}'}/k} \otimes \mcF \migisurj \Omega_{\mbP^{\mr{log}'}/k} \otimes (\mcF/ \mcL)
\end{align}
is $\mcO_\mbP$-linear and injective.
Moreover, since $\mr{deg}(\mcL) = \mr{det}(\Omega_{\mbP^{\mr{log}'}/k} \otimes (\mcF/ \mcL)) \left(= r' -1 \right)$, this morphism turns out to be an isomorphism.
This means that the triple $(\mcF, \nabla_\mcF, \mcL)$ forms a dormant $\mr{GL}_2$-oper of level $\N$ on $\msP'$ (cf. Remark \ref{errc3}, (ii)).
Hence, the pull-back of this data via the log \'{e}tale morphism $\phi^\mr{log}$ defines  a dormant $\mr{GL}_2$-oper
\begin{align}
\msF^{\heartsuit}_{\phi} := \phi^{\mr{log}*} (\mcF, \nabla, \mcL)
\end{align}
of level $\N$ on $\msX$.

\SSP
\bpr \label{Peoid}
The level-$\N$ dormant $\mr{PGL}_2$-oper
 $\msF^{\heartsuit \Rightarrow \spadesuit}_{\phi}$
 on $\msX$ represented  by $\msF_\phi^\heartsuit$ is of radii $\vec{\rho}$.
\epr
\begin{proof}
The problem  is the computation of the radii of $\msF^\spadesuit_{\phi}$.
Let us take $i \in \{ 1, \cdots, r \}$, and  choose  $j \in \{1, \cdots, r' \}$ such that  $\sigma'_j = \phi (\sigma_i)$.
Also, choose a local function $t$ on $\mbP$ defining $\sigma'_j$.
This local function allows us to identify 
the formal neighborhood $\widehat{D}_{\sigma'_{j}}$ of $\sigma'_j$ in $\mbP$ with $\mr{Spec} (k[\![t]\!])$.
Since the ramification index of $\phi$ at $\sigma_i$ is $\lambda_i$,
the formal neighborhood $\widehat{D}_{\sigma_i}$ of $\sigma_i$ in $X$ may be identified with $\mr{Spec}(k[\![t^{1/\lambda_i}]\!])$ and the restriction of $\phi$ to $\widehat{D}_{\sigma_i}$ may be identified with the morphism $\mr{Spec}(k[\![t^{1/\lambda_i}]\!]) \migi \mr{Spec}(k[\![t]\!])$ induced by the natural inclusion $k[\![t]\!] \migiincl k[\![t^{1/\lambda_i}]\!]$.
The $\breve{\mcD}_{k[\![t]\!]}^{(\N-1)}$-module corresponding to the restriction of $(\mcF, \nabla)$ to $\widehat{D}_{\sigma'_j}$ is isomorphic to $(k[\![t]\!], \widehat{\nabla}_0) \oplus (t^{-1} \cdot k[\![t]\!], \widehat{\nabla}_0)$.
It follows that the pull-back of $(\mcF, \nabla)$ to $\phi^\mr{log}$ restricted to $\widehat{D}_{\sigma_i}$ is isomorphic to   $(k[\![t^{1/\lambda_i} ]\!], \widehat{\nabla}_0) \oplus ( (t^{1/\lambda_i})^{-\lambda_i} \cdot k[\![t^{1/\lambda_i}]\!], \widehat{\nabla}_0)$ (which is also isomorphic to  $(k[\![s]\!], \widehat{\nabla}_0) \oplus (k[\![s]\!], \widehat{\nabla}_{\overline{\lambda}_i})$ by putting $s := t^{1/\lambda_i}$).
Hence, 
the exponent of $\msF_\phi^{\heartsuit}$  at $\sigma_i$ coincides with $[0, \overline{\lambda}_i]$, which implies 
$\rho_{\msF_\phi^{\heartsuit \Rightarrow \spadesuit}, i} = \rho_i$.
This completes the proof of this proposition.
\end{proof}
\SSP

Since $\tau_0$ is invariant under  pull-back by automorphisms of $\mbP$,
the isomorphism class of $\msF^{\heartsuit}_{\phi}$ does not depend on the choice of the representative $\phi$ of $[\phi]$.
Hence,  the above proposition implies that  the assignment $[\phi] \mapsto \msF_\phi^\spadesuit$ gives 
a well-defined map of sets
\begin{align} \label{eavbk}
\Upsilon_{\msX, \vec{\lambda}} : \mr{Cov}^\mr{tame}_{\msX, \vec{\lambda}}\migi {^\N}\mcO p_{2,  \msX, \vec{\rho}}^{^\mr{Zzz...}}.
\end{align}

\LSP
\subsection{Tamely ramified endomorphisms of a $3$-pointed projective line} \label{SS042e}

Denote  by 
$[0]$, $[1]$, and $[\infty]$ the $k$-rational points of   $\mbP$ determined by the values $0$, $1$, and $\infty$ respectively.
After ordering  the three points $ [0], [1], [\infty]$, we obtain 
 a unique (up to isomorphism) $3$-pointed proper smooth  curve 
\begin{equation} \label{1051}
\msP := (\mbP/k, \{ [0], [1], [\infty] \})
\end{equation}
 of genus $0$ over $k$.
 In particular, we obtain a log curve $\mbP^\mr{log}$ over $k$.

Next, let us take a triple $(\rho_0, \rho_1, \rho_\infty)$ of elements of $(\mbZ /p^{\N} \mbZ)^{\times}/\{ \pm 1 \}$.
There exists  the triple of integers $(\lambda_0, \lambda_1, \lambda_\infty)$   satisfying the following conditions:
\begin{itemize}
\item[(a)]
$2  \cdot \rho_x = \lambda_x$ as elements of  $(\mbZ/p^\N \mbZ)/\{ \pm 1 \}$ and 
$0 < \lambda_x < p^{\N}$
 for every  $x=0,1, \infty$;
\item[(b)]
The sum $\lambda_0 + \lambda_1 + \lambda_\infty$ is odd $< 2 \cdot p^{\N}$.
\end{itemize}
Let us write $\mcO_{\mbP}^+ := \mcO_{\mbP} (\lambda_0 \cdot [0]+ \lambda_1 \cdot [1] + \lambda_{\infty} \cdot [\infty])$.
Note that there is a unique $\mcD_{\mbP^{\mr{log}}}^{(\N -1)}$-module structure  $\nabla^+$ on $\mcO_{\mbP}^+$ whose restriction to    $U := \mbP \setminus \{ [0], [1], [\infty] \}$ coincides with  the trivial $\mcD_{U}^{(\N-1)}$-module structure  on $\mcO_{\mbP}^+ |_U = \mcO_U$ (cf. Remark \ref{roafos}).

The following assertion is a special case of ~\cite[\S\,3, Theorem 3.3, (ii)]{Oss1}; we shall prove it by a relatively elementary argument.

\SSP
\bpr \label{Loo09}
Let $\phi : \mbP \migi \mbP$ be a tamely ramified covering classified by
$\mr{Cov}^{\mr{tame}}_{\msP, (\lambda_0, \lambda_1, \lambda_\infty)}$.
Then, 
 the points $\phi ([0]), \phi ([1]), \phi ([\infty])$ are mutually distinct.
\epr
\begin{proof}
First,  we shall suppose that the set $\phi (\{ [0], [1], [\infty] \})$ consists of one point.
By considering the fiber of $\phi$ over this point, we see that  $\mr{deg}(\phi) \geq \lambda_0 +\lambda_1 + \lambda_\infty$.
It follows from the Riemann-Hurwitz formula that 
\begin{align}
-2  \left(= 2 \cdot (\text{genus of $\mbP$}) -2 \right)= -2 \cdot \mr{deg}(\phi) + \sum_{x=0,1,\infty} (\lambda_x -1) \leq -(\lambda_0 + \lambda_1 + \lambda_\infty) -3.
\end{align}
Thus, we obtain a contradiction.
Next,  suppose that $\phi (\{[0], [1], [\infty]  \})$ consists of two  points.
After possibly applying a linear transformation, we may assume that these two elements coincide with $\{[0],  [\infty]\}$, and that $\{ [0] \} \subseteq \phi^{-1}([0])$ and $\{ [1], [\infty] \} \subseteq \phi^{-1}([\infty])$.
Hence, $\phi$ defines a tamely ramified  covering of $\mbG_m \left(= \mbP \setminus \{ [0], [\infty] \}\right)$.
Recall that the tame fundamental group $\pi^\mr{tame}_1 (\mbG_m)$ of $\mbG_m$ is isomorphic to $\widehat{\mbZ}^{p'}$, the maximal prime-to-$p$ quotient of $\widehat{\mbZ} := \varprojlim_{n \in \mbZ_{\geq 0}}\mbZ/n\mbZ$.
A topological generator $\sigma$ of $\pi^{\mr{tame}}_1 (\mbG_m)$ 
acts on the fiber over the point near $[0]$ as a cyclic permutation.
On the other hand, $\sigma^{-1}$ acts on the fiber over the point near $[\infty]$ as a product of two disjoint  cyclic permutations.
This is a contradiction.
Hence, the the image $\phi (\{ [0], [1], [\infty] \})$ consists of three  points.
This completes the proof of this assertion.
\end{proof}
\SSP

\begin{rema} \label{R453}
Because of  Proposition \ref{Loo09}, each
element of $\mr{Cov}^{\mr{tame}}_{\msP, (\lambda_0, \lambda_1, \lambda_\infty)}$  has a unique representative 
$\phi : \mbP \migi \mbP$ satisfying $\phi ([x])=[x]$ for every $x = 0, 1,\infty$.
Following ~\cite{ABEGKM} (or ~\cite{BEK}), we call such a covering  a {\it dynamical Belyi map}.

Now, let us take a dynamical Belyi map $\phi : \mbP \migi \mbP$ classified by $\mr{Cov}^{\mr{tame}}_{\msP, (\lambda_0, \lambda_1, \lambda_\infty)}$;
this covering corresponds to 
a representation of the tame fundamental group $\pi_1^{\mr{tame}}(\mbP \setminus \{ [0], [1], [\infty]\})$ (which is obtained as a quotient of the profinite completion $\widehat{\Pi}_{0, 3}$ of the group $\Pi_{0, 3} := \langle \gamma_0, \gamma_1, \gamma_\infty \, | \, \gamma_0 \gamma_1 \gamma_\infty =1 \rangle$).
Hence, by  
the above proposition, 
$\phi$ determines  
 three cyclic permutations $\sigma_0$, $\sigma_1$, and $\sigma_\infty$ (in the
 symmetric group $\mfS_d$ of $d$ letters for some $d \geq 1$)
   of orders $\lambda_0$, $\lambda_1$, and $\lambda_\infty$, respectively, satisfying 
$\sigma_0 \circ \sigma_1 = \sigma_\infty$.
A trivial elementary argument shows that this condition implies  the following inequalities:
\begin{align} \label{eapahks}
|\lambda_0 -\lambda_1| < \lambda_\infty <  \lambda_0 + \lambda_1.
\end{align}
These inequalities also can be obtained by  the inequality $\mr{deg}(\phi)\left(= \frac{\lambda_0 + \lambda_1 + \lambda_\infty -1}{2} \right) \geq \lambda_0, \lambda_1, \lambda_\infty$.

Conversely, suppose that a  subgroup of $\mfS_d$ generated by three  cyclic permutations $\sigma_0$, $\sigma_1$, $\sigma_\infty$ with $\sigma_0 \circ \sigma_1 = \sigma_\infty$ has order prime to $p$.
Then, 
the assignment  $\gamma_x \mapsto \sigma_x$ ($x = 0,1,\infty$)  induces a representation $\pi_1^{\mr{tame}}(\mbP \setminus \{ [0], [1], [\infty]\}) \migi \mfS_d$ because 
the surjection   $\widehat{\Pi}_{0, 3} \migisurj \pi_1^{\mr{tame}}(\mbP \setminus \{ [0], [1], [\infty]\}) $ becomes bijective after taking their maximal prime-to-$p$ quotients.
In particular, the corresponding  tamely ramified covering is  classified by $\mr{Cov}^\mr{tame}_{\msP, (\lambda_0, \lambda_1, \lambda_\infty)}$.
See ~\cite{BEK} for the study concerning  the relationship between such cyclic permutations and  dynamical Belyi maps in characteristic $p$.
\end{rema}

\LSP
\subsection{Dormant $\mr{PGL}_2$-opers on a $3$-pointed projective line} \label{SS0jj42}

In what follows,  we shall  prove that the map $\Upsilon_{\msX, \vec{\lambda}}$ defined in (\ref{eavbk}) becomes a bijection in the case where $\msX = \msP$.
To do this, we construct its inverse map.
We first prove the following proposition.

\SSP
\bpr \label{Pioe}
Let $\msF^\spadesuit$ be a dormant $\mr{PGL}_2$-oper of level $\N$ on $\msP$.
Then, there exists a dormant $\mr{GL}_2$-oper $\msF^\heartsuit := (\mcF, \nabla, \mcL)$   of level $\N$  on $\msP$ 
 satisfying $\msF^{\heartsuit \Rightarrow \spadesuit} \cong \msF^\spadesuit$ and  $\mr{det}(\mcF, \nabla) \cong (\mcO_{\mbP}^+, \nabla^+)$,
Moreover,  such a $\mr{GL}_2$-oper   is uniquely determined up
  isomorphism.
\epr
\begin{proof}
First, we prove  the existence portion.
Let us take a dormant $\mr{GL}_2$-oper $\msF^\heartsuit := (\mcF, \nabla, \mcL)$ of level $\N$ on $\msP$ with $\msF^{\heartsuit \Rightarrow \spadesuit} \cong \msF^\spadesuit$.
Since $\mcL \cong \Omega_{\mbP^{\mr{log}}} \otimes (\mcF/\mcL)$ (cf. (\ref{fai752})),  the following equalities hold: 
\begin{align}
\mr{deg}(\mcO^+_{\mbP} \otimes \mr{det}(\mcF)^\vee) &= \lambda_0 + \lambda_1 + \lambda_\infty - \mr{deg}(\mcL) - \mr{deg} (\mcF/\mcL)  \\
& = \lambda_0 + \lambda_1 + \lambda_\infty  +1 - 2 \cdot  \mr{deg}(\mcL). \notag
\end{align}
In particular, the degree of $\mcO^+_{\mbP} \otimes \mr{det}(\mcF)^\vee$
 is even.
Hence, it follows from Lemma \ref{Lkg09}  that
there exists a  dormant  $\mcD_{\mbP^\mr{log}}^{(\N-1)}$-bundle $(\mcN, \nabla_\mcN)$  such that  $\mcN$ is a line bundle with  $\mcN^{\otimes 2} \cong (\mcO_{\mbP}^+, \nabla^+)\otimes\mr{det}(\mcF, \nabla_\mcF)^\vee$.
By putting  $(\mcF', \nabla', \mcL') := \msF^\heartsuit_{\otimes (\mcN, \nabla_\mcN)}$, 
we have
\begin{align}
\mr{det}(\mcF', \nabla') \cong (\mcN, \nabla_\mcN)^{\otimes 2} \otimes \mr{det}(\mcF, \nabla) \cong  (\mcO^+_{\mbP}, \nabla^+).
\end{align}
Thus, $(\mcF', \nabla', \mcL')$ specifies the required $\mr{GL}_2$-oper of level $\N$.

Next, we shall prove the uniqueness portion.
Let $\msF^\heartsuit_i := (\mcF_i, \nabla_i, \mcL_i)$ ($i=1,2$) be dormant $\mr{GL}_2$-opers  of level $\N$ on $\msP$ satisfying the required conditions.
Since $\msF_1^{\heartsuit \Rightarrow \spadesuit} = \msF_2^{\heartsuit \Rightarrow \spadesuit}$,
there exists  a dormant $\mcD_{\mbP^\mr{log}}^{(\N-1)}$-module $(\mcN, \nabla_\mcN)$ such that $\mcN$ is a line bundle  and $\msF_2^{\heartsuit} \cong (\msF_1^{\heartsuit})_{\otimes (\mcN, \nabla_\mcN)}$.
If $\nabla^\mr{triv}$ denotes the  trivial $\mcD_{\mbP^\mr{log}}^{(\N-1)}$-module structure on $\mcO_\mbP$, then we have
\begin{align}
(\mcO_\mbP, \nabla^\mr{triv}) &\cong  \mr{det}(\mcF_2, \nabla_2) \otimes (\mcO_\mbP^+, \nabla^+)^\vee \\
&\cong  \mr{det} ((\mcN, \nabla_\mcN) \otimes (\mcF_1, \nabla_1)) \otimes (\mcO_\mbP^+, \nabla^+)^\vee\notag \\
&\cong 
(\mcN, \nabla_\mcN)^{\otimes 2} \otimes \mr{det}(\mcF_1, \nabla_1) \otimes (\mcO_\mbP^+, \nabla^+)^\vee\notag \\
&\cong (\mcN, \nabla_\mcN)^{\otimes 2}.\notag
\end{align} 
Since $\mr{Pic}(\mbP)\cong [\mcO_\mbP (1)] \cdot \mbZ$,
the line bundle $\mcN$ may be identified with $\mcO_\mbP$.
Moreover, by the uniqueness portion of Proposition \ref{Lkg09},
$\nabla_\mcN$ coincides with $\nabla^\mr{triv}$ via a fixed identification $\mcN = \mcO_\mbP$.
Thus, 
we have $(\msF_1^{\heartsuit})_{\otimes (\mcN, \nabla_\mcN)} \cong \msF_1^{\heartsuit}$, which implies  that
$\msF^\heartsuit_2$ is isomorphic to $\msF^\heartsuit_1$.
This completes the proof of the uniqueness portion.
\end{proof}
\SSP

Now, let us take a dormant $\mr{PGL}_2$-oper $\msF^\spadesuit$  on $\msP$ of  level $\N$ and radii $(\rho_0, \rho_1, \rho_\infty)$.
 Also, let  $\msF^\heartsuit := (\mcF, \nabla, \mcL)$ be the dormant $\mr{GL}_2$-oper of level $\N$ resulting from Proposition \ref{Pioe} applied to  $\msF^\spadesuit$.
In particular, we have  $\mr{deg}(\mcL) = \frac{\lambda_0 + \lambda_1 + \lambda_\infty +1}{2}$
and $\mr{deg}(\mcF/\mcL) = \frac{\lambda_0 + \lambda_1 + \lambda_\infty -1}{2}$.
Denote by $\tau$ the $\mcO_\mbP$-linear morphism $F^{(\N)*}_{\mbP/k}(\mcF^\nabla) \migi \mcF$
 extending the ($\mcO_{\mbP^{(\N)}}$-linear)  inclusion $\mcF^\nabla \migiincl \mcF$;
 the morphism $\tau$ is compatible with $\nabla^\mr{can}_{\mcF^\nabla}$ (cf. Remark \ref{eRfi})  and $\nabla$.
We shall write $\mcL^\natural := \mcL \cap \mr{Im}(\tau)$.
Since the restriction of $\tau$ to $U := \mbP \setminus \{ [0], [1], [\infty] \}$ is an isomorphism,
the quotient sheaf $\mcL/\mcL^\natural$ is a torsion sheaf supported on $\{ [0], [1], [\infty] \}$.

\SSP
\ble \label{LLdjw3}
The length of  $\mcL/\mcL^\natural$ at the marked point $x \in \{[0], [1], [\infty] \}$ is  $\lambda_x$.
 Moreover, the natural inclusion $\mcL/\mcL^\natural \migiincl \mr{Coker}(\tau)$ is an isomorphism.
\ele
\begin{proof}
Let us fix $x \in \{[0], [1], [\infty] \}$, and 
choose  a local function  $t$ defining  $x$.
The restriction  of $(\mcF, \nabla)$ to the formal neighborhood $\widehat{D}_x$ of $x$ may  be expressed as
$(k[\![t]\!], \widehat{\nabla}_{\overline{a}}) \oplus (k[\![t]\!], \widehat{\nabla}_{\overline{b}})$
for some integers $a, b$
  with $0 \leq b \leq a \leq p^{\N}-1$.
The radius of $\nabla$ at $x$ coincides with $\rho_x$ by assumption, so the equality  $2 \cdot  \rho_x =  \overline{a-b}$  of  elements in  $(\mbZ/p^\N \mbZ)/\{ \pm 1 \}$ holds.
Since $\mr{det}(\mcF, \nabla) \cong (\mcO_{\mbP}^+, \nabla^{+})$,
a computation of the lengths at $x$ of 
$\mr{det}(\mcF, \nabla)$ and $(\mcO_{\mbP}^+, \nabla^+)$ implies  
  $a + b \equiv \lambda_x$ (mod $p^{\N}$).
Hence, since $a + b \leq 2 \cdot p^{\N}$, it follows that $a +b$ is either $\lambda_x$ or $\lambda_x + p^{\N}$.
We shall prove the claim that
$a + b = \lambda_x$.
Suppose, on the contrary,  that $a + b = p^{\N}+ \lambda_x$.
The condition that $ 0 \leq a- b \leq p^{\N}-1$ and $2 \cdot  \rho_x = \overline{a-b}$ in $(\mbZ/p^\N \mbZ)/\{ \pm 1 \}$ implies 
  that $a - b$ is either $\lambda_x$ or $p^{\N}-\lambda_x$.
Since $a + b$ and $a -b$ have the same parity,
the equality   $a - b = p^{\N}-\lambda_x$ holds.
Hence,  we have 
\begin{align}
2 \cdot  a = (a + b)+ (a-b) = (p^\N + \lambda_x) + (p^\N -\lambda_x) =  2 \cdot p^{\N}.
\end{align}
 This implies the equality  $a = p^{\N}$, which is a contradiction.
This proves the claim.

By the claim proved just now,
we have $a + b = \lambda_x$ and $a - b = \lambda_x$.
In particular, $a = \lambda_x$ and $b= 0$, which means that
the restriction  of $(\mcF, \nabla)$ to $\widehat{D}_x$ is isomorphic to $(k[\! [t]\!], \widehat{\nabla}_{\overline{\lambda}_x}) \oplus (k[\! [t]\!], \widehat{\nabla}_0)$. 
 It follows that the restriction of $\mr{Im}(\tau)$ to $\widehat{D}_x$
 coincides with $t^{\lambda_x}\cdot k[\! [t]\!] \oplus k[\! [t]\!] \left( \subseteq k[\! [t]\!]^{\oplus 2}\right)$.
On the other hand, according to the proof in  Proposition \ref{C001}, the inclusion $\mcL \migiincl \mcF$ corresponds, after choosing a suitable trivialization $\Gamma (\widehat{D}_{x}, \mcL |_{\widehat{D}_{x}}) \isom k[\![t]\!]$,
to the $k[\![t]\!]$-linear morphism $k[\![t]\!] \migi k[\![t]\!]^{\oplus 2}$  given by $1 \mapsto (u_1, u_2)$ for some $u_1, u_2 \in k[\![t]\!]^\times$.
By taking account of this observation, we see that
 the length of $\mcL/\mcL^\natural$ at $x$ coincides with $\lambda_x$.
 Moreover, since  the length of $\mr{Coker}(\tau)$ at $x$ is  $a +  b = \lambda_x$,
 the inclusion $\mcL/\mcL^\natural \migiincl \mr{Coker}(\tau)$ turns out to be  an isomorphism.
 \end{proof}

\SSP
\ble \label{L009}
There exists an $\mcO_{\mbP^{(\N)}}$-linear isomorphism $\gamma :  \mcO_{\mbP^{(\N)}}^{\oplus 2} \isom \mcF^\nabla$.
\ele
\begin{proof}
Since $\mcF^\nabla$ is a rank $2$ vector bundle on the genus-$0$ curve $\mbP^{(\N)}$, it is isomorphic to the direct sum of two line bundles.
Let us fix an isomorphism  $\gamma : \mcO_{\mbP^{ (\N)}}(a) \oplus \mcO_{\mbP^{ (\N)}} (b) \isom \mcF^\nabla$, where $a$ and $b$ are some  integers with $a \geq b$.
The pull-back of $\gamma$ via $F^{(\N)}_{\mbP/k}$ defines an isomorphism 
 $\gamma^F : \mcO_{\mbP}(p^\N \cdot a) \oplus \mcO_{\mbP}(p^\N \cdot  b) \isom F^{(\N)*}_{\mbP/k}(\mcF^\nabla)$.
It follows from  Lemma \ref{LLdjw3} that
\begin{align}
\lambda_0 + \lambda_1 + \lambda_\infty &= \mr{length}(\mcL/\mcL^\natural) \\
& = \mr{length}(\mr{Coker}(\tau)) \notag \\
&= \mr{deg}(\mcF) - \mr{deg}(F^{\N*}(\mcF^\nabla)) \notag \\
&= \lambda_0 + \lambda_1 + \lambda_\infty + p^N (a+b).\notag
\end{align}
This implies $b = -a$.
Next, let us consider the composite
\begin{align} \label{Er45}
\mcO_{\mbP} (p^\N \cdot a) \xrightarrow{v \mapsto (v, 0)}\mcO_{\mbP} (p^\N \cdot a) \oplus \mcO_{\mbP}(-p^\N \cdot a) \xrightarrow{\gamma^F} F^{(\N) *}_{\mbP/k}(\mcF^\nabla) \xrightarrow{\tau} \mcF \migisurj \mcF/\mcL.
\end{align}
We shall suppose that $a >0$.
The assumption  $\lambda_0 + \lambda_1 + \lambda_\infty < 2 \cdot p^\N$ implies
\begin{align}
\mr{deg}(\mcO_{\mbP}(p^\N \cdot  a)) \geq p^\N  > \frac{\lambda_0 + \lambda_1 + \lambda_\infty -1}{2} = \mr{deg}(\mcF/\mcL).
\end{align}
Hence, the composite (\ref{Er45}) must be the zero map.
It follows that the image $\mcI$ of the inclusion into the first factor $\mcO_{\mbP}(p^\N \cdot a) \migiincl \mcO_{\mbP} (p^\N \cdot a) \oplus \mcO_{\mbP}(-p^\N \cdot  a)$ is contained in $\mcL \left(\subseteq \mcF \right)$ via $\tau \circ \gamma^F$.
But, since $\mcI$ is closed under 
$\nabla_{\mcO_{\mbP^{(\N)}}(a)}^\mr{can} \oplus \nabla_{\mcO_{\mbP^{(\N)}}(b)}^\mr{can}$,
the line subbundle $\mcL$ must be  closed under $\nabla$.
This contradicts the fact that  $\mr{KS}_{\msF^\heartsuit}^1$ (cf. (\ref{eraopaz98}))
is an isomorphism.
Thus, we conclude that $a = 0$, i.e., $\gamma$ defines an isomorphism $\mcO_{\mbP}^{\oplus 2} \isom \mcF^\nabla$.
\end{proof}
\SSP

Let $\gamma$ be as asserted in the above lemma.
For convenience, we occasionally  use the notation $X$ to denote the underlying projective line  $\mbP$.
The pull-back $\gamma^F : \mcO_{X}^{\oplus 2} \isom F^{(\N)*}_{X/k}(\mcF^\nabla)$ of $\gamma$ by $F^{(\N)}_{X/k}$
induces a trivialization  $\mbP (\gamma^F) : \mbP \times_k X \isom \mbP (F^{(\N)*}_{X/k}(\mcF^\nabla))$ of the $\mbP$-bundle $\mbP (F^{(\N)*}_{X/k}(\mcF^\nabla))$ associated to $F^{(\N)*}_{X/k}(\mcF^\nabla)$.
The sheaf  $\mcL^\natural$, regarded as a line bundle of $F^{(\N)*}_{X/k}(\mcF^\nabla)$ via $\tau$,   defines a global section $\sigma : X \migi \mbP (F^{(\N)*}_{X/k}(\mcF^\nabla))$.
Thus, we obtain the composite
\begin{align}
\phi_{\msF^\spadesuit} : X \xrightarrow{\sigma}\mbP (F^{(\N)*}_{X/k}(\mcF^\nabla)) \xrightarrow{\mbP (\gamma^F)^{-1}} \mbP \times_k X \xrightarrow{\mr{pr}_1} \mbP.
\end{align}

\SSP
\ble \label{L008}
The morphism $\phi_{\msF^\spadesuit} : X \migi \mbP^1$ defines 
a tamely ramified covering classified by  the set $\mr{Cov}^\mr{tame}_{\msP,  (\lambda_0, \lambda_1, \lambda_\infty)}$.
\ele
\begin{proof}
Let $x$ be a $k$-rational point of $X$.
To complete the proof, we shall describe the morphism $\phi_{\msF^\spadesuit}$ restricted to   the formal neighborhood $\widehat{D}_x$ of $x$.

First, suppose that $x \in \{ [0], [1], [\infty]\}$.
As observed in the proof of Lemma \ref{LLdjw3}, 
 the $\breve{\mcD}_{k[[t]]}^{(\N-1)}$-module corresponding to the restriction of $(\mcF, \nabla)$ to $\widehat{D}_x$ is isomorphic to $(t^{-\lambda_x} \cdot k[\![t]\!], \widehat{\nabla}_0)\oplus (k[\![t]\!], \widehat{\nabla}_0) \left(\cong (k[\![t]\!], \widehat{\nabla}_{\overline{\lambda}_x})\oplus (k[\![t]\!], \widehat{\nabla}_0) \right)$.
According to the proof of Proposition \ref{C001}, the inclusion $\mcL \migiincl \mcF$  restricted to $\widehat{D}_x$ may be identified, after choosing a suitable trivialization $\Gamma (\widehat{D}_x, \mcL |_{\widehat{D}_x}) \isom k[\![t]\!]$,
with the $k[\![t]\!]$-linear morphism $k[\![t]\!] \migi t^{-\lambda_x} \cdot k[\![t]\!]  \oplus k[\![t]\!]$  determined   by $1 \mapsto (t^{-\lambda_x} \cdot 1, u)$ for some $u \in  k[\![t]\!]^\times$.
Under this identification,
the inclusion $\mcL^\natural \migiincl F_{X/k}^{(\N)*}(\mcF^\nabla)$ restricted to $\widehat{D}_x$ corresponds to the inclusion $t^{\lambda_x} \cdot k[\![t]\!] \migi k[\![t]\!]^{\oplus 2}$ determined  by $t^{\lambda_x} \cdot 1 \mapsto (1, t^{\lambda_x} \cdot u)$.
This implies that the restriction of $\phi_{\msF^\spadesuit}$ to $\widehat{D}_x$ arises from  the $k$-algebra endomorphism of $k[\![t]\!]$ given  by $t \mapsto t^{\lambda_x} \cdot u$.
Hence, the ramification index  of $\phi_{\msF^\spadesuit}$ at $x$ coincides with $\lambda_x$.

Next, suppose that  $x \in X \setminus \{ [0], [1], [\infty] \}$.
The inclusion $\mcL \migiincl \mcF$ restricted to $\widehat{D}_x$ may be described, after choosing suitable trivializations of $\mcL |_{\widehat{D}_x}$ and $\mcF |_{\widehat{D}_x}$,
as the $k[\![t]\!]$-linear morphism $k[\![t]\!] \migi k[\![t]\!]^{\oplus 2}$ determined by $1 \mapsto (1, t \cdot v)$ for some $v \in k[\![t]\!]^\times$.
Hence, by the same reason as above, 
the ramification  index of $\phi_{\msF^\spadesuit}$ at $x$ turns out to be  $1$, which means that
$\phi_{\msF^\spadesuit}$  is  \'{e}tale  at $x$.
This completes the proof of this lemma.
\end{proof}
\SSP

Note  that the resulting element $[\phi_{\msF^\spadesuit}] \in \mr{Cov}^\mr{tame}_{\msP,  (\lambda_0, \lambda_1, \lambda_\infty)}$ does not depend on the choice of the trivialization $\gamma : \mcO_{\mbP^{(\N)}}^{\oplus 2} \isom \mcF^\nabla$ because 
$\gamma$ is uniquely determined up to forward composition with an element of  $\mr{Aut}_{\mcO_{\mbP^{(\N)}}} (\mcO_{\mbP^{(\N)}}^{\oplus 2}) \left(= \mr{PGL}_2 (k)\right)$.
Thus, the assignment $\msF^\spadesuit \mapsto [\phi_{\msF^\spadesuit}]$ determines a well-defined map of sets  
${^\N}\mcO p^{^\mr{Zzz...}}_{2, \msP, (\rho_0, \rho_1, \rho_\infty)} \migi \mr{Cov}^\mr{tame}_{\msP,  (\lambda_0, \lambda_1, \lambda_\infty)}$.
One may verify that this map  specifies, by construction,   the inverse to 
the map $\Upsilon_{\msP, (\lambda_0, \lambda_1, \lambda_\infty)}$.
Thus, we have obtained the following assertion.

\SSP
\bpr \label{P00f45}
The assignments $[\phi] \mapsto \msF^\spadesuit_{\phi}$ (i.e, $\Upsilon_{\msP, (\lambda_0, \lambda_1, \lambda_\infty)}$) and  $\msF^\spadesuit \mapsto [\phi_{\msF^\spadesuit}]$ constructed above 
give a bijective correspondence 
\begin{align}
{^\N}\mcO p^{^\mr{Zzz...}}_{2, \msP, (\rho_0, \rho_1, \rho_\infty)}\cong\mr{Cov}^\mr{tame}_{\msP,  (\lambda_0, \lambda_1, \lambda_\infty)}.
\end{align}
\epr

\LSP
\subsection{Correspondence between tamely ramified coverings and dormant $\mr{PGL}_2$-opers} \label{SE42}

The following proposition  (together with Proposition \ref{P00f45}) may be regarded as a variant of the rigidity assertion for dynamical Belyi maps proved  in ~\cite[\S\,2, Lemma 2.1]{LiOs2}.

\SSP
\bpr \label{P00245}
Let $(\rho_0, \rho_1, \rho_\infty)$ be   an element of $((\mbZ/p^\N \mbZ)/\{\pm 1 \})^{\times 3}$.
Then, a dormant $\mr{PGL}_2$-oper  on $\msP$ of level $\N$ and  radii $(\rho_0, \rho_1, \rho_\infty)$ is, if it exists,  uniquely determined.
That is to say,  the following inequality holds: 
\begin{align}
\sharp ({^\N}\mcO p^{^\mr{Zzz...}}_{2, \msP, (\rho_0, \rho_1, \rho_\infty)}) 
\leq 1.
\end{align}
\epr
\begin{proof}
Suppose that  ${^\N}\mcO p^{^\mr{Zzz...}}_{2, \msP, (\rho_0, \rho_1, \rho_\infty)} \neq \emptyset$.
By the canonical morphism 
$\mcD_{\mbP^\mr{log}}^{(0)} \migi \mcD_{\mbP^\mr{log}}^{(\N-1)}$,   each 
element of ${^\N}\mcO p^{^\mr{Zzz...}}_{2, \msP, (\rho_0, \rho_1, \rho_\infty)}$ 
 induces a dormant $\mr{PGL}_2$-oper of level $1$ and radii 
 $(\overline{\rho}_0, \overline{\rho}_1, \overline{\rho}_\infty)$, where $\overline{\rho}_x$ (for $x = 0,1, \infty$) denotes the image of $\rho_x$ via the natural surjection $(\mbZ/p^\N \mbZ)/\{ \pm 1 \} \migisurj \mbF_p/\{ \pm 1 \}$.
  Hence,   ${^1}\mcO p^{^\mr{Zzz...}}_{2, \msP, (\overline{\rho}_0, \overline{\rho}_1,  \overline{\rho}_\infty)}$ must be nonempty.
This fact together with a comment in Remark \ref{peec} implies  that $(\overline{\rho}_0, \overline{\rho}_1,  \overline{\rho}_\infty) \in (\mbF_p^\times/\{ \pm 1 \})^{\times 3}$, or equivalently, $(\rho_0, \rho_1, \rho_\infty) \in ((\mbZ/p^\N \mbZ)^\times/\{ \pm 1 \})^{\times 3}$.
In  particular, there exists a triple of integers  $(\lambda_0, \lambda_1, \lambda_\infty)$ associated to $(\rho_0, \rho_1, \rho_\infty)$ satisfying the conditions (a), (b)  described   in \S\,\ref{SS042e}. 

Now, suppose that there exist two dormant $\mr{PGL}_2$-opers   $\msF_1^\spadesuit$,  $\msF^\spadesuit_2$ on $\msP$ of level $\N$ and  radii $(\rho_0, \rho_1, \rho_\infty)$.
For each $j=1,2$,
denote by $\msF^\heartsuit_j := (\mcF_j, \nabla_j, \mcL_j)$ the dormant  $\mr{GL}_n$-oper of level $\N$  resulting from Proposition \ref{Pioe} applied to $\msF^\spadesuit_j$.
The inclusion $\mcF_j^\nabla \migiincl \mcF_j$ extends to 
 an $\mcO_{\mbP}$-linear morphism $\tau_j : F^{(\N)*}_{\mbP/k} (\mcF_j^\nabla) \migi \mcF_j$.
Let us write $\mcL_j^\natural := \mcL_j \cap \mr{Im}(\tau_j)$ and write
$\iota_j$ for the natural isomorphism $\mcL_j /\mcL_j^\natural \isom \mr{Coker}(\tau_j)$ (cf. Lemma  \ref{LLdjw3}).
Also, denote by $\overline{\nabla}_j$ the $\mcD_{\mbP^\mr{log}}^{(0)}$-module structure on $\mcF_j$ induced from $\nabla_j$.
The triple $\overline{\msF}_j^{\heartsuit} := (\mcF_j, \overline{\nabla}_j, \mcL_j)$ defines a dormant $\mr{GL}_2$-oper in the classical sense, i.e., of level $1$.
In particular, it induces a dormant $\mr{PGL}_2$-oper $\overline{\msF}_i^{\heartsuit \Rightarrow \spadesuit}$  on $\msP$ of level $1$ and radii 
 $(\overline{\rho}_0, \overline{\rho}_1, \overline{\rho}_\infty)$.
Recall from ~\cite[Chap.\,I, \S\,4, Theorem 4.4]{Mzk2}  (cf. Remark \ref{peec}) that   dormant $\mr{PGL}_2$-opers of level $1$ on $\msP$ are  completely determined by their radii.
It follows that $\overline{\msF}_1^{\heartsuit \Rightarrow \spadesuit}  = \overline{\msF}_2^{\heartsuit \Rightarrow \spadesuit}$.
By the uniqueness assertion in Proposition \ref{Pioe},
there exists an isomorphism of $\mr{GL}_2$-opers $\alpha : \overline{\msF}_1^{\heartsuit}  \isom  \overline{\msF}_2^{\heartsuit}$.
This isomorphism 
is restricted to an isomorphism $\alpha |_{\mcL_1} : \mcL_1 \isom \mcL_2$, which  induces, via taking the respective quotients,  an isomorphism $\alpha |_{\mcL_1/\mcL^\natural_1} : \mcL_1/\mcL_1^\natural \isom \mcL_2/\mcL^\natural_2$.
The composite $\alpha |_{\mr{Coker}(\tau_1)} := \iota_2 \circ \alpha |_{\mcL_1/\mcL^\natural_1} \circ  \iota_1^{-1}$ specifies an isomorphism
$\mr{Coker}(\tau_1) \isom \mr{Coker}(\tau_2)$.

In what follows, we shall prove the commutativity of the following square diagram:
\begin{align} \label{Eed3df}
\vcenter{\xymatrix@C=56pt@R=36pt{
\mcF_1 \ar[r]^-{\alpha} \ar[d]_-{\pi_1} & \mcF_2\ar[d]^{\pi_2} 
\\
\mr{Coker}(\tau_1) \ar[r]_-{\alpha |_{\mr{Coker}(\tau_1)}} & \mr{Coker}(\tau_2),
}}
\end{align}
where $\pi_j$ ($j=1,2$) denotes the natural projection $\mcF_j \migisurj \mr{Coker}(\tau_j)$.
Let us take $j \in \{1, 2 \}$,  $x \in \{ [0], [1], [\infty] \}$.
Also, choose  a local function $t$ on $\mbP$ defining $x$.
Denote by $\widehat{D}_x$ the formal neighborhood of $x$ in $\mbP$, which may be identified with $\mr{Spec}(k[\![t]\!])$.
Fix an identification of 
the restriction of $(\mcF_j, \nabla_j)$ to $\widehat{D}_x$ with   $(k[\! [t]\!], \widehat{\nabla}_{\overline{\lambda}_x}) \oplus (k[\! [t]\!], \widehat{\nabla}_0)$ (cf. the proof of Lemma \ref{LLdjw3}).
As observed  in  the proof of  Proposition \ref{C001},
the inclusion $\mcL \migiincl \mcF$ corresponds, after choosing a suitable trivialization
$\Gamma (\widehat{D}_x, \mcL |_{\widehat{D}_{x}}) \isom k[\! [t]\!]$,
to the $k[\! [t]\!]$-linear morphism $k[\! [t]\!] \migi k[\! [t]\!]^{\oplus 2}$ given by 
$1 \mapsto (u_j, 1)$ for some $u_j \in k[\! [t]\!]^\times$.
The isomorphism $\alpha$ restricted to  $\widehat{D}_x$ defines an  automorphism of  $(k[\! [t]\!], \widehat{\nabla}_{\overline{\lambda}_x}) \oplus (k[\! [t]\!], \widehat{\nabla}_0)$.
But, since $\rho_x \neq 0$ in $\mbF_p/\{\pm1 \}$ (or equivalently,  $\overline{\lambda}_x \neq 0$),  it follows from Proposition \ref{P0045}, (ii), that this automorphism may be expressed as $\mu_{v} \oplus \mu_w$ for some $v, w \in k[\! [t]\!]^\times$, where $\mu_{(-)}$ denotes the endomorphism of $k[\! [t]\!]$ given  by multiplication by  $(-)$.
Since $\alpha |_{\mcL_1}$ is obtained as the restriction of $\alpha$,
the restriction of  $\alpha |_{\mcL_1}$ to $\widehat{D}_x$ may be expressed as $\mu_w$, and the equality $v u_1 = w u_2$ holds.
Hence, for each $(g, h) \in k[\! [t]\!]^{\oplus 2} \left(= \mcF_1 |_{\widehat{D}_x} \right)$,
we have
\begin{align}
(\pi_2\circ\alpha) ((g, h)) & = \pi_2 ((vg, wh)) \\
& =  \left(\frac{vg}{u_2} \cdot u_2, \frac{vg}{u_2} \cdot 1 \right)  \text{mod} \ \mr{Im}(\tau_2) \notag \\
& = \iota_2 \left( \frac{vg}{u_2} \ \mr{mod} \ \mcL^\natural_2\right) \notag  \\
& = (\iota_2 \circ \alpha |_{\mcL_1/\mcL^\natural_1}) \left(\frac{vg}{w u_2}  \ \mr{mod} \ \mcL^\natural_1 \right)  \notag \\
& = (\alpha |_{\mr{Coker}(\tau_1)}\circ\iota_1) \left(\frac{g}{u_1}  \ \mr{mod} \ \mcL^\natural_1 \right)\notag \\
& =  \alpha |_{\mr{Coker}(\tau_1)} \left(\left(\frac{g}{u_1} \cdot u_1, \frac{g}{u_1} \cdot 1 \right) \text{mod} \ \mr{Im}(\tau_1) \right) \notag \\
& = (\alpha |_{\mr{Coker}(\tau_1)} \circ \pi_1) ((g, h)). \notag
\end{align}
This shows the desired commutativity of  (\ref{Eed3df}).

Moreover, the commutativity of (\ref{Eed3df}) just proved  implies that $\alpha$ is restricted, vis $\iota_1$ and $\iota_2$,  to an isomorphism $\alpha' : F^{(\N)*}_{\mbP/k}(\mcF_1^\nabla) \isom F^{(\N)*}_{\mbP/k}(\mcF_2^\nabla)$.
Since $\mcF_1^\nabla \cong \mcF_2^\nabla \cong \mcO_{\mbP^{(\N)}}^{\oplus 2}$ (cf. Lemma \ref{L009}) and $k$ is perfect, the morphism 
\begin{align}
\mr{Hom}_{\mcO_{\mbP^{(\N)}}} (\mcF^\nabla_1, \mcF_2^\nabla) \migi \mr{End}_{\mcO_\mbP} (F^{(\N)*}_{\mbP/k}(\mcF_1^\nabla), F^{(\N)*}_{\mbP/k}(\mcF_2^\nabla))
\end{align}
 arising  from  pull-back by  $F^{(\N)}_{\mbP/k}$ is bijective.
In particular, $\alpha'$ comes from an isomorphism $\mcF_1^\nabla \isom \mcF_2^\nabla$, and hence $\alpha'$ is compatible with the respective $\mcD_{\mbP^\mr{log}}^{(\N-1)}$-actions $\nabla^\mr{can}_{\mcF^\nabla_1}$, $\nabla^\mr{can}_{\mcF^\nabla_2}$ (cf.  Remark \ref{eRfi}).
Since $\nabla_j$ is the unique $\mcD_{\mbP^\mr{log}}^{(\N-1)}$-module structure on $\mcF_j$ extending $\nabla_{\mcF_j^\nabla}^\mr{can}$ via $\tau_j$,
the isomorphism $\alpha$, being an extension of $\alpha'$, preserves the $\mcD_{\mbP^\mr{log}}^{(\N-1)}$-action.
It follows that $\alpha$ defines an isomorphism  of level-$\N$ $\mr{GL}_2$-opers  $\msF^\heartsuit_1 \isom \msF^\heartsuit_2$, which induces  the equality  $\msF_1^\spadesuit = \msF_2^\spadesuit$.
This completes the proof of this proposition.
  \end{proof}
\SSP

\begin{rema} \label{esap77}
The proof of the above proposition shows that
${^\N}\mcO p^{^\mr{Zzz...}}_{2, \msP, (\rho_0, \rho_1, \rho_\infty)} = \emptyset$ unless
$(\rho_0, \rho_1, \rho_\infty) \in ((\mbZ/p^\N \mbZ)^\times / \{ \pm 1 \})^{\times 3}$.
As mentioned  in Remark \ref{peec}, this fact for $\N =1$ was already verified in ~\cite{Mzk2}.
\end{rema}
\SSP

We shall write
\begin{align}
\mr{Cov}^{\mr{tame}}_{\msP, +} \ \left(\text{resp.,} \ \mr{Cov}^{\mr{tame}}_\msP  \right)
\end{align}
for the set of equivalence classes of finite, separable,  and tamely ramified coverings $\phi : \mbP \migi \mbP$  satisfying the following conditions:
\begin{itemize}
\item
The set of ramification points of $\phi$ coincides with $\{ [0], [1], [\infty] \}$;
\item
If $\lambda_x$ ($x=0, 1, \infty$) denotes the ramification index of $\phi$ at $[x]$, then 
$\lambda_0$, $\lambda_1$, $\lambda_\infty$ satisfy the inequality $\lambda_0 + \lambda_1 + \lambda_\infty < 2 \cdot p^\N$ (resp., $\lambda_0$, $\lambda_1$, $\lambda_\infty$ are all odd and satisfy the inequality $\lambda_0 + \lambda_1 + \lambda_\infty < 2 \cdot p^\N$).
\end{itemize}
Here, the equivalence relation is defined in such a way that
two  coverings $\phi_1, \phi_2 : \mbP \migi \mbP$ are  equivalent if there exists an element $h \in  \mr{PGL}_2 (k) \left(= \mr{Aut}_k (\mbP)  \right)$ with $\phi_2 = h \circ \phi_1$.
Since the identity morphism $\mr{id}_\mbP$ of $\mbP$ defines a tamely ramified covering with ramification indices $(1,1,1)$,  the set $\mr{Cov}^{\mr{tame}}_{\msP, +}$ (resp., $\mr{Cov}^{\mr{tame}}_{\msP}$) is nonempty.
By applying some  of the  results proved so far, we  obtain the following assertion.

\SSP
\bt[cf. Theorem \ref{Re345}] \label{C00d3e}
The assignment $\phi \mapsto \msF^\spadesuit_\phi$ gives a $4$-to-$1$ (resp., a $1$-to-$1$, i.e., bijective) correspondence
\begin{align}
\Upsilon_{\msP, +} : \mr{Cov}^{\mr{tame}}_{\msP, +} \migisurj {^\N}\mcO p^{^\mr{Zzz...}}_{2, \msP} \ \left(\text{resp.,} \ \Upsilon_{\msP} :  \mr{Cov}^{\mr{tame}}_{\msP} \isom  {^\N}\mcO p^{^\mr{Zzz...}}_{2, \msP}\right).
\end{align}
In particular,  the set ${^\N}\mcO p^{^\mr{Zzz...}}_{2, \msP}$ is  finite and  admits  an inclusion   into the following set:
\begin{align} \label{Eran3}
\mbB_\N^\circledast :=\left\{ (\lambda_0, \lambda_1, \lambda_\infty) \in \mbB^{\times 3}  \, \big| \, 
\lambda_0 + \lambda_1 + \lambda_\infty < 2 \cdot p^\N  \ \text{and} \ |\lambda_0 - \lambda_1| \leq \lambda_\infty \leq \lambda_0 + \lambda_1 \right\},
\end{align}
where $\mbB$ denotes the set of positive odd integers $a$ with $p \nmid a$, $a < p^\N$.
\et
\begin{proof}
First, we shall consider the former assertion.
Let us take a dormant $\mr{PGL}_2$-oper $\msF^\spadesuit$ of level $\N$ classified by ${^\N}\mcO p^{^\mr{Zzz...}}_{2, \msP}$, and denote by    $(\rho_0, \rho_1, \rho_\infty)$  the radii of  $\msF^\spadesuit$.
The result of  Proposition \ref{P00245} shows that  $\msF^\spadesuit$ is the unique dormant $\mr{PGL}_2$-oper on $\msP$ of level $\N$ and radii $(\rho_0, \rho_1, \rho_\infty)$.
Now, let us choose a triple of integers $\vec{\lambda} := (\lambda_0, \lambda_1, \lambda_\infty)$ associated to $(\rho_0, \rho_1, \rho_\infty)$ as defined in \S\,\ref{SS042e}.
As observed in Remark \ref{R453}, this triple satisfies the inequalities in  (\ref{eapahks}).
It follows that the triples
\begin{align}
\vec{\lambda}_0 &:= (\lambda_0, p^\N - \lambda_1, p^\N -\lambda_\infty), \\
\vec{\lambda}_1 &:= (p^\N - \lambda_0, \lambda_1, p^\N -\lambda_\infty),\notag \\
\vec{\lambda}_\infty &:= (p^\N - \lambda_0, p^\N-\lambda_1, \lambda_\infty),  \notag
\end{align}
respectively, satisfy conditions  (a) and (b) in \S\,\ref{SS042e}, and conversely, each triple  of integers satisfying (a) and (b) is one of the four tripes $\vec{\lambda}$, $\vec{\lambda}_0$, $\vec{\lambda}_1$, $\vec{\lambda}_\infty$.
Thus, by Proposition \ref{P00f45}, the preimage of the element defined by $\msF^\spadesuit$ via  $\Upsilon_{\msP, +}$ coincides with  $\{[\phi_{\vec{\lambda}}], [\phi_{\vec{\lambda}_0}], [\phi_{\vec{\lambda}_1}], [\phi_{\vec{\lambda}_\infty}] \}$, where, for a triple $\vec{\lambda}' := (\lambda'_0, \lambda'_1, \lambda'_\infty)$,  
$[\phi_{\vec{\lambda}'}]$  denotes
a unique (up to equivalence) covering classified by 
  $\mr{Cov}_{\msP, +}^{\mr{tame}}$ whose ramification index at $[x]$ ($x = 0, 1, \infty$)  is   $\lambda'_x$.
  This proves the non-resp'd assertion.
  Also, the resp'd assertion follows from the fact that only one of the four  triples     
   $\vec{\lambda},  \vec{\lambda}_0, \vec{\lambda}_1, \vec{\lambda}_\infty$ satisfies the condition that all factors   are odd.
  
  Finally, the latter assertion follows from the resp'd portion of the former assertion (and its proof).
  This completes the proof of this theorem.
\end{proof}
\SSP

\begin{rema}
In the case of $\N=1$,
we know that the embedding ${^1}\mcO p^{^\mr{Zzz...}}_{2, \msP} \migiincl \mbB_1^\circledast$ resulting from Theorem \ref{C00d3e} is bijective (cf. ~\cite[Introduction, \S\,1.2,  Theorem 1.3, (2)]{Mzk2}).
The resulting correspondence ${^1}\mcO p^{^\mr{Zzz...}}_{2, \msP} \cong  \mbB_1^\circledast$ allows us to translate dormant $\mr{PGL}_2$-opers into edge-colorings  on trivalent  graphs, as well as  lattice points of a rational polytope  (cf. ~\cite{LO}, ~\cite{Wak2}).
However, at the time of writing the present paper, the author does not know much about the image of this map for a general $\N$.
\end{rema}

\SSP

\begin{rema} \label{R056}
In ~\cite[Chap.\,II, \S\,2, Definition 2.2]{Mzk2}, S. Mochizuki introduced the notion of a(n) {\it (dormant) $\M$-connection} (for each nonnegative integer $\M$) on a flat $\mbP^1$-bundle.
Here, we  recall its definition briefly. 
Let $\msX := (X, \{ \sigma_i \}_i)$ be   as in (\ref{ea34j}).
Also, let 
$(\mcP, \nabla)$ be a flat $\mbP$-bundle on $X^\mr{log}$, i.e.,  a $\mbP$-bundle $\mcP$ on $X$ equipped with a logarithmic connection $\nabla$ (with respect to the log structure of $X^\mr{log}$).
Denote by $W_{\M+1}$  the ring of Witt vectors with coefficients in $k$ of length $\M+1$.
Then,  a dormant $\M$-connection  on $(\mcP, \nabla)$ (of prescribed radii)  is defined as  a crystal in $\mbP$-bundles on the log crystalline site  $\mr{Crys} (X^\mr{log}/W_{\M+1})$  inducing  $(\mcP, \nabla)$ via reduction module $p$ and satisfying some other  conditions.
The condition of being a dormant $\M$-connection  is described in terms of $p^{\M+1}$-curvature in the sense of ~\cite[Chap.\,II, the discussion in \S\,2.1]{Mzk2}, which 
 is different from (but closely related to) our definition of $p^{\M+1}$-curvature.
According to ~\cite[Chap.\,IV, \S\,2, Theorem 2.3]{Mzk2},
there exists a bijective correspondence between the elements in $\mr{Cov}_{\msP}^{\mr{tame}}$ and the set of dormant $(\N-1)$-connections on 
a torally indigenous bundle (i.e., a $\mr{PGL}_2$-oper of level $1$)
 on $\msP$.
By combining this fact with Theorem \ref{C00d3e}, we see that 
 each  dormant $\mr{PGL}_2$-oper  on $\msP$ of level $\N$ and radii $\vec{\rho} \in ((\mbZ/p^\N\mbZ)^\times/\{ \pm1\})^{\times 3}$
  may be uniquely extended to a dormant $(\N-1)$-connection of the same radii.
\end{rema}

\LSP

\end{document}